%% file: JCV.tex
\documentclass[11pt]{article}
\usepackage[utf8]{inputenc}
\usepackage[T1]{fontenc}
\usepackage{amsthm,amsmath,amssymb,amsfonts,a4wide}
\usepackage[usenames]{color}
\definecolor{plum}  {rgb}{.6,0,.6}
\definecolor{forest}  {rgb}{0,.4,0.2}
\definecolor{midnight}  {rgb}{0,0,.7}
\usepackage[colorlinks=true,linkcolor=midnight,pdftex,citecolor=forest]{hyperref}
\usepackage{multicol,multirow,latexsym}
\usepackage{graphicx,subfigure}
\usepackage{mathrsfs,dsfont}
\usepackage{tikz}
\usepackage{natbib}
\usepackage{xcolor}
\usepackage{algpseudocode}
\usepackage{algorithm}
\usepackage{enumitem}
\usepackage{hyperref}
\usepackage{nameref,cleveref}
\usepackage[normalem]{ ulem }
\usepackage{soul,xr}
\numberwithin{equation}{section}
\theoremstyle{plain}
\newtheorem{theorem}{Theorem}%[section]
\newtheorem{definition}{Definition}
\newtheorem{remark}{Remark}
\newtheorem{proposition}{Proposition}
\newtheorem{lemma}{Lemma}

\newenvironment{assumption}[2]{\noindent\smallskip\\\textbf{Assumption} \textbf{#1} \emph{{#2}}}
\renewcommand{\L}{\mathbb{L}}

\newcommand{\E}{\ensuremath{\mathbb{E}}}
\renewcommand{\P}{\ensuremath{\mathbb{P}}}
\newcommand{\R}{\ensuremath{\mathbb{R}}}
\newcommand{\ZZ}{\ensuremath{\mathbb{Z}}}
\newcommand{\fX}{\mathrm{f}_X}
\newcommand{\hfX}{\tilde{\mathrm{f}}_X}
\newcommand{\nX}{n_X}

\newcommand{\var}{\mbox{Var}}

\newcommand{\argmax}{\mathop{\mathrm{arg\,\max}}}
\newcommand{\1}{\mathds{1}}

\newcommand{\frod}{\hat{f}_{\hat{h}}}
\newcommand{\CDRodeo}{\textsc{CDRodeo}}
\newcommand{\Rodeo}{\textsc{Rodeo}}
\newcommand{\Hhp}{\mathcal{H}_{\text{hp}}}
\newcommand{\HhpRev}{\mathcal{H}_{\text{hp}}^{\text{Rev}}}
\newcommand{\HhpDir}{\mathcal{H}_{\text{hp}}^{\text{Dir}}}
\newcommand{\Ehp}{\mathcal{E}_{\text{hp}}}
\newcommand{\Atilde}{\widetilde{\mathcal{A}}_n}
\newcommand{\rel}{\mathcal{R}}
\newcommand{\Bernf}[1]{\mathcal{B}\text{ern}_{\bar{f}}(#1)}
\newcommand{\Bernfabs}[1]{\mathcal{B}\text{ern}_{\vert\bar{f}\vert}(#1)}
\newcommand{\BernZ}[1]{\mathcal{B}\text{ern}_{\bar{Z}}(#1)}
\newcommand{\BernZabs}[1]{\mathcal{B}\text{ern}_{\vert\bar{Z}\vert}(#1)}
\newcommand{\f}[1]{\hat{f}_{#1}(w)}
\newcommand{\seuil}[1]{\lambda_{#1}}
\newcommand{\dX}{\tilde{\delta}_X}
\newcommand{\hX}{h_X}
\newcommand{\KX}{\mathcal{K}}
\newcommand{\fXK}{\tilde{\mathrm{f}}_{X}^{\mathcal{K}}}
\newcommand{\fXKi}[1]{\tilde{\mathrm{f}}_{X,#1}^K}
\newcommand{\Un}{~ \mathcal{U}_1}
\newcommand{\Up}{~ \mathcal{U}'_1}
\newcommand{\supp}{\mathrm{supp}}
\newcommand{\fbar}[1]{\bar{f}_{#1}(w)}
\newcommand{\Cs}{\text{C}_{\sigma}}
\newcommand{\CEZA}{\text{C}_{E\vert \bar{Z} \vert}}

\newcommand{\PP}[1]{\mathbb{P}\left(#1\right)}
\newcommand{\CgfA}{\text{C}_{\gamma \vert f\vert}}

\newcommand{\CvX}{\text{C}_{\text{var}_X}}
\newcommand{\CbiasX}{\text{C}_{\text{bias}_X}}
\newcommand{\CEbar}{\text{C}_{\bar{\text{E}}}}
\newcommand{\CEZ}{\text{C}_{E\bar{Z}}}

\newcommand{\condZbar}{\text{cond}_{\bar{Z}}}
\newcommand{\BZA}[1]{\mathcal{B}\text{ern}_{\vert\bar{Z}\vert}(#1)}
\newcommand{\UUn}{\ \mathcal{U}}
\newcommand{\BZ}[1]{\mathcal{B}\text{ern}_{\bar{Z}}(#1)}
\newcommand{\CgZA}{\text{C}_{\gamma \vert \bar{Z} \vert}}
\newcommand{\hirr}{h_{\text{irr}}}
\newcommand{\tirr}{t_{\text{irr}}}
\newcommand{\ibeta}{\beta^{-1}}
\newcommand{\Act}{\mathcal{A}\text{ct}}
\newcommand{\DZ}[1]{\Delta_{Z,#1}}

\newcommand{\MD}[1]{\text{M}_{\Delta}}
\newcommand{\CMD}{\text{C}_{\text{M}\Delta}}
\newcommand{\Cl}{\text{C}_{\lambda}}
\newcommand{\tminimax}{t(\Am,\CAm)}
\newcommand{\Am}{A}
\newcommand{\CAm}{C_\Am}
\newcommand{\biasbar}[1]{\bar{B}_{#1}}
\newcommand{\jm}{j_\Am}
\newcommand{\hint}[1]{h^{(\text{int},#1)}}
\newcommand{\hu}[1]{h^{[#1,u]}}
\newcommand{\huvar}[1]{h^{\lbrace #1,u\rbrace}}
\newcommand{\CBbar}{\text{C}_{\bar{B}}}
\newcommand{\logb}{\log_{\frac{1}{\beta}}}

\newcommand{\AfXmin}{$\boldsymbol{\mathcal{L}_X}$}
\newcommand{\AestfX}{$\boldsymbol{\mathcal{E}\fX}$}
\newcommand{\Monoton}{$\boldsymbol{\mathcal{C}}$}
\newcommand{\Monotonbis}{$\boldsymbol{\mathcal{M}}$}
\newcommand{\K}{\mathrm K}
\newcommand{\rcp}{\tfrac{7}{4}}  %ancien \red{C_p}
\newcommand{\X}{\widetilde{X}}

\Crefname{thm}{Theorem}{Theorems}
\crefname{thm}{theorem}{theorems}
\Crefname{assumption}{Assumption}{Assumptions}
\crefname{assumption}{assumption}{assumptions}
\Crefname{proposition}{Proposition}{Propositions}
\crefname{proposition}{proposition}{propositions}
\Crefname{lemma}{Lemma}{Lemmas}
\crefname{lemma}{lemma}{lemmas}
\Crefname{rmk}{Remark}{Remarks}
\crefname{rmk}{remark}{remarks}
\makeatletter
\newcommand{\myref}[1]{\cref{#1}\mynameref{#1}{\csname r@#1\endcsname}}
\newcommand{\Myref}[1]{\Cref{#1}\mynameref{#1}{\csname r@#1\endcsname}}

\graphicspath{{../}}

\date{\today}

\begin{document}

\title{Adaptive greedy algorithm for moderately large dimensions in kernel conditional density estimation}

\author{
{\sc Minh-Lien Jeanne Nguyen},\\[2pt]
 Mathematical Institute,\\ University of Leiden\\
Niels Bohrweg 1, 2333 CA Leiden, Netherlands\\
m.j.nguyen@math.leidenuniv.nl\\[7pt]
{\sc Claire Lacour},\\[2pt]
 LAMA, CNRS\\ Univ Gustave Eiffel,
        Univ Paris Est Creteil\\ 
        F-77447 Marne-la-Vall\'ee, France\\
        claire.lacour@univ-eiffel.fr \\[7pt]
{\sc Vincent Rivoirard} \\[2pt]
CEREMADE, CNRS, UMR 7534\\
       Universit\'e Paris-Dauphine, PSL University\\
       75016 Paris, France\\
       Vincent.Rivoirard@dauphine.fr\\[7pt]
}

\maketitle

\begin{abstract}
This paper studies the estimation of the conditional density $f(x,\cdot)$  of $Y_i$ given $X_i=x$, from the observation of an 
i.i.d.~sample $(X_i,Y_i)\in \R^d$, $i\in \{1,\dots,n\}.$
We assume that $f$ depends only on $r$ unknown components with  typically $r\ll d$.
%}
We provide an adaptive fully-nonparametric strategy based on kernel rules to estimate $f$. To select the bandwidth of our kernel rule, we propose a new fast iterative algorithm inspired by the Rodeo algorithm \citep{LW06} to detect the sparsity structure of $f$. More precisely, in the minimax setting, 
% if $f$ is $\beta$-H\"olderian (with $\beta$ also unknown), we prove   that 
our pointwise estimator, which is adaptive to both the regularity and the sparsity,  achieves the quasi-optimal rate of convergence. Our results also hold for density estimation. The computational complexity of our method is only $O(dn \log n)$. 
A deep numerical study shows nice performances of our approach.

\medskip

\noindent\textbf{Keywords:} \textit{Conditional density, Sparsity, Minimax rates, Kernel density estimators, Greedy algorithm.}
\end{abstract}

%Changements de notations $a$ devient $1+\epsilon$ ? $\UUn$ devient $\UUn_w$ ?
%%%%%%%%%%%%%%%%%%%%%%%%%%
%%%%%%%%%%%%%%%%%%%%%%%%%%
\section{Introduction}
%%%%%%%%%%%
\subsection{Motivations}
Consider $W=(W_1,\dots,W_n)$ a sample of a couple $(X,Y)$ of multivariate random vectors: 
for $i\in\{1,\dots,n\}$, $$W_i=(X_i,Y_i),$$
with $X_i$ valued in $\R^{d_1}$ and $Y_i$ in $\R^{d_2}$. We denote $d:=d_1+d_2$ the joint dimension.
 We assume that the marginal distribution of $X$ and the conditional distribution of $Y$ given $X$ are absolutely continuous with respect to the Lebesgue measure, and we denote by $\fX$ the marginal density of $X$  (and more generally by ${\rm f}_Z$ the density of any random vector $Z$). Let us define $f:\R^d\rightarrow\R_+$ such that for any $x\in\R^{d_1}$, $f(x,\cdot)$ is the conditional density of $Y$ conditionally on $X=x$:
$$f(x,y)dy=d\P_{Y|X=x}(y).$$ 
In this paper, we aim at estimating the conditional density $f$ at a set point $w=(x,y)$   in $\R^d$.

{Estimating a conditional density may be done in any regression framework,  i.e.  as soon as we observe a (possibly multidimensional)  response $ Y $ associated with a (possibly multidimensional) covariate $X$. 
The regression function $ \E [Y | X = x] $ is often studied, but this mean is in fact a summary of the entire distribution and may lose information (think in particular to the case of an asymmetric or multimodal distribution).}
Thus the problem of estimating the conditional distribution is considered in various application fields: 
meteorology, insurance, medical studies, geology, astronomy. See \cite{Jeanne1} and references therein. Moreover, the
ABC methods (Approximate Bayesian Computation)  are actually dedicated to find a conditional distribution (of the parameter given observations) in the case where the likelihood is not computable but simulable: see 
 \cite{IzbickiLeePospisil18} (and references therein) where the link between conditional density estimation and ABC  is studied.

Several nonparametric methods have been proposed for estimating a conditional density:  \cite{HBG96} and \cite{FYT96} have improved the seminal Nadaraya-Watson-type estimator of \cite{rosenblatt69} and 
\cite{LinchengZhijun85}, as well as \cite{DgZ03} who introduced another weighted kernel estimator. For these kernel estimators, different methods have been advocated to tackle the bandwidth selection issue: bootstrap approach \citep{BH01} or  cross-validation  variants, see \cite{FYim04, HGI10}, %\cite{efromovich10b}, 
\cite{IchimuraFukuda10}.
Later, adaptive-in-smoothness estimators have been introduced:  \cite{brunelcomtelacour07} with piecewise polynomial representation, \cite{Chagny} with wraped base method, \cite{CohenLePennec} with penalized maximum likelihood estimator, \cite{BLR16} with Lepski-type methods and \cite{sart17} with tests-based histograms.
%Sugiyama et al 2010

All above references do not really deal with the curse of dimensionality. From a theoretical point of view, the minimax rate of convergence for such nonparametric statistical problems is known to be $n^{-s/(2s+d)}$ (possibly up to a logarithmic term), where $s$ is the smoothness of the target function. This illustrates that estimation gets increasingly hard when $d$ is large. Moreover the computational complexity of above methods is often intractable as soon as $d$ is larger than 3 or 4. 
A first answer to overcome this limitation is to consider single-index models as \cite{BouazizLopez10} or semi-parametric models as \cite{FPYZ09}, but this implies a strong structural assumption.
A more general advance has been made by \cite{HRL04} who assume that some components of $X$ can be irrelevant, i.e. that they do not contain any information about $Y$ and should be dropped before conducting inference. Their cross-validation approach allows them  to obtain a minimax rate for a $r_1$-dimensional 
$C^2$ function, where $r_1$ is the number of relevant $X$-components. \cite{efromovich10b} has improved these non-adaptive results by using thresholding and Fourier series and achieves the minimax rate  $n^{-s/(2s+r_1)}$ without any knowledge of $r_1$ nor $s$. Note that above rates were established for the $\L^2$-loss whereas we shall consider the pointwise loss.
Moreover these combinatorial approaches make their computation cost prohibitive when both $n$ and $d$ are large.
In the same framework, \cite{Shiga2015} assume that the dependence of $Y$ on the relevant components is additive. Another way is paved by \cite{OtneimTjostheim18} who estimate the dependence structure in a Gaussian parametric way while estimating marginal distributions nonparametrically.
More recently, \cite{IzbickiLee16, IzbickiLee17} have proposed two attractive methodologies using orthogonal series estimators 
in the context of  an eventual smaller unknown intrinsic dimension of the support of the conditional density. In particular, the Flexcode method originally proposes to transfer successful procedures for high dimensional regression  to the conditional density estimation  setting by interpreting the coefficients of the orthogonal series estimator as regression functions, which allows to adapt to data with different features (mixed data, smaller intrinsic dimension, relevant variables) in function of the regression method. 
However, the optimal tuning parameters depend in fact on the unknown intrinsic dimension. Furthermore, optimal minimax rates are not achieved, revealing the specific nature of the problem of conditional density estimation, more intricate, in full generality, than regression.  %\textcolor{green}{Tartiner un peu plus le papier o\`u il y a Flexcode?}
%%%%%%%%%%%
\subsection{Objectives, methodology and contributions}\label{sec:OMC}
In this paper, we wish to estimate the conditional density $f$ by assuming that only $r\in \{0,\ldots,d\}$ components are {\it relevant}, i.e.  that there exists a subset $\mathcal{R}\subset\lbrace1,\dots,d\rbrace$  with cardinal $r$, such that for any fixed $\lbrace z_j\rbrace_{j\in\mathcal{R}}$,  the function $\lbrace z_k\rbrace_{k\in\mathcal{R}^c}\mapsto f(z_1,\dots,z_d)$ is constant on the neighborhood of $w$, with $\mathcal{R}^c=\{1,\ldots,d\}\setminus\mathcal{R}$. We denote $f_{\mathcal{R}}$ the restriction of $f$ to the relevant directions.  Assuming that $f$ is $s$-H\"olderian, our goal is to provide an  estimation procedure such that it achieves the best adaptive rate. The meaning of {\it adaptation} is \emph{twofold} in this paper; the first meaning corresponds to adaptation with respect to the smoothness, which is the classical meaning of adaptation. The second one corresponds to adaptation with respect to the sparsity. So, our goal is to propose an optimal procedure in this context, meaning that it does not depend on the knowledge of $s$ and $r$, and even $\mathcal{R}$. Furthermore, for practical purposes in moderate large dimensions, it should be implemented with low computational time.

For this purpose, we consider a particular kernel estimator depending on a bandwidth $h\in\R_+^d$ to be selected. To circumvent the curse of dimensionality, we consider an iterative algorithm on a special path of bandwidths inspired by the  {\it \textsc{Rodeo}} procedures proposed by \cite{LW06} and \cite{LW08} for nonparametric regression, \cite{LLW07} for density estimation and \cite{Jeanne1} for conditional density estimation. More precisely, our new procedure, called RevDir \CDRodeo{}, is a variation of the \CDRodeo{} proposed by \cite{Jeanne1} (and called Direct \CDRodeo{} in the sequel). Each iteration step of this new algorithm is based on comparisons between partial derivatives of our kernel rule, denoted $Z_{hj}$, and specific thresholds $\lambda_{hj}$, respectively defined in \eqref{defZ} and \eqref{threshold}. Let us mention that for variable selection in the regression model with
very high ambient dimension, \cite{CD12} used similar ideas to select the relevant variables by comparing some quadratic functionals
of empirical Fourier coefficients to prescribed significance levels. Consistency
of this (non-greedy) procedure is established by \cite{CD12}.

We establish that, up to a logarithmic term whose exponent is positive but as close to 0 as desired, RevDir \CDRodeo{} achieves the rate $((\log n)/n)^{s/(2s+r)}$, which is the optimal adaptive minimax rate on H\"older balls $\mathcal{H}_d(s, L)$, when the conditional density depends on $r$ components.  When $r$ is much smaller than $d$, this rate is much faster than the usual rate $((\log n)/n)^{s/(2s+d)}$ achieved by classical kernel rules. 
Furthermore, unlike previous \textsc{Rodeo}-type procedures, our procedure is adaptive with respect to both the smoothness and the sparsity.
To the best of our knowledge, our RevDir \CDRodeo{}  procedure is the first algorithm achieving quasi-minimax rates for conditional density estimation in this setting where both sparsity and smoothness are unknown. 
We lead a deep numerical study of parameters tuning of the algorithm. Then the numerical performances are presented for several examples of conditional densities. 
In particular RevDir \CDRodeo{} is able to tackle the issue of sparsity detection. Moreover, for each relevant component, reconstructions are satisfying.
%
%Furthermore, tuning RevDir \CDRodeo{}  is easy (see Section \ref{sssectAssumpResult}) \textcolor{red}{Phrase sur les r\'esultats de simulation} and
Finally, we show that the total worst-case complexity of the RevDir \CDRodeo{}  algorithm is only $O(d n\log n).$ 
This last result is very important for modern statistics where many problems deal with very large datasets.
%%%%%%%%%%%%%%%%%
\subsection{Plan of the paper and notation}
The plan of the paper is the following. First we describe in Section \ref{sec:method} the estimation procedure. We give heuristic ideas based on the minimax approach and explain why some modifications of the Direct \CDRodeo{} procedure are necessary. Then a detailed presentation of our algorithm is provided in Section~\ref{sec:RevDir}. Next, the main result is stated in Section \ref{sec:resul}. The complexity of the algorithm is computed in Section~\ref{sec:complexity}. After tuning the method, the latter  is illustrated via simulations and examples in Section \ref{sec:simus}. The proofs are gathered in Section \ref{sec:proofs}.

\bigskip 
In the sequel, we adopt the following notation. Given two functions $\phi,\psi:\R^d\to\R$, two integers $j,k$, two vectors $h$ and $h'$, two real numbers $a$ and $b$, we denote
\begin{itemize} 
\setlength\itemsep{0em}
\item[-] $\|\phi\|_q=\left(\int |\phi(u)|^qdu\right)^{1/q}$ the $\L_q$ norm of $\phi$ for any $q\geq 1$;
\item[-] $\phi\star\psi$ the convolution product $u\mapsto\int_{\R^d} \phi(u-v)\psi(v)dv$; 
\item[-]  $\partial_j \phi$ the partial derivative of $\phi$ with respect to the direction $j$ (or $\frac{\partial}{\partial u_j} \phi$ when there is ambiguity on the variable);
\item[-] $j:k$ the set of integers from $j$ to $k$;
\item[-] $|A|$ the cardinal of a set $A$,
\item[-] $h\preceq h'$ the partial order on vectors defined by: $h_k\leq h'_k$, for $k\in 1:d$.
\item[-] $a \lesssim b$ (respectively $a\approx b$) means that  the inequality (respectively the equality) is satisfied up to a constant.
\end{itemize}

%%%%%%%%%%%%%%%%%%%%%%%%%%
%%%%%%%%%%%%%%%%%%%%%%%%%%
\section{Estimation procedure}\label{sec:method}
%%%%%%%%%%%%%%%%%%%%%%%%%%
%%%%%%%%%%%%%%%%%%%%%%%%%%
As mentioned in Introduction, the goal of this paper is to provide an estimator of the conditional density achieving pointwise adaptive minimax rates, where the meaning of adaptation is twofold as explained in Section~\ref{sec:OMC}. %Hereafter, we fix the estimation point $w=(x,y)$.

Our estimation procedure follows the kernel methodology. We use a specific family of kernel estimators \citep{BLR16}, called hereafter the \textsc{Blr} estimators and detailed in Section~\ref{kernel:sec}. The selection of the bandwidth is introduced with heuristic considerations and detailed in Section~\ref{sel:section} in the spirit of \Rodeo{} \citep{LW08,Jeanne1}. After presenting advantages and limitations of the latter, we propose a new algorithm called RevDir \CDRodeo{}. 

\color{black}

%%%%%%%%%%%%%%%%%%%%%%
\subsection{Kernel rule}\label{kernel:sec}
We use the \textsc{Blr} family of kernel estimators as it presents some significant advantages explained below. The \textsc{Blr} family is defined as follows. Let $K:\R\rightarrow\R$ be a kernel function, namely $K$ satisfies $\int_{\R} K(t)dt =1$. Then, given a bandwidth $h=(h_1, \dots,h_d)\in(0,1]^d$, the estimator of $f(w)$ associated with $K$ and $h$ is defined by 
\begin{equation}\label{Deffh}
\f{h}:=\frac{1}{n} \sum\limits_{i=1}^n \frac{1}{\hfX\left(X_i\right)}\K_h(w-W_i), 
\end{equation}
where for any $v\in\R^d$,
$$\K_h(v)=\prod_{j=1}^{d}h_j^{-1}K(v_j/{h_j})$$ and
$\hfX$ is an estimator of $\fX$, built from a sample $\tilde X$ not necessarily independent of $W$. 
\begin{remark}\label{rmkKDE}
Note that (non conditional) density estimation is a special case of this problem studied, as it corresponds to the setting where $d_1=0$ and $\fX\equiv1 \ (\equiv\hfX)$. In this case, $\f{h}$ is the usual kernel density estimator.
\end{remark}
Several arguments justify the choice of the \textsc{Blr} family, rather than the intensively studied family expressed as a ratio of two density estimates of ${\rm f}_W$ and $\fX$, following:
$$f(x,y)=\frac{{\rm f}_W(x,y)}{\fX(x)}.$$  
Indeed, this last decomposition takes into account the characteristics (smoothness, sparsity) of ${\rm f}_W$ and $\fX$ instead of those of our target $f$. More precisely, an irrelevant component of the conditional density may be relevant for both the joint density ${\rm f}_W$ and the marginal density $\fX$ and it occurs in particular when a component of $X$ is independent of $Y$. Similarly, the smoothness of $f$ can be different from those of the functions ${\rm f}_W$ and $\fX$, which potentially would deteriorate the rates of convergence.

Conversely, the \textsc{Blr} estimators estimate $f$ more directly: in particular, their expectations can be written as the usual kernel regularization of $f$: under some mild assumptions on $K$ and $f$ and with $\hfX=\fX$,
{\small
\begin{equation}\label{egaliteKh}
\E[\f{h}]=\iint \frac{1}{\fX(u)}\K_h(w-(u,v)){\rm f}_W(u,v)dudv= \int \K_h(w-z)f(z)dz=(\K_h\star f)(w).
\end{equation}
}
 % $\fX$ is unknown, see \Cref{sssectAssumpResult} for the theoretical requirement about the choice of the estimator $\hfX$.
 
 %Equality~\eqref{egaliteKh} shows that the selection of $h$ will be essentially dictated by the intrinsic properties of the conditional density~$f$.

%%%%%%%%%%%%%%%%%%%%%%%%%%%%%%%%
\subsection{Selection of the bandwidth}\label{sel:section}
The principal issue in kernel rules is the choice of the bandwidth.  In particular, we consider a $d$-dimensional bandwidth,  instead of a scalar one which would be easier and faster to select but would also deteriorate the performances of the estimator.
%%%%%%%%%%%%
\subsubsection{Heuristic minimax arguments}\label{sec:heuristic}
We consider $\mathcal{H}_{d\vert r}(s, L)$ the set of functions of $\mathcal{H}_{d}(s, L)$ with at most $r$ relevant components, and its associated (squared) pointwise minimax risk
$$\underset{{\hat{T}_n}}{\inf}\underset{f\in\mathcal{H}_{d\vert r}(s, L)}{\sup} \E[(\hat{T}_n(w)-f(w))^2],$$ where the infimum is taken over all estimators of $f$ built from the sample $W$.

In the case of kernel rules, let us denote $h^*$ the \emph{minimax bandwidth} minimizing this risk.
We can decompose the squared risk in bias and variance terms:
\begin{equation}\label{biais-variance}
R(h^*):=\E[(\f{h^*}-f(w))^2]=B^2(h^*)+\var(\f{h^*}).
\end{equation}
For any bandwidth $h\in (0,1]^d$, the usual respective upper bounds  for the bias and variance are typically
\begin{equation}\label{eEquBiasfh}
B^2(h):=\Big(\E[\f{h}]-f(w)\Big)^2\lesssim \sum\limits_{j\in\rel} h_j^{2s}
\end{equation}
and
\begin{equation}\label{eEquVarfh}
\mathrm{V}_h:= \var(\f{h})\lesssim \frac{1}{n \prod_{j=1}^dh_j}.
\end{equation}
The minimizer $h^*$ on $(0,1]^d$ of the minimax risk is then of the form:
\begin{equation}
h^*_j=\left\lbrace
\begin{array}{l}
n^{-1/(2s+r)} \text{ for } j\in{\mathcal R},\\ 
1 \text{ for }j\notin{\mathcal R}.
\end{array} 
\right.
\end{equation}
Given this bandwidth, which depends on $s$, $r$ and ${\mathcal R}$, and given a sharp estimator $\hfX$, the \textsc{Blr} estimator achieves the minimax rates $n^{-\frac{s}{2s+r}}$. In the literature of conditional density estimation, to the best of our knowledge, no method provides theoretical results achieving the \emph{twofold} adaptive rates. Besides, the smoothness-adaptive procedures of bandwidth selection are based on optimization over $d$-dimensional grids of bandwidths, thus require intensive computation, even in moderately high dimension as the grid grows exponentially fast with the dimension. 

The principle of \Rodeo{}, and of its derived versions \citep{LW06, Jeanne1}, is to progressively build a monotonous path of bandwidths through the bandwidths grid. The construction of this path is based on tests at each iteration to decide if a bandwidth component has a convenient level or still has to be multiplied by an iterative step factor.  %to ensure a reasonable running time. 
The tests rely on the partial derivatives of the estimator with respect to the components of the current bandwidth: for $h\in(0,1]^d$ and $j\in 1:d$,
 \begin{equation}\label{defZ}
Z_{hj}:= \frac{\partial}{\partial h_j} \hat{f}_h(w).
\end{equation}
The main idea is to use $Z_{hj}$ as a proxy of $\frac{\partial}{\partial w_j} f$, relying on the natural intuition that the more $f$ is varying, the smaller the bandwidth is needed to fit the curve.  It is consistent with the minimax bandwidth level $h^*_j=1$ for irrelevant $j$ and the flatness of the curve in such a direction. 
Using the \textsc{Blr} family of conditional density estimators, the $Z_{hj}$'s are well defined as soon as the kernel $K$ is  $C^1$.
%, choice made throughout the paper. 
They are straightforwardly expressed, thus easily implementable, by using the following equation:
\begin{equation}\label{formulaZhj}
Z_{hj}=-\frac{1}{n h_j^2}
 \sum\limits_{i=1}^n \frac{1}{\hfX(X_i)} J\left(\tfrac{w_j-W_{ij}}{h_j}\right)
\prod\limits_{k\neq j}^{d}h_k^{-1}K\left(\frac{{w}_k-W_{ik}}{h_k}\right),
\end{equation}
where $J$ denotes the function $t\mapsto K(t) + tK'(t)$. Note that, under the condition  $\hfX=\fX$, if $j$ is an irrelevant component,
\begin{equation} 
\E[ Z_{hj}]=0,
\end{equation}
 which is expected in view of \eqref{defZ}  (see Lemma~\ref{lmZbar} in Appendix or Lemma~6 of \cite{Jeanne1} for a rigorous proof). The tests involved in the \Rodeo{} procedure consist in comparing $|Z_{hj}|$ to a threshold $\seuil{hj}$.
The threshold is chosen as follows:~
\begin{equation}\label{threshold}
\seuil{hj}:=\Cl\sqrt{\frac{(\log n)^{a}}{n h_j^2 \prod_{k=1}^d h_k}},
\end{equation}
with $\Cl= 4\Vert J\Vert_2\Vert K\Vert_2^{d-1}$ and an hyperparameter $a>1$. It is determined by Bernstein's concentration inequalities to ensure that with high probability $Z_{hj}$ is close to its expectation: $|Z_{hj}-\E[Z_{hj}]|\leq \frac12\seuil{hj}$.
%$\E[Z_{hj}]$ up to a distance $\frac12\seuil{hj}$.  
The hyperparameter $a$ quantifies the degree of high probability.
This definition is justified by following heuristic arguments. With $B(h)=\E[\f{h}]-f(w)$,
$$\frac{\partial}{\partial h_j} B(h)=\frac{\partial}{\partial h_j}\E[\f{h}]=\E\left[\frac{\partial}{\partial h_j}\f{h}\right]=\E[Z_{hj}].$$
If the upper bound of \eqref{eEquBiasfh} is tight and since, with large probability, $Z_{hj}\approx \E[Z_{hj}]$, we obtain, for $j\in \mathcal{R}$
$$|Z_{hj}|\approx h_j^{s-1}.$$
We stop the algorithm when $|Z_{hj}|\approx\seuil{hj}$ since for this bandwidth $h$, we expect
$${h_j}^{s-1}\approx \seuil{hj}\approx\frac{1}{h_j\sqrt{n\prod_{k=1}^dh_k}} \text{ (up to the logarithmic term)},$$
which corresponds to the minimax bandwidth $h^*$ which satisfies the minimax trade-off:
$${h^*_j}^{2s}\approx  \frac{1}{n \prod_{j=1}^dh^*_j},$$
for $j\in \mathcal{R}$. 
%%%%%%%%%%%%
\subsubsection{Initialization of the algorithm and variants of    \CDRodeo{}}\label{Init}
The previous paragraph explains quantities involved in the algorithm, its main ideas and the stopping criterion. We now study the initialization of the algorithm. We describe several alternatives.
\paragraph{Direct \CDRodeo{}  algorithm.}
The natural idea consists in  initializing the bandwidth at a large enough level and then decreasing the components of the bandwidth until 
$\vert Z_{hj}\vert\leq \seuil{hj}$. The detailed procedure is stated in \Cref{AlgoDirect}. \\
\begin{algorithm*}[h!]
\caption{Direct \CDRodeo{}  algorithm \setcounter{enumi}{-1}  \label{AlgoDirect}}
\begin{itemize}
\item[$\mbox{ }$ \hspace{1cm}\bf Given] a starting bandwidth $h^{(0)}=(h_0,\ldots,h_0)$ with $h_0>0$, the decreasing iterative step factor $\beta\in(0,1)$, a hyperparameter $a>1$, the activation of all components.
\item[\bf While] {
there are still active components,\\
 for all active component $j$, we test if $|Z_{hj}|$ is large (with respect to a threshold $\lambda_{hj}$ defined in \eqref{threshold}):
}
\begin{itemize}
\item[-] If $|Z_{hj}|>\seuil{hj}$, then $h_j\leftarrow\beta h_j$, and $j$ remains active.
\item[-] Else, $j$ is deactivated and $h_j$ remains unchanged for the next steps of the path.
\end{itemize}
\item[\bf Output] The loop stops when either all components are deactivated or the bandwidth is too small $\left(\prod_{j=1}^d h_j<\frac{\log n}{n}\right)$, then the final bandwidth is selected and denoted $\hat{h}$.
\end{itemize}
\end{algorithm*}

This procedure, called Direct \CDRodeo{}, has been deeply studied by \cite{Jeanne1}. % following \cite{??}.
Two cases can be distinguished for a component $h_j$. Either $h_j$ is selected at the first iteration, or when $\vert Z_{hj}\vert\approx \seuil{hj}$.

%We denote the initial large bandwidth $h^{(0)}=(h_0,\ldots,h_0)$.
 In the first case, remark that testing $\vert Z_{h^{(0)}j}\vert\leq \seuil{h^{(0)}j}$ corresponds to testing the hypothesis $\vert \E[Z_{h^{(0)}j}]\vert \leq \frac12 \seuil{h^{(0)}j}$, which is satisfied for any irrelevant component $j$: for any $h$,
$\E[ Z_{hj}]=0.$
%
% \begin{proof} For $j\in \rel^c$, assuming $\hfX=\fX$ (which {\color{blue}is relaxed in the proof of the main results} if $\hfX$ estimates $\fX$ well enough),  \red{hyp. à préciser ?}
%$$
%\E[ Z_{hj}]=\frac{-1}{h_j^2}\int_{\R^d}J\left(\frac{w_j-u_j}{h_j}\right)
%\left[\prod\limits_{k\neq j}^{d}h_k^{-1}K\left(\frac{{w}_k-u_k}{h_k}\right)\right]\frac{{\rm f}_W(u)}{{\rm f}_X(u_{1:d_1})}du.$$
%Remark that for $k\in\mathcal{R}^c$ and $k\not=j$, $\int_{\R} h_k^{-1}K\left(\frac{{w}_k-u_k}{h_k}\right)du_k=1$. Thus,
%$$
%\E[ Z_{hj}]=\frac{-1}{h_j^2}\int_{u_\mathcal{R}\in\R^r}\left[\prod\limits_{k\in\mathcal{R}}h_k^{-1}K\left(\frac{{w}_k-u_k}{h_k}\right)\right]f_\mathcal{R}(u_\mathcal{R})du_\mathcal{R}\underbrace{ \int_{u_j\in\R} J\left(\frac{w_j-u_j}{h_j}\right)d{u_j}}_{=0}=0.
%$$
%\end{proof}
So, with high probability the irrelevant bandwidth components are selected at the initialization level $h_0$, \textit{i.e.} as large as allowed by the procedure, in line with the minimax approach.

\begin{figure}[h!]
\begin{tikzpicture}[scale=2.2]
%%% zones %%%%
\fill[color=gray!10]
(4,4) -- (4, 2.56)%plot [domain=4:5] (\x, {4*(0.64)^(\x-4))})-- (5,2.56)
-- (9,2.56)
-- (9,4) -- (4,4) ;

\fill[color=gray!60]
(4,2.56) -- (4,1.873)%plot [domain=5:5.7] (\x, {4*(0.64)^(\x-4))})-- (5.7,1.873) 
-- (9,1.873)
-- (9,2.56) --(5,2.56) ;

\fill[color=gray!10]
(4,1.873) -- (4,0.87) %plot [domain=5.7:7] (\x, {4*(0.64)^(\x-4))})-- (7,1.05) 
-- (9,0.87)
-- (9,1.873) -- (5.7,1.873);

\fill[color=gray!60]
 (4,0.87) -- (4,0) --
 (9,0) -- (9,0.87) -- cycle;
 
%% nom des zones %%
  \draw (5,4.6) node[draw,fill=gray!10] {$\quad$};
   \draw (6.1,4.6) node { zone where $|Z_{hj}|> \lambda_{hj}$};
  
% \draw[gray] (3,1.5) -- (4,2.21);
  % \draw[gray] (3,1.5) -- (4,0.44);
   \draw (5,4.4) node[draw,fill=gray!60] {$\quad$};
  \draw (6.1,4.4) node { zone where $|Z_{hj}|\leq \lambda_{hj}$};

%% cadre general %%
\draw[->, thick] (4,0) -- (9.5,0);
\draw (4,4.5) node[above] {$h_j$};
\draw [->, thick] (4,-0.1) -- (4,4.5);
\draw(8.3,0) node[below] {Direct CDRodeo  progression};
\draw (4-0.2,0) node {$0$};
\draw [double] (9,0) -- (9,4.5);
\draw (9,4.55) node[above right] {output};

%% premier chemin %%
\draw (4,4) node[ left, blue] {$h_0$};
\draw (4,4) node[ blue] {$\bullet$};
\draw (4,3.66) node[ left, blue] {$\beta h_0$};
%\draw (4.2,3.66) node[ blue] {$\bullet$};
\draw[blue, dashed] (4,3.66) -- (4.2,3.66);
\draw (4-0.2,3.46) node[blue] {$\cdot$};
\draw (4-0.2,3.41) node[blue] {$\cdot$};
\draw (4-0.2,3.36) node[blue] {$\cdot$};
\draw[blue, thick] [domain=4:5] plot (\x, {4*(0.64)^(\x-4))});
\draw[->, blue, thick]  (5, 2.56) -- (9,2.56);
\draw(9,2.56)  node[right, blue] {$ \hat h_j\gg  h_j^*$ };

%% deuxieme chemin%%%
\draw (4,1.7) node[ above left,cyan] {$h_0'$};
\draw[cyan, dashed] (4, 1.8) -- (5.78,1.8);
\draw (5.78,1.8) node[ cyan] {$\bullet$};
\draw (4,1.55) node[ left, cyan] {$\beta h_0'$};
\draw[cyan, dashed] (4, 1.55) -- (6.12,1.55);
\draw (4-0.2,1.35) node[cyan] {$\cdot$};
\draw (4-0.2,1.3) node[cyan] {$\cdot$};
\draw (4-0.2,1.25) node[cyan] {$\cdot$};
%\draw (6.12,1.55) node[ cyan] {$\bullet$};
\draw[black, dashed] [domain=5:5.76] plot (\x, {4*(0.64)^(\x-4))});
\draw[cyan, thick] [domain=5.78:7.42] plot (\x, {4*(0.64)^(\x-4))});
\draw[->, cyan, thick]  (7.42, 0.87) -- (9,0.87);
\draw(9,0.87)  node[right, cyan] {$\hat h_j'\approx h_j^*$ };

\end{tikzpicture}
\caption{Two bandwidth paths for Direct CDRodeo with two different initializations, when $h_j \mapsto |Z_{hj}|/\lambda_{hj}$ is not monotonuous (larger than 1  in the  lightgray zone and smaller than 1 in the darkgray zone). Starting with a large $h_0$, the algorithm stops when $|Z_{hj}|$  becomes smaller than $\lambda_{hj}$ and provides a too large output bandwidth.
Starting with $h_0'$, the algorithm can provide the optimal bandwidth $h_j^*$.
Observe that the area where $\E[Z_{hj}]\geq \seuil{hj}$  is unknown, so $h_0'$ is intractable.}\label{Tikz}
\end{figure}
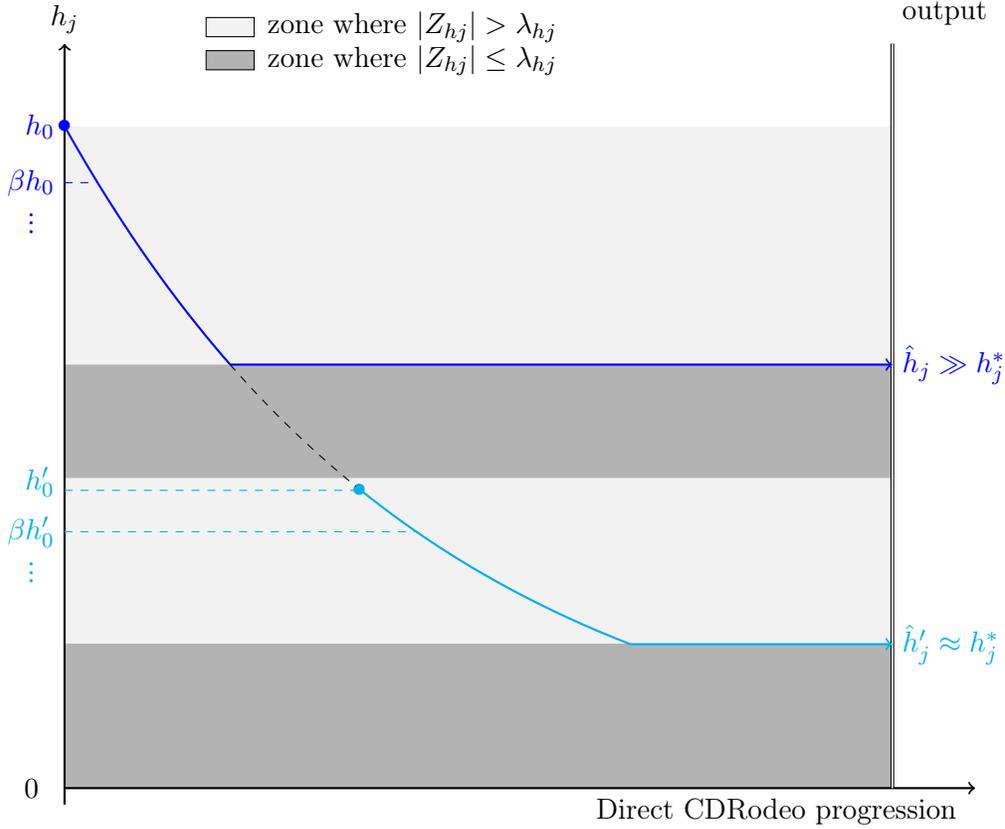

In the second case, the component $j$ is selected after a few iterations, and $\vert Z_{hj}\vert\approx\seuil{hj}$ (where the approximation is due to the discretization in $\lbrace \beta^k h_0, k\in\mathbb{N}\rbrace^d$). Thus with high probability: $\frac12 \seuil{hj} \lesssim\vert\E[Z_{hj}]\vert \lesssim \frac34\seuil{hj}$. 
For a relevant component $j$, for $s$ an integer larger than~1, \cite{Jeanne1} proved that
$$|\E[Z_{hj}]|\approx h_j^{s-1},$$
if the derivative satisfies $\vert \partial_j^{s} f \vert >0$ on the neighborhood of the evaluation point $w$.
%
%Thus the bandwidth selected by Direct \CDRodeo{} satisfied:
%\begin{equation}\label{fenetreminimax}
%{h_j}^{s-1}\approx \seuil{h^*j}\approx\frac{1}{h_j^*\sqrt{n\prod_{k=1}^dh_k^*}} \text{ (up to the logarithmic term)},
%\end{equation}
%\textcolor{red}{probl\`eme de notations} which is the required trade-off for the minimax rate over the class of $C^{s}$-functions. In particular, the minimax bandwidth $h^*$ satisfies for $j\in\rel$,
%\begin{equation}\label{fenetreminimax}
%{h_j^*}^{s-1}\approx \seuil{h^*j}\approx\frac{1}{h_j^*\sqrt{n\prod_{k=1}^dh_k^*}} \text{ (up to the logarithmic term)}.
%\end{equation}.
The assumption is quite restrictive. In particular, it excludes any density that is locally a polynomial of order smaller than $s$. Moreover, $s$ has to be an integer.\\
%The main purpose of this assumption is to ensure the lower bound $\vert \E[Z_{hj}]\vert \gtrsim h_j^{s-1}$ for any $h_j$ smaller than the selected bandwidth (the upper bound being easy to prove with standard assumptions for kernel estimators). But the real requirement is to have an increasing derivative of the bias $\vert \E[Z_{hj}] \vert$ \emph{i.e.} bias convexity/concavity.
%Then, for any bandwidth value $h$ potentially selected by \CDRodeo{}:
%\begin{align*}\label{eqBiasIntegralZhj}
%\left(B(h)\right)^2
%&= \left(-\sum\limits_{j\in\rel} \int_{u=0}^1 h_j \E[Z_{(uh)j}] du \right)^2
%\approx \left(\sum\limits_{j\in\rel}  h_j \E[Z_{hj}]\right)^2
%\approx\left( \sum\limits_{j\in\rel}  h_j \seuil{hj} \right)^2
%\approx \frac{1}{n\prod_{k=1}^dh_k} \text{ (up to the logarithmic term)}\text{\red{utiliser la notation approx log}} \\
%&\approx V(h),
%\end{align*}
%which leads to the bias-variance trade-off ensuring  minimax rates.\\
When this assumption is not satisfied, Direct \CDRodeo{} may stop with a too large bandwidth. Indeed, remember it begins with a large initial bandwidth in order to select large irrelevant bandwidth components, but the relevant components have to be selected much smaller. Between these two bandwidth levels, 
$\E[Z_{hj}]$ may have a change of sign, thus vanishes briefly before becoming larger (in absolute value) than $\seuil{hj}$ again. We have illustrated this problem in Figure~\ref{Tikz} where we show that the initialization $h_0$ is not convenient. 
%results in a too early stopping whereas it would have been preferable to initialize as $h_0'$ for this relevant component but it also reduces the irrelevant components to $h_0'$ instead of $1$, leading also to a larger variance and thus, a non-optimal trade-off.
 %Note also the area where $\E[Z_{hj}]\geq \seuil{hj}$  is unknown, so $h_0'$ is intractable.

In view of this issue,  we consider in the following some variations to the {\it Direct} \CDRodeo{} procedure. 

\color{black}

\paragraph{A Reverse \CDRodeo{}  algorithm.}
The first variation which could be considered is the {\it Reverse} \CDRodeo{} procedure in the same spirit as \cite{LLW07} (see Section~4.2 therein). We start with a  small bandwidth and use a sequence of non-decreasing  bandwidths to select the optimal value, still by comparing the $Z_{hj}$'s with the $\seuil{hj}$'s.  More precisely, instead of decreasing the bandwidth components by multiplied them by the factor $\beta$ when $\vert Z_{hj}\vert \geq \seuil{hj}$, the reverse algorithm increases them by \emph{dividing} them by $\beta$ when $\vert Z_{hj}\vert < \seuil{hj}$. 
Note that with this second test, it does not matter if $Z_{hj}$ vanishes. As illustrated by \cite{LLW07}, this approach is very useful for image data. However, the choice of the initial bandwidth is very sensitive. In particular, assume that $f$ has a very low regularity and has only one relevant component, say the first one for instance. In this case, if $h^*$ is the ideal bandwidth, $h^*_1$ has to be as small as possible, i.e. $h^*_1=1/n$ (up to a logarithmic term). 
Therefore, since $\mathcal R$ is unknown, the initialization of the  bandwidth must be not larger than  $h_{0,\text{rev}}=(1/n,\ldots,1/n)$.  However, such a small bandwidth
 leads to instability problems. In particular, the variance of $\f{h_{0,\text{rev}}}$ is of order $n^{d-1}$ (see \Cref{eEquVarfh}).
%Let us first explain in details the Direct \CDRodeo{} algorithm (the one addressing conditional density estimation) and its limitations, before describing our method and its improvements.

%%%%%%%
\subsubsection{Our method: the RevDir \CDRodeo{} procedure}\label{sec:RevDir}
\begin{figure}
\begin{center}
\begin{tikzpicture}
\draw[->] (-0.5,0) -- (12,0);
\draw (12,0) node[right] {iterations};
\draw [->] (0,-2) -- (0,3);
\draw (0,3) node[above] {$h$};
\draw(3,3) node {Reverse Step};
\draw [-] (6,-2) -- (6,3);
\draw(9,3) node {Direct Step};
\draw [-] (11,-2) -- (11,3);
\draw (-0.5,0) node[ left] {$h_0$};
\draw (0.5,0) node[below] {-$1$};\draw (0.5,0) node {$|$};
\draw (1,0) node[below] {-$2$};\draw (1,0) node {$|$};
\draw (1.5,0) node[below] {-$3$};\draw (1.5,0) node {$|$};
\draw (6.5,0) node[below] {$0$};\draw (6.5,0) node {$|$};
\draw (7,0) node[below] {$1$};\draw (7,0) node {$|$};
\draw (7.5,0) node[below] {$2$};\draw (7.5,0) node {$|$};
\draw (11,0) node[below left] {end};
\draw[red, very thick] [domain=0:3] plot  (\x, {3^(\x/3)-1});
\draw[red, very thick] (3,2) -- (11.1,2);
\draw[red] (11,2) node[right] {$\hat h_j$};
\draw[red] (3,0) node[below] {$t_j$};
%\draw[red] (3,0) node {$|$};
\draw[red, dashed] (3,-0.1) -- (3,2);
\draw[blue, very thick] (0,0) -- (6.5,0);
\draw[blue, very thick] [domain=6.5:9] plot  (\x, {3^(-(\x-6.5))-1});
\coordinate (A) at (9,-0.93585);
\coordinate (B) at (11.1,-0.93585);
\draw[blue, very thick] (A) -- (B) ;
\draw[blue] (B) node[right] {$\hat h_k$};
\draw[blue] (9,0) node[above] {$t_k$};
\draw[blue, dashed] (9,0.1) -- (A);
%legende
\node[draw,text width=2cm] at (1.5,-1.5){\textcolor{red}{$j\in \Act^{(-1)}$}
\textcolor{blue}{$k\in \Act^{(0)} $}};
\end{tikzpicture}
\caption{The two patterns of bandwidth path: the components $j\in \Act^{(-1)}$ with a deactivation time $t_j\leq 0$ in red, and in blue the components $k\in \Act^{(0)}$ with a deactivation time $t_k\geq 0$.
}\label{figrodeo}
\end{center}
\end{figure}
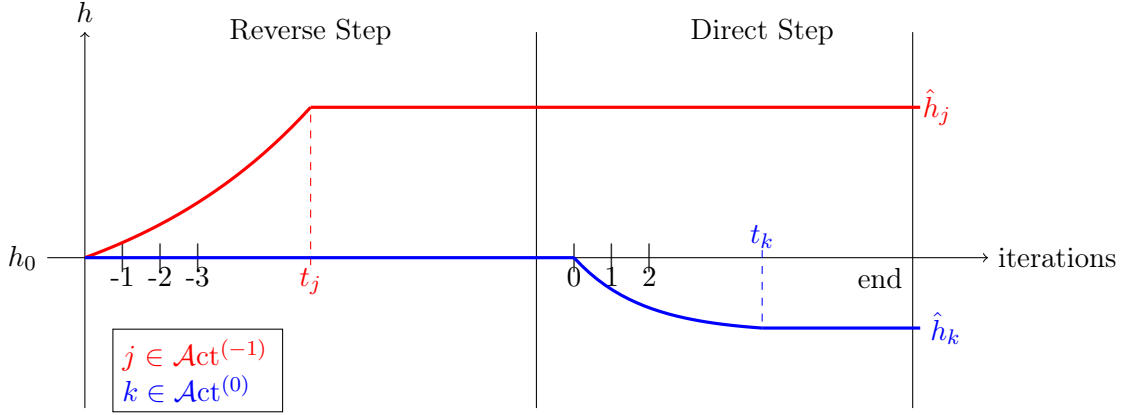

{
In view of the analysis led in Section~\ref{Init}, % shows that, to circumvent previous issues, we have to combine Direct and Reverse \CDRodeo{} procedures, leading to the {\it RevDir} \CDRodeo{} procedure. 
we propose to give the option for each bandwidth component to either increase or decrease.  The procedure is precisely described by Algorithm \ref{algo}. 
The initial bandwidth can then be chosen at an intermediate level (and we show later that $h_0$ has to be chosen larger than the relevant components of the minimax bandwidth), then
our procedure comprises the two following steps: 
%after fixing the initial bandwidth whose components are all equal to $h_0$, and  

\begin{enumerate}
\item The first step consists in the execution of a Reverse \CDRodeo{} procedure to increase the bandwidth components that need to be increased (including the irrelevant ones).
\item The second step executes a Direct \CDRodeo{} procedure on the other bandwidth components. 
\end{enumerate}
}
\begin{algorithm*}[h!]
\caption{RevDir \CDRodeo{}  algorithm \label{algo}}
\begin{enumerate}
\item \textit{Input:} the estimation point $w$, the observations $W$, the bandwidth decreasing factor $\beta\in (0,1)$,  the bandwidth initialization value $h_0>0$, a tuning parameter $a>1$.
\item \textit{Initialization:}{\begin{enumerate}		
		\item[$\triangleright$] Initialize the trial bandwidth: \textit{for }$k\in 1:d$, $H_k^{(0)}\leftarrow h_0$.
		\item[$\triangleright$] Determine which variables are active for the Reverse Step or for the Direct Step: 
		\begin{itemize}
		\item[] $\Act^{(-1)} \leftarrow\lbrace k\in 1:d, \, \vert Z_{H^{(0)}k} \vert \leq \lambda_{{H^{(0)}k}} \rbrace$
		\item[] $\Act^{(0)} \leftarrow \lbrace1:d\rbrace \setminus \Act^{(-1)}$
		\end{itemize}
		%\item Initialize the interation counter: $t=0$.
	\end{enumerate}
	}
%%%%%%%%%%%%%%%%	
\item	\textit{Reverse Step:}
\begin{itemize}
	\item[$\triangleright$] Initialize the counter: $t\leftarrow -1$
	\item[$\triangleright$] Initialize the current bandwidth: $\hat{h}^{(-1)}\leftarrow H^{(0)}$
	\item[$\triangleright$] While ($\Act^{(t)}\neq\emptyset $) \& ($\max \hat{h}_k^{(t)}\leq \beta) : $
	\begin{enumerate}
		\item[$\blacktriangleright$] Set the current trial bandwidth:
		$H^{(t)}_k=\left\lbrace \begin{tabular}{@{}ll}
		$\ibeta \hat{h}_k^{(t)}$ &$\text{ if }k\in\Act^{(t)}$  \\ 
		$\hat{h}_k^{(t)}$ &$\text{ else.}$  \\ 
		\end{tabular}\right. $
		\item[$\blacktriangleright$]  Set the next active set: 
		$\Act^{(t-1)}\leftarrow \lbrace k\in\Act^{(t)}, \, \vert Z_{H^{(t)}k} \vert \leq \lambda_{{H^{(t)}k}} \rbrace$
		\item[$\blacktriangleright$]  Update the current bandwidth:
		$\hat{h}^{(t)}_k\leftarrow \left\lbrace \begin{tabular}{@{}ll}
		$ H_k^{(t)}$ &$\text{ if }k\in\Act^{(t-1)}$  \\ 
		$\hat{h}_k^{(t)}$ &$\text{ else.}$  \\ 
		\end{tabular}\right. $
		\item[$\blacktriangleright$] Initialize the next bandwidth:  $\hat{h}^{(t-1)}\leftarrow\hat{h}^{(t)}$
		\item[$\blacktriangleright$] Decrement the counter: $t\leftarrow t-1$
	\end{enumerate}
\end{itemize}
%%%%%%%%%%%%%%%%%	
%\item	\textcolor{red}{\textit{Reverse Step:}
%\begin{itemize}
%	\item[$\triangleright$] Initialize the counter: $t\leftarrow 0$
%	\item[$\triangleright$] Initialize the current bandwidth: $\hat{h}^{(0)}\leftarrow H^{(0)}$
%	\item[$\triangleright$] While ($\Act^{(t-1)}\neq\emptyset $) \& ($\max \hat{h}_k^{(t)}\leq \beta) : $
%	\begin{enumerate}
%		\item[$\blacktriangleright$] Decrement the counter: $t\leftarrow t-1$
%		\item[$\blacktriangleright$] Set the current trial bandwidth:
%		$H^{(t)}_k=\left\lbrace \begin{tabular}{@{}ll}
%		$\ibeta \hat{h}_k^{(t-1)}$ &$\text{ if }k\in\Act^{(t)}$  \\ 
%		$\hat{h}_k^{(t-1)}$ &$\text{ else.}$  \\ 
%		\end{tabular}\right. $
%		\item[$\blacktriangleright$]  Set the next active set: 
%		$\Act^{(t-1)}\leftarrow \lbrace k\in\Act^{(t)}, \, \vert Z_{H^{(t)}k} \vert \leq \lambda_{{H^{(t)}k}} \rbrace$
%		\item[$\blacktriangleright$]  Set the current bandwidth:
%		$\hat{h}^{(t)}_k\leftarrow \left\lbrace \begin{tabular}{@{}ll}
%		$ H_k^{(t)}$ &$\text{ if }k\in\Act^{(t-1)}$  \\ 
%		$\hat{h}_k^{(t-1)}$ &$\text{ else.}$  \\ 
%		\end{tabular}\right. $
%		%\item[$\blacktriangleright$] Initialize the next bandwidth:  $\hat{h}^{(t-1)}\leftarrow\hat{h}^{(t)}$
%	\end{enumerate}
%\end{itemize}
%}
%%%%%%%%%%%%%%%
\item \textit{Direct Step:}
	\begin{itemize}
	\item[$\triangleright$] Initialize the current bandwidth: $\hat{h}^{(0)}\leftarrow \hat{h}^{(t)}$
	\item[$\triangleright$] Reinitialize the counter: $t\leftarrow 0$
	\item[$\triangleright$] While $\left(\Act^{(t)}\neq\emptyset \right) \& \left(\prod\limits_{k=1}^d \hat{h}_k^{(t)}\geq \frac{(\log n)^{1+a}}{n}\right)$: 
	\begin{itemize}
		\item[$\blacktriangleright$] Increment the counter: $t\leftarrow t+1$
		\item[$\blacktriangleright$] Set the current active set: 
		$\Act^{(t)}\leftarrow \lbrace k\in\Act^{(t-1)}, \, \vert Z_{\hat{h}^{(t-1)}k} \vert > \lambda_{{\hat{h}^{(t-1)}k}} \rbrace	$
		\item[$\blacktriangleright$] Set the current bandwidth: 
		$\hat{h}^{(t)}_k\leftarrow \left\lbrace \begin{tabular}{@{}ll}
		$ \beta.\hat{h}_k^{(t-1)}$ &$\text{ if }k\in\Act^{(t)}$  \\ 
		$\hat{h}_k^{(t-1)}$ &$\text{ else.}$  \\ 
		\end{tabular}\right. $
	\end{itemize}
\end{itemize}
%%%%%%%%%
\item \textit{Output:} $\hat{h}\leftarrow \hat{h}^{(t)}$ (and compute $\hat{f}_{\hat{h}}(w)$).
\end{enumerate}
\end{algorithm*}
The output bandwidth of the algorithm is denoted by $\hat h$, and the estimator of $f$ by $\hat f:=\hat f_{\hat{h}}$. 
Figure \ref{figrodeo} illustrates the two kinds of path for the bandwidth components. If the component belongs to $\Act^{(-1)}$ (resp. $\Act^{(0)}$), it is deactivated during the Direct Step (resp. the Reverse Step) and has to be chosen larger (resp. smaller) than the initial bandwidth value $h_0$.
Note that the RevDir procedure generalizes both the Direct and Reverse procedures in function of the choice of $h_0$. Indeed, if we set $h_0=1$, the RevDir procedure behaves as a Direct procedure with the same initialization. Conversely, setting $h_0=1/n$ brings us back on the Reverse procedure. {Nonetheless, the purpose of our approach is to provide a better tuning of $h_0$, as discussed in the next section, to solve the initialization issue of the Direct and Reverse procedures. }
%However, its choice is less critical as for the Direct \CDRodeo{} procedure as soon as it is larger than $h_j^*$ for any $j\in{\mathcal R}$.

%\textcolor{red}{\rule{15cm}{1mm}}

%
%{\color{forest}
%The major goal of Direct  \CDRodeo{} is to select a bandwidth $h$ satisfying
%$B^2(h)\approx \mathrm{V}_h$, 
%which is ensured (under some specific assumptions) by stopping the procedure when $\vert Z_{hj}\vert \approx \seuil{hj}$ for a well chosen threshold.
%}
%
%
%{\color{forest}
%In the next paragraph, we explain heuristically how this algorithm is able to simultaneously detect the irrelevant components and select the suitable bandwidth for relevant and irrelevant components.
%}
%
%\Rodeo{}-type algorithms have mainly been studied with decreasing path, setting the iterative step factor $\beta<1$. We call these decreasing algorithms \emph{Direct} in the sequel. This one-way path limits the choice for initialization of the bandwidth, which affects the convergence rates of these methods.

%%%%%%%%%%%%%%%%%%%%%%%%%%
%%%%%%%%%%%%%%%%%%%%%%%%%%
\section{Theoretical results}\label{sec:resul}
%%%%%%%%%%%%%%%%%%%%%%%%%%
\subsection{Sparsity and smoothness classes of functions}
\label{sec:sparsesmooth}
This section is devoted to the theoretical results satisfied by the RevDir \CDRodeo{}  procedure. We consider a kernel function $K:\R\rightarrow \R$ of class ${C}^1$, with compact support denoted $\supp(K)$. We shall also assume that $K$ is
of order $p$, \textit{i.e.}: for $\ell\in 1:p-1$, $\int_{\R} t^\ell K(t)dt=0$.
Taking a kernel of order $p$ is usual for the control of the bias of the estimator.
Then, we define the neighborhood $\UUn$ of the point $w\in\R^{d}$   as follows: 
 $$\UUn:=\left\lbrace u\in\mathbb{R}^{d}: w-u\in\left(\supp(K)\right)^{d}\right\rbrace.$$
 In the sequel, we denote
 $$\Vert f\Vert_{\infty,\UUn}:=\sup_{x\in\UUn}|f(x)|.$$
\begin{remark}\label{remU}
The size of $\UUn$ is fixed. But $\UUn$ could be chosen so that its size goes to 0. In this case, we have to modify the stopping rule of the Reverse Step, namely $\max \hat{h}_k^{(t)}\leq \beta$, to force $\max \hat{h}_k^{(t)}\overset{n\rightarrow \infty}{\longrightarrow} 0$. For instance, if we impose $\max \hat{h}_k^{(t)}\leq \frac1{\log n}$, the rates of convergence of our estimate would typically be deteriorated by a logarithmic term.
\end{remark}
The notion of relevant components has already been introduced in Section~\ref{sec:OMC} but subsequent results only need the function $f$ to be locally sparse, so we shall consider the following definition depending on $\UUn$.
\begin{definition}\label{def:relevant}
We denote $\mathcal{R}$ the subset of $\lbrace0,\dots,d\rbrace$ with cardinal $r$ such that for any fixed $\lbrace z_j\rbrace_{j\in\mathcal{R}}$, the function $\lbrace z_k\rbrace_{k\in\mathcal{R}^c}\mapsto f(z_1,\dots,z_d)$ is constant on $\mathcal{U}$. We call \textit{relevant} any component in $\mathcal{R}$.
\end{definition}
The previous definition means that on $\UUn$, $f$ depends only on $r$ of its $d$ variables. %The notion of relevant component is local: it depends on the point $w$ where $f$ is estimated. 
\begin{remark}  In fact, \CDRodeo{} detects more complex sparsity structures. In particular, $\rel^c$ could be enlarged to the components which are polynomial of degree smaller than the order $p$ of the kernel, namely it suffices to consider $f(z)=z_j^l g(z_{-j})$ with $ l \in 0: p-1$,  $z_{-j}=(z_1, \dots, z_{j-1}, z_{j+1}, \dots, z_d)$ and where $g$ is an arbitrary function.  Then $j$ is considered as an irrelevant component by both our algorithm and in the bias-variance trade-off.  Indeed, assume that $\hfX=\fX$  for the sake of simplicity.  Then  for the algorithm, easy computations leads to $\E[Z_{hj}]=0$  as an irrelevant bandwidth component. Then our algorithm behaves exactly as if $j$ were irrelevant and select a large $\hat{h}_j$ (with high probability).  For the bias-variance trade-off, the bias for $f$ is proportional to the bias for $g$ (multiplied by a term that does not depend on $h$):
$$(\K_h\star f - f) (w)=w_j^l \:(\K_h^{(-j)}\star g-g) (w_{-j}), \qquad \K_h^{(-j)}(z_{-j}):=\prod_{k\neq j}K_h(z_k),$$
exactly as if $j$ were irrelevant and $f(z)= cg(z_{-j})$, for $c$ a constant. Then since only the variance depends on $h_j$,
%Since there is no bias in direction $j$, 
the bias-variance trade-off chooses a large value for $h_j$. \\
In particular there is no need for preliminary linear variable selection as suggested in Section~6.1 of \citep{LW08}.
\end{remark}
In the sequel %,
% we consider the minimax point of view and 
we derive rates on  H\"older balls defined as follows.
\begin{definition}\label{def:holder}
Let $L>0$ and  $s>0$. We say that the conditional density $f$ belongs to the  H\"older ball of smoothness $s$ and radius $L$, denoted $\mathcal{H}_d(s, L)$, if $f$ is of class $C^q$ and if it satisfies for all $z\in \mathcal{U}$  and  for all  $t \in \R$ such that $z+te_k\in \mathcal{U}$
$$\left|\partial_k^{q}f(z+te_k) - \partial_k^{q}f(z)\right|\leq L |t|^{s- q},$$
%\begin{eqnarray*}
%&&\mathcal{H}_d(s, L)=\Big\{ f: \R^d\to \R \text{ of class }C^q \text{ such that  for all } z\in \mathcal{U}, \\
%&& \text{ for all }t \in \R \text{ such that  }z+te_k\in \mathcal{U}\quad \left|\partial_k^{q}f(z+te_k) - \partial_k^{q}f(z)\right|\leq L |t|^{s- q}
%\Big\}\\
%\end{eqnarray*}
where  $q=\lceil s-1 \rceil=\max\{l\in\mathbb{N}: l<s\}$ and $e_k$ is the vector where all coordinates are null except the $k$th one which is equal to 1. 
\end{definition}
We investigate adaptive results in terms of sparsity and smoothness properties on H\"older balls $\mathcal{H}_d(s, L)$, with $s>1$. 
Adaptation means that our procedure will not depend on the knowledge of $\mathcal{R}$ and $(s,L)$. 
The condition on $s$ means that $f$ has to be at least $C^1.$ This technical assumption is related to our methodology based on derivatives of $\f{h}$ as proxies of derivatives of $f$ to detect relevant components.
%%%%%%%%%%%%%%%%%%%%%%%%%%%%%%%%
%\subsection{Tuning the RevDir \CDRodeo{} procedure}%\label{sectionRodeoParameter}
%%%%%%%%%%%%%%%%%%%%%%%%%%%%%%
\subsection{Range of the algorithm inputs and assumptions}\label{sssectAssumpResult}
The RevDir \CDRodeo{} procedure depends on three tuning parameters, namely $h_0$, $\beta$ and $a$.
In the sequel, we take $\beta\in(0,1)$. {
%\textcolor{forest}{ The smaller the decay factor $\beta$, the bigger the step size, thus the faster the procedure but the larger the approximation error. }
Since $\beta$ is an exponential decay factor, its value has no influence on rates of convergence (up to the constant factor). 
%\textcolor{magenta}{But of course, the larger $\beta$, the more accurate the procedure, but the larger the computational time. }
%
%In practice, we set $\beta$ close to $1$.
%
 The parameter $a$ will be assumed to be larger than 1. Its value does not affect the main polynomial factor $n^{-\frac{s}{2s+r}}$ of the rate of convergence but only the logarithmic factor: the smaller $a$, the smaller the exponent of the logarithmic factor. }%In practice, $a$ will be larger but close to 1.
See Section~\ref{sec:calib} for a detailed analysis of the practical choices for $a$ and $\beta.$
Finally, to initialize the procedure, we take $h_0$ such that 
\begin{equation}\label{ineq-initialisation}
\Cl^{2/d}\left(\frac{(\log n)^a}{n}\right)^{\frac{1}{d(2p+1)}}\leq h_0\leq 1, 
\end{equation}
where $\Cl$, only depending on the kernel $K$, is defined in Section~\ref{sec:heuristic}. Note in particular that the lower bound does not depend on any unknown value, and thus can be implemented as the bandwidth initialization.
Besides, observe that each component of the minimax bandwidth for estimating $f$ on $\mathcal{H}_d(s, L)$ is of order $n^{-1/(2s+r)}$ for relevant components and are constant for irrelevant ones. So, if $s\leq p$, as assumed in Theorem~\ref{thm}, then $h_0$ is larger than all relevant components of the optimal bandwidth, as required by the RevDir \CDRodeo{} procedure.

%%%%%%%%%%%%%%%%%%%%%%%%%%%%%%%%%
%\subsection{\CDRodeo{} RevDir parameters choice.} \label{sectionRodeoParameter}
%\paragraph{Kernel $K$.}\label{KC} 
% We choose the kernel function $K:\R\rightarrow \R$ of class $C^1$, with compact support and of order $p$, \textit{i.e.}: for $\ell=1,\dots,p-1$, $\int_{\R} t^\ell K(t)dt=0$.\\
%Taking a kernel of order $p$ is usual for the control of the bias of the estimator.
%\paragraph{Parameter $\beta$.} Let  $\beta\in(0,1)$ be the decreasing factor of the bandwidth. The larger $\beta$, the more accurate the procedure, but the longer the computational time. From the theoretical point of view, it remains of little importance, as it only affects the constant terms. In practice, we set it close to $1$.
%
%
%
%\paragraph{Bandwidth initialization.} 
%We choose
%$$1\geq h_0\geq\Cl^{1/d}\left(\frac{(\log n)^a}{n}\right)^{\frac{1}{d(2p+1)}}$$
%with
%$\Cl:= 4\Vert J\Vert_2\Vert K\Vert_2^{d-1}$ (where we recall the definition of $J:t\mapsto K(t) + tK'(t)$).
%
%
%\paragraph{Threshold $\seuil{h,j}$.} For any bandwidth $h\in(\R_+^*)^d$ and for $j=1:d$, we set the threshold as follows:
%\begin{equation}
%\seuil{hj}:=\Cl\sqrt{\frac{(\log n)^{a}}{n h_j^2 \prod_{k=1}^d h_k}},
%\end{equation}
%with
%$\Cl= 4\Vert J\Vert_2\Vert K\Vert_2^{d-1}$.
%%(where we recall the definition of $J:t\mapsto K(t) + tK'(t)$).
%The expression is obtained by using concentration inequalities on $Z_{hj}$. 
%
%~\\
%Hereafter, unless otherwise specified, the parameters are chosen as described in this section. 
To derive rates of convergence for $\hat f(w)$, we need three assumptions. The first two ones are related to $\fX$, the density of the $X_i$'s. We recall that the evaluation point is $w = (x, y)$.
\begin{assumption}{\AfXmin\label{AfXmin}}
{[Lower bound on $\fX$]\\
The density $\fX$ is bounded away from 0 in the neighborhood of $x$: $$
\delta:=\inf_{u \in \mathcal{U}_1}\fX(u) >0,$$
where $\mathcal{U}_1:=\left\lbrace u\in\mathbb{R}^{d_1}: x-u\in\left(\supp(K)\right)^{d_1}\right\rbrace.$
}
\end{assumption}
\begin{remark}\label{remU1}
 Similarly to Remark~\ref{remU}, the size of $ \mathcal{U}_1$ is fixed but it could decrease to 0 if we modify the stopping rule of the Reverse Step.
\end{remark}
This assumption  is classical in the regression setting or for conditional density estimation. Indeed,  if $\fX$ is equal or close to 0 in the neighborhood of $x$, we  have no or very few observations to estimate the distribution of $Y$ given $X=x$. Thus, this assumption is required in all of the aforementioned works about conditional density estimation.

\bigskip

The next assumption specifies that we can estimate $\fX$ very precisely.
\begin{assumption}{\AestfX}
{[Estimation of $\fX$]\\
The estimator of $\fX$ in \eqref{Deffh} satisfies the following two conditions: 
\begin{enumerate}[label=Condition (\roman*)]
  \item a positive lower bound: $\dX:=\inf\limits_{u\in \Un} \hfX(u)>n^{-1/2}$, \label{fXtildemin}
  \item a concentration inequality in local sup norm: \label{fXtildeAccuracy}
    $$\P\left(\sup\limits_{u\in\Un}\left\vert \fX(u)-\hfX(u)\right\vert > M_{X} \frac{(\log n)^{\frac{a}2}}{\sqrt{n}}\right)
  \leq  \exp(-(\log n)^{1+\frac{a-1}2}),$$
  \end{enumerate}
  with $M_{X}:= \frac{\delta\Vert J\Vert_2 \Vert K\Vert_2^{d-1}}{4\Vert f\Vert_{\infty,\mathcal{U}}\Vert J\Vert_1 \Vert K\Vert_1^{d-1}}$.  
  }
\end{assumption}
\begin{remark}
For the simpler problem of density estimation, since $\fX\equiv 1\equiv\hfX$,  Assumption~\AestfX{} is obviously satisfied.
\end{remark}
{
This $\sqrt{n}$-rate can be achieved either by restricting $\fX$ to a parametric class,  or by assuming we have at hand a larger sample of $X$. In particular, the following proposition provides precise conditions to satisfy Assumption \AestfX{} using a well-tuned kernel density estimator $\hfX$.}
Furthermore, $\hfX$, the estimator provided by the proof of Proposition~\ref{propfXtilde}, is easily implementable.
\begin{proposition}\label{propfXtilde} Given a sample $\tilde{X}$ with same distribution as $X$ and of size $\nX=n^c$ with  $c>1$,
if $\fX$ is of class ${C}^{p'}$ with $p'\geq \frac{d_1}{2(c-1)}$, there exists an estimator  $\hfX$ which satisfies Assumption \AestfX.
\end{proposition}
To prove Proposition~\ref{propfXtilde}, we build $\hfX$ as a truncated kernel estimator with a fixed bandwidth, but other methods can be used in practice, as, for instance, a Rodeo algorithm for density estimation. Actually
any reasonable nonparametric estimator would have a rate of convergence in sup norm of the form $n_X^{-\beta}$ (typically $\beta=p'/(2p'+d_1)$) up to a logarithmic term. Then Condition (ii) of Assumption~\AestfX{} is verified as soon as $n_X^{-\beta}\leq n^{-1/2}$ and we need $c\geq 1+d_1/(2p')$. Then, observe that if $\fX$ is of class ${ C}^{\infty}$, then we just need $c=1$ and we can take $\tilde X=X$. If we know that $\fX$ is at least of class ${ C}^1$ but its precise smoothness is unknown, taking $c\geq  1+d_1/2$ is sufficient to satisfy assumptions of Proposition~\ref{propfXtilde}. 
% which imposes a condition $n_X\geq n^{c}$ as in Proposition~\ref{propfXtilde}. Here $c\geq  1+d_1/2$ is sufficient for any density of class $C^1$. 

\bigskip

The next assumption is necessary to control the bias.
\begin{assumption}{\Monoton\label{Monoton}}
{%[Convexity]
\\
For all $j\in\rel$, for all $h$ and $h'\in[\frac1n, 1]^d$ such that $h\preceq h'$, 
$\vert\mathbb{E}{[}\bar Z_{h,j}{]}\vert 
\leq \vert \mathbb{E}{[}\bar Z_{h',j}{]} \vert $,
where $\bar Z_{h,j}$ is defined as $Z_{h,j}$ in \eqref{defZ} but with true $\fX$ replacing $\hfX$.}
\end{assumption}

{
Let us comment Assumption~\Monoton{}.  First observe that it is verified by the sharp bound of $\E[\bar Z_{hj}]$ over the class $\mathcal{H}_d(s, L)$ ($s>1$): denoting $M_j$ the pseudo-kernel defined by $M_j(z)= J(z_j)\prod_{k\neq j}K(z_k),$
\begin{align}\label{convMj}
\vert \E[\bar Z_{hj}]\vert&= \vert \tfrac{\partial}{\partial h_j} ( \K_h\star f - f)(w)\vert
%&=-\frac1{h_j}  M_{h}^j \star f(w)\\
=\vert \tfrac{1}{h_j}\int M_j(z) [f(w-h.z) -f(w)]dz\vert 
%&=-\frac{c_K}{h_j}\left(\sum_kh_k^q(\partial_k^q f)(w)+o(\|h\|^q)\right)\\
\lesssim h_j^{-1} \sum_{k=1}^d h_k^s,
\end{align} 
the last inequality coming from Taylor expansion and the $s$-Hölder smoothness of $f$.  
Then for all $h, h'\in[\frac1n, 1]^d$ such that $h\preceq h'$,  $h_j^{-1}\sum_{k=1}^d h_k^s \leq {h'_j}^{-1}\sum_{k=1}^d {h'_k}^s. $

Assumption~\Monoton{} is named after convexity or concavity, as it requires monotony of $\E[\bar Z_{hj}]$, which is the derivative of the bias (after removing the potential perturbations of the pre-estimator $\hfX$ by replacing $\hfX$ by $\fX$).  The absolute value in the assumption is simply a way to cover both cases (convexity and concavity), since in fact $\E[\bar Z_{hj}]\rightarrow 0$ as $h\rightarrow 0$ (at least in the scope of our results: when the smoothness $s$ is larger than $1$).  In the context of the algorithm, this assumption prevents $\E[Z_{hj}]$ from vanishing temporarily and thus the algorithm from stopping prematurely.  Ensuring that, the bias-variance trade-off is achieved. 

Note that otherwise,  the non convexity of the squared bias would reverberate on the squared risk,  making its minimization much harder, especially when we target greedy algorithms to avoid a computationally intensive optimization over all bandwidths.}

%\textcolor{magenta}{
%\sout{Assumption~\Monoton{} est equivalente a la monoticit\'e de $\E[\bar Z_{hj}]$ car $\E[\bar Z_{hj}]$ tend vers 0 lorsque $h$ tend vers 0,  $\E[\bar Z_{hj}]$ ne change plus de signe. Cette monotonie concorde avec la monotonie des majorants sharps habituels. Donc on a  vraiment la convexite de B(h) (ou de -B(h)). Si le biais n'etait pas convexe, le risque ne le serait peut-etre pas, et la recherche du minimum du risque serait rendue complique par l'existence potentielle de plusieurs minima locaux.}} 

\begin{remark}
If $f$ is smooth enough so that $\frac{\partial^p }{\partial h_j^p} f(h)\not=0$ with $p$ such that $\int u^pK(u)du\not=0$, then Assumption~\Monoton{} is not required. Nevertheless, the procedure cannot be adaptive in this case. See \cite{Jeanne1}.
\end{remark}

%%%%%%%%%%%%%%%%%%%%%%%%
\subsection{Main result}\label{sec:mainresult}
We now derive the main result of our paper proved in Section~\ref{sec:proofs} in which we show that $\hat h$ is closed to the ideal bandwidth $h^*$ defined in Section~\ref{sec:heuristic}. Thus our algorithm is able to both detect the irrelevant components and select the minimax bandwidth for relevant and irrelevant components.

\begin{theorem}\label{thm} For any $r\in 0:d$, $1<s\leq p$ and $L>0$, if $f$ has only $r$ relevant components  and belongs to $\mathcal{H}_d(s, L)$, then
under Assumptions \AfXmin, \AestfX and \Monoton, the pointwise risk of the RevDir \CDRodeo{}  estimator $\frod(w)$ is bounded as follows: for any $l\geq 1$, for $n$ large enough, 
\begin{equation}\label{eq:rate}
\E\left[
\left\vert
\frod(w)-f(w)
\right\vert^l
\right]^{1/l}
\leq  \text{C} 
 	\left(\frac{(\log n)^a}{n}\right)^{\frac{s}{2s+r}} 
\end{equation}
where $\text{C}$ only depends on $d,r,K, \beta,\delta, L, s,\|f\|_{\infty,\UUn}$. 
%{\color{violet} Sûre que ça ne dépend pas de $l$}
\end{theorem}
We can compare the obtained rate with the classical pointwise  adaptive minimax rate for estimating a $s$-regular $r$-dimensional density, which is $((\log n)/n)^{s/(2s+r)}$ (see  \cite{rebelles15}). Our procedure achieves this rate up to the term $(\log n)^{s(a-1)/(2s+r)}$. In Section~\ref{sssectAssumpResult}, we specify that any value $a>1$ is suitable. So, our procedure is nearly optimal. Actually, we need $a>1$ to ensure that for $n$ large enough,
\begin{equation}\label{a=1}
(\log n)^{a-1}\geq \frac{\Vert f\Vert_{\infty, \UUn}}{\delta}
\end{equation}  but if an upper bound (or a pre-estimator) of $\frac{\Vert f\Vert_{\infty, \UUn}}{\delta}$ were known, we could obtain the similar result with $a=1$, and our procedure would be rate-optimal without any additional logarithmic term. 
Remember that the term $(\log n)^{s/(2s+r)}$ is the price to pay for adaptation with respect to the smoothness (see  \cite{tsybakov98}). Theorem~\ref{thm} shows that, in our setting, there is no additive price for not knowing the sparsity, i.e. the value of $r$. This result is new for conditional density estimation.
\begin{remark} 
%\textcolor{red}{Claire est ok avec ce qui suit}
Assumption~\Monoton{} allows for a sharp control of the bias of our estimate and is only used in Section~\ref{sec:eMajBias}.  Refining the decomposition of the term $\bar B_h$ in 
\eqref{eDecompBiashint} shows that we can relax Assumption~\Monoton{}.
This is done in the supplementary file  \citep{JCV-supp} where Assumption~\Monoton{} is  replaced by Assumption~\Monotonbis.
The price to pay is an extra logarithmic term $(\log n)^{\frac{2s}{2s+r}}$ in the upper bound~\eqref{eq:rate}. 
% by replacing it with the following one: For all $j\in\rel$, for all $h$ and $h'\in(\mathbb{R}^*_+)^d$ such that $h\preceq h'$, 
%$h_j\vert\mathbb{E}{[}\bar Z_{h,j}{]}\vert 
%\leq h'_j\vert \mathbb{E}{[}\bar Z_{h',j}{]} \vert $. Observe that $h_j\vert\mathbb{E}{[}\bar Z_{h,j}{]}\vert$ has a simple expression, namely
%$$h_j\vert\mathbb{E}{[}\bar Z_{h,j}{]}\vert=\left|\int M_j(z) [f(w-h.z) -f(w)]dz\right|.$$
%The price to pay is an additive logarithmic term in the minimax rate of convergence (the right hand side of \eqref{eq:rate}).
 \end{remark}
%}

%As already noticed in \citep{Jeanne1}, some additional logarithmic factors could be unavoidable due to 
%the unknown sparsity structure, which needs to be estimated. 
%\textcolor{violet}{Il faudrait plutôt mettre en lumière le gain du facteur $(\log n)^{d-r}$ par rapport à \citep{Jeanne1}, obtenu grâce à la procédure RevDir, surtout que je ne pense plus que cette phrase soit vraie, puisque maintenant si $a=1$ on obtient la puissance de log optimal (cf \cite{rebelles15}).\\ 
% À potentiellement rajouter : le commentaire sur le $a$ dans la puissance du log (ci-dessous, que j'ai repris de \citep{Jeanne1}).}
%%%%%%%%%%%%%%%%%%%%%%%%%%%%
\subsection{Algorithm complexity}\label{sec:complexity}
We now discuss the complexity of \CDRodeo{}. %
%without taking into account the pre-computation cost of $\hfX$ at the points $X_i$, $i=1:n$ (used for computing the $Z_{hj}$). %but a fast procedure for $\fX$ is required, to avoid losing \CDRodeo{} computational advantages. \textcolor{blue}{typiquement un rodeo}
Regarding the computation cost of $\hfX$, the estimator built for the proof of Proposition~\ref{propfXtilde} has complexity $O(d_1 n^c)$ but in practice we use a \textsc{Rodeo} estimator with the same sample size $n$, which has a complexity $O(d_1 n\log n)$ for each computation of $\hfX(X_i)$ which causes an additional cost in $O(d_1 n^2\log n)$ {(applying following \Cref{complexity})}.

Regarding the main part of the algorithm, during the Reverse Step, $|\Act^{(-1)}|$ components are updated, and, for fixed $h$,
the computation of all $Z_{hj}$'s  and the comparisons to the thresholds $\seuil{hj}$ need $O(|\Act^{(-1)}|n)$ operations.
In the same way, during the Direct Step, $|\Act^{(0)}|$ components are updated and each update needs $O(|\Act^{(0)}|n)$ operations.
Since the number of updates is at worse of order $\log(n)$ (because of the stopping conditions), and $|\Act^{(-1)}|+|\Act^{(0)}|\leq d$, we obtain the following proposition. More details can be found in the proof (see Section~\ref{proofcomplexity}).
%
%
%
%Number of updates of the bandwidth (for all components): 
%\begin{itemize}
%\item Reverse Step : $\frac{|\Act^{(-1)|}{d(2p+1)}\log(n)$ because of the stopping condition and the value of $h_0$
%\item Direct Step: $\frac{1}{\log (1/\beta)}\log(n)$ because of the stopping condition
%\end{itemize}
%
%\medskip
%Each time : 
%\begin{itemize}
%\item
%computation of $Z_{hj}$ : $O(n?)$ operations
%\item 
%computation of the threshold and comparison  : $O(d)$ operations
%\end{itemize}
%
%Using that $|\Act^{(-1)}|+|\Act^{(0)}|\leq d$, the sum of these two steps leads to the following proposition.

\begin{proposition}\label{complexity}
Apart from the computation of $\hfX$, the total worst-case complexity of RevDir \CDRodeo{}  algorithm is 
$$O(d n\log n).$$
\end{proposition}

Notice that for classical methods with optimization on a bandwidths grid, the complexity is of ordrer $dn|H|^d$, where $|H|$ denotes the size of the grid for each component. In practice, the grid has to include at least $\log n$ points, which leads to a computational cost $O(d n(\log n)^d)$. For illustration, $d=5$ and $n=10^5$, the ratio of complexities {$\frac{d n(\log n)^d}{d n\log n}$ is already larger than}
$1.7 \times 10^4$.
%%%%%%%%%%%%%%%%%%%%%%%%%
%%%%%%%%%%%%%%%%%%%%%%%%%
\section{Simulations}\label{sec:simus}
This section is devoted to the numerical analysis of our algorithms.  
In Section~\ref{sec:examples}, we first describe the three examples on which we test \CDRodeo{}. Then we calibrate its parameters in Section~\ref{sec:calib}.  
{We finally look at its numerical performances in Section~\ref{sec:perf}: we first analyse the behavior of \CDRodeo{} for different examples then assess the sparsity detection by adding an increasing number of irrelevant components. In particular, our analysis relies on the fact that the behavior of \CDRodeo{} is easily explainable from the bandwidth it selects.}

%%%%%%%%%%%%%%%%%%%%%%%%%
%Note also that considering the multivariate nonparametric setting, we consider large sample sizes (of the order of  several thousands) in order to assess adequately the performances of \CDRodeo{}. Remember that

%{\textcolor{forest}{Sout?} \textcolor{forest}{on en parle déjà à un autre endroit: 3e paragraphe de 4.3.2, mais sans la quantification avec le epsilon}\textcolor{magenta}{\indent Simulations will be performed in the nonparametric setting and for (moderately) large dimensions (typically for $d$ between $2$ and $13$), so, to assess the performances of \CDRodeo{} conveniently, samples sizes will be of order of several thousands. Indeed considering minimax rates, for a precision $\epsilon$ and a smoothness $s$, the sample size has to be larger than $\epsilon ^{-\frac{2s+d}{s}}$. 
 %Therefore, running times cannot be very moderate (see Proposition~\ref{complexity}) \textcolor{forest}{??}. It is important to note even if $n$ is large, our algorithm performs in reasonable running times.}
%Since simulations will be performed for moderately large dimensions (typically for $d$ between 2 and 8), we will take samples whose size is in the order of several thousand.
%}
%%%%%%%%%%%%%%%%%%%%%%%%%
\subsection{Examples}\label{sec:examples}
% \textcolor{red}{Ces 4 lignes sont sans doute inutiles} Before describing the tested examples,
%we fix the notations: we consider a  $n\times d$ matrix $W$ of $n$ i.i.d.~observations of dimension $d$. We decompose the observation $W_i$ ($i$th row of $W$)  as $W_i=(X_i, Y_i)$, where $X_i$ is the auxiliary vector of dimension $d_1$ and $Y_i$ the response vector of dimension $d_2$ (so that $d=d_1+d_2$).
%
We describe 3 examples. For this purpose, we denote ${\mathcal N}(a,b)$ the Gaussian distribution with mean $a$ and variance $b$, ${\mathcal U}_{[a,b]}$ the uniform distribution on the compact set $[a,b]$ and ${\mathcal IG}(a,b)$ the inverse-gamma distribution with parameters $(a,b)$. 
\begin{itemize}
\item\underline{Example (a):} We consider $d_2=2$ response variables and $d_1\in 1:4$ auxiliary variables with the following hierarchical structure:
$$Y_{i2}\sim {\mathcal IG}(4,3),\quad Y_{i1}|Y_{i2}\sim{\mathcal N}(0,Y_{i2}),\quad X_{ij}|Y_i\stackrel{iid}{\sim}  {\mathcal N}(Y_{i1},Y_{i2}),$$
which leads to the following conditional density (derived in \cite[Chapter IV, Section~5.a]{Jeanne-these}):
$$f:(x,y)\mapsto \1_{\{y_2>0\}}\frac{\sqrt{d_1+1}}{\sqrt{2\pi}\Gamma(4+\frac{d_1}{2})}(\beta_1(x))^{4+\frac{d_1}{2}}y_2^{-(5+\frac{d_1+1}{2})}
e^{-\frac{\beta_1(x)}{y_2}-\frac{\left(y_1-\frac{\sum_{j=1}^{d_1}x_j}{d_1+1}\right)^2}{\left(\frac{2y_2}{d_1+1}\right)}}
$$
with $\beta_1(x):=\frac{1}{2}\Big(6+\sum_{j=1}^{d_1}x_j^2-\frac{(\sum_{j=1}^{d_1}x_j)^2}{d_1+1}\Big)$.

This example is an usual Bayesian model (see for example \citep{RaynalMPR18}) where one of the tasks is to retrieve the posterior distribution $f$ of the mean $Y_{i1}$'s and the variance $Y_{i2}$'s given the normal observations $X_i$'s, which is exactly what our method performs in this paper. 

\item\underline{Example (b):} We consider $d_2=1$ response variable and  $d_1\in 1:12$ auxiliary variables with the following hierarchical structure:
$$X_{ij}\stackrel{iid}{\sim} {\mathcal N}(0,1),\quad Y_{i1}|X_i \sim{\mathcal N}(3X_{i1}^3,0.5^2).$$
In this case, the conditional density is then
$$f:(x,y)\mapsto \sqrt{\frac{2}{\pi}}e^{-2(y-3x_1^3)^2}.$$
\item\underline{Example (c):} We consider $d_2=1$ response variable and  $d_1\in 1:12$ auxiliary variables with the following hierarchical structure:
$$X_{ij}\stackrel{iid}{\sim} {\mathcal U}_{[-1,1]},\quad Y_{i1}|X_i \sim{\mathcal N}(3X_{i1}^3,0.5^2).$$
In this case, the conditional density is then
$$f:(x,y)\mapsto \sqrt{\frac{2}{\pi}}e^{-2(y-3x_1^3)^2}\1_{\{x\in[-1,1]^{d_1}\}}.$$
\end{itemize}
Example (a), in which $r=d$, will be used as reference for estimation without sparsity structure and will illustrate the estimation difficulty when we have to face with the curse of dimensionality. Examples (b) and (c) circumvent the curse of dimensionality given their sparsity structure: $r=2$ ($Y_i$ is scalar and depends only on $X_{i1}$). Note that Example (c) is discontinuous, whereas our method rather targets $C^1$-functions.

%%%%%%%%%%%%%%%%%%%%%%%%%
\subsection{Calibration...}\label{sec:calib}

In this section, we focus on the calibration of the threshold $\seuil{hj}$ and the decay factor $\beta$, whereas some other parameters are fixed: in particular, we are using the Gaussian kernel, and the initialization value $h_0$ is chosen as the lower bound provided in \eqref{ineq-initialisation}:
\begin{equation}\label{eqInit h0}
h_0=\Cl^{2/d}\left(\frac{(\log n)^a}{n}\right)^{\frac{1}{d(2p+1)}}
\end{equation}
 with $\Cl= 4\Vert J\Vert_2\Vert K\Vert_2^{d-1}$. The choice of the threshold is quite sensitive since it influences the bias-variance trade-off, and intensive simulations have been performed to determine the convenient tuning, while the decay factor (which only quantifies the size of the step) rather impacts the running times of the procedure. 
 
 We determine each parameter separately, since their respective impact is rather independent.  Moreover, to avoid the influence of a chosen $\hfX$ (and its peculiar specificities), the calibration is run with known $\fX$ (which is plugged as input of the algorithm).

%%%%%%%%%%%%%%%%%%
\subsubsection{... of the threshold}\label{sec:a}
Since the calibration of $\beta$ is not done yet, we fix for this section $\beta=0.9$. 

Given Definition \eqref{threshold}, two parameters influence the threshold: $a$ and $\Cl$, but they are clearly redundant. Therefore only the calibration of $a$ will be performed while we take the theoretical value of $\Cl$.

We compare on a grid of values of $a$ the absolute error of our estimator, i.e. $\left\vert
\frod(w)-f(w)
\right\vert$. 
We abbreviate it AE in the following.
Several settings are considered, each corresponding to a separate graph. In particular, we consider for each example a variety of sample sizes ($n\in\{10 \   000 ; 50 \  000 ; 100\ 000 ; 200 \  000\}$) and $X$ of different dimensions ($d_1\in 1:d_{\max}$, with $d_{\max}=6$ in Example (a) and $d_{\max}=9$ in Examples (b) and~(c)). Moreover, in each graph, we consider $3$ samples (in the graphs with different line types) and several evaluation points $\lbrace w^k\rbrace _{k=1:16}$ randomly drawn according to the joint distribution ${\rm f}_W$ (the $16$ pastel curves in the graphs). Note that, to refine our selection of $a$, we add a logarithmic grid to the standard grid of integers, after observing that the AE minimizers increase sublinearly with $d_1$.
\begin{figure}[h]
\includegraphics[width=\linewidth, trim=0cm 0.5cm 0cm 1cm,clip]{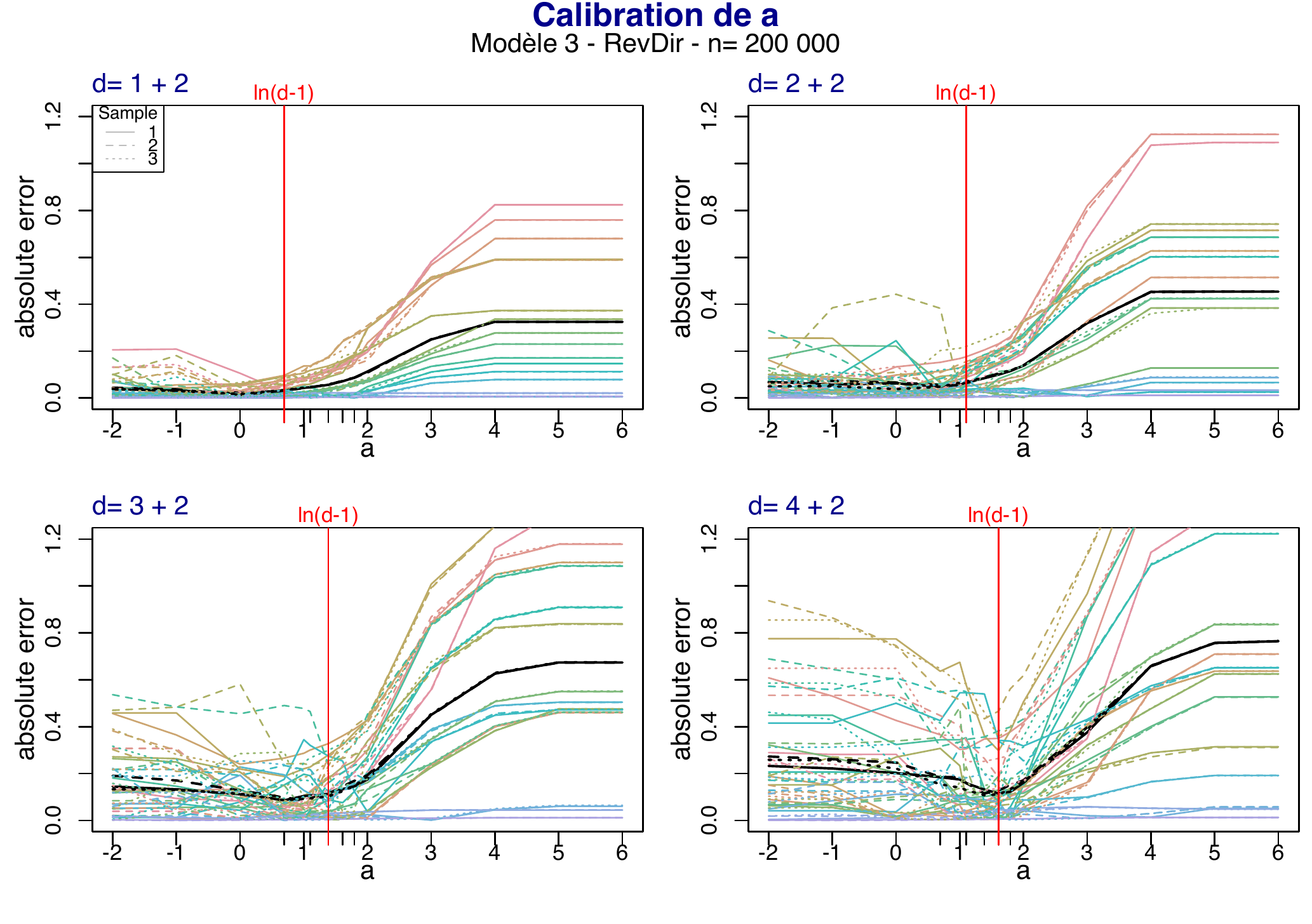} 
\vspace{-0.5cm}
\caption{\textbf{Illustration of the tuning of $a$ for Example (a) with $n=200 \  000$ for growing dimensions.} In each subgraph: AE curves in function of $a$ for $16$ evaluation points $w^k$ (the warmer the pastel color, the larger $f(w^k)$) given $B=3$ samples (differentiated by line type) at fixed dimension (specified top left in the form $d=d_1+d_2$). In black lines: the average per sample of the $16$ pastel curves. The vertical straight red line: our final choice.}
\label{fig Calib a}
\end{figure}
We simply provide here one case (Example (a) with $n=200 \   000$ in Figure~\ref{fig Calib a}), but the whole set of figures can be found in supplementary material (see \cite{JCV-supp}). %where we have considered Example (a) and the case $n=200 \   000$.

For ease of interpretation, the average per sample over the different evaluation points has been added in thicker black line. Then, our goal is to determine this minimizer as a function of the varying parameters mentioned above.
A good point is that the minimizers do not seem to depend on the sample size (cf the whole set of figures).
However the effect of the dimension is more sensitive. First note that the larger $a$, the larger the thresholds $\seuil{hj}$, thus the larger $\hat{h}$. 
We observe the chaotic behavior of \CDRodeo{} for small values of $a$ (especially for large dimension and small sample size) and, for large values of $a$, the superposition of the curves built from different samples, meaning low variance but large bias of the estimators. This corresponds to the usual phenomenon of under- and over-smoothing.
%Consequently, large values of $a$ lead to a large bias but a small variance, and vice versa.
%In particular, the behavior of \CDRodeo{} is rather chaotic for too small values (especially for large dimension and small sample size), while the behavior is almost deterministic for large values (but the error increases quickly due to the bias term).
%We first observe that the $3$ average curves are very close, highlighting the stability of our procedure. This can be observed on Figure~\ref{fig Calib a} and on the graphs of Appendix for which a large range of small values of $a$ lead to satisfying results. Note however that for too small values of $a$, we can observe a chaotic behavior of the errors in particular for small sample sizes and when the dimension is large. These values have to be avoided, in particular for samples presenting significative different absolute errors. Conversely, for large values of $a$, the stochastic error vanishes and each sample has the same behavior. However, for such values of $a$, the risk increases quickly due to the bias term.

Finally, a good trade-off is achieved by the tuning
$$a=\log(d-1),$$
and all the following simulations will be implemented with this choice.
%Considering  the whole set of simulations, the choice $a=\log(d-1)$ satisfies a good trade-off and all the following simulations will be implemented with this choice.

%\textcolor{red}{Premi\`ere partie sans doute trop n\'egative. Deuxi\`eme partie ok mais plus tard}This tuning is not optimal for Example~(c), but  still, performs better than Example~(a): our method copes rather well with the discontinuity and finds the relevant components of Example~(c), whereas the curse of dimensionality makes fully non-sparse Example~(a) hard to estimate (see also Section~\ref{sssect SimuSparse}).
%%%%
\subsubsection{... of the step size}\label{sec:beta}
\begin{figure}[t]
\includegraphics[width=\linewidth, trim=0cm 0.3cm 0cm 1.5cm,clip]{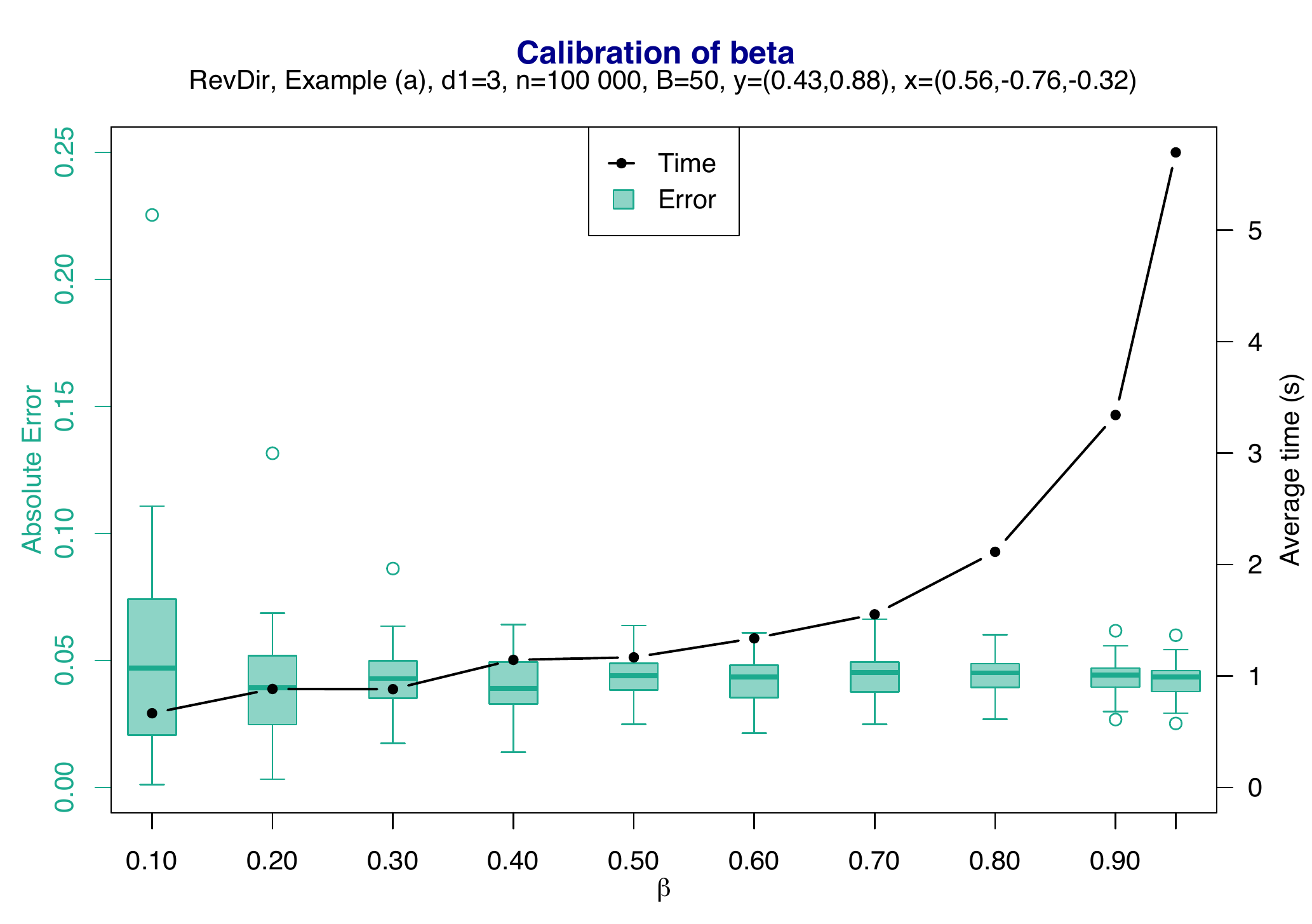} 
\vspace{-0.7cm}
\caption{\textbf{Illustration of the calibration of $\beta$.} For Example~(a) with $d_1=3$, given $B=50$ samples of size $n=100 \  000$, boxplots of the AEs and their average running times (in black) in function of $\beta$.}\label{FigSelectBeta}
\end{figure}
{
Let us now tune the step parameter, namely $\beta$ the multiplicative decay factor of the bandwidth.
As one can expect, the calibration of $\beta$ is a compromise between running times and estimation sharpness: %\textcolor{magenta}{the closer to $1$ the parameter $\beta$, the finer the associated exploration grid of bandwidths and the sharper our procedure but also the slower}
{the smaller the parameter $\beta$, the bigger the step size leading to a faster procedure but a larger approximation error}.

In Figure~\ref{FigSelectBeta}  (corresponding to Example (a) with $d_1=3$ and $n=100 \  000$), we put in perspective the boxplots (built given $50$ samples) of the AEs with their mean running times.
As one could expect with a multiplicative factor, the computational time increases exponentially fast with $\beta$: in particular, the running time explodes when $\beta\geq 0.9$.
%The comparison between Figure~\ref{FigSelectBeta} and Figure ?? in \cite{JCV-supp} shows very similar behaviors of the computational time, meaning that the dimension has a low impact. %But, as expected, we observe that absolute errors are smaller when the dimension is small. 
%Note also instability issues for too small $\beta$'s: in particular, 
Conversely the smaller $\beta$, the larger standard deviation of the boxplots, therefore $\beta$ should not be taken too small. 
%Note that at fixed value of $n$, the optimal bias does not converge to $0$ when $\beta\rightarrow 1$ to achieve the bias-variance trade-off,  but should rather be approximately equal to the variance, as it can be observed in Figure~\ref{FigSelectBeta}.

To sum up, the range of values satisfying a good compromise is quite large. To fix the parameter, we take $$\beta=0.8,$$
and all the following simulations will be implemented with this choice. 
}
%%%%
%%%%%%%%%%%%%%%%%%%%%%%%%
\subsection{Numerical performances}\label{sec:perf}
%%%
In this section, we assess the performances of our procedure according to two directions: we first visualize how our procedure reconstructs functions, then we focus on the sparsity detection, the key property of our algorithm to circumvent the curse of dimensionality.
\subsubsection{Reconstructions:  direction-by-direction visualization and estimation of $\fX$}
We first focus on a global visualization of the estimation of the function $f$: in particular,  we are interested in the performances of our estimator evaluated on a grid. 
Two kinds of estimates are considered: one in which the true $\fX$ is plugged, the other in which $\hfX$ is estimated by our procedure with the following methodology.

\paragraph{Density estimation: a RevDir \CDRodeo{} procedure for the input $\lbrace \hfX(X_i)\rbrace_{i=1:n}$.} 
{
First, for the sake of practicality, we use the same sample to compute $\hfX$ and $\frod$. Note that there is no requirement of independence in the theoretical results. 

We use the RevDir \CDRodeo{} procedure, since it can perfectly be used for estimating standard densities (cf Remark~\ref{rmkKDE}). 
Since our method is pointwise, we need to compute $\fX(X_i)$ for each $i\in1:n$.
%, which would mean a running time multiplied by $n=100\ 000$, our sample size in the following simulation. 
Note that sparsity structures are rarer in standard densities, for which all variables are of interest, than in conditional densities.
Therefore the straightforward estimation of $\fX$ is limited by the dimension of $X$ due to the curse. 
To circumvent this fact, we propose to add conditioning, artificially, by decomposing $\fX$ as follows:
$$\fX(x)={\rm f}_{X_1}(x_1) \prod\limits_{j=2}^{d_1} f_{X_j \vert X_{1:(j-1)}}(x_{1:j}).$$
Notice that the $n$ estimates $\lbrace \tilde{\rm f}_{X_1}(X_{i,1})\rbrace_{i=1}^n$ are needed as input to compute the $\lbrace \tilde{f}_{X_{2}\vert X_{1}}(X_{i,1:2})\rbrace_{i=1}^n$, which are needed to compute the $\lbrace \tilde{f}_{X_{3}\vert X_{1:2}}(x_{1:3})\rbrace_{i=1}^n$, and so on. 
%Then, we estimate each (standard or conditional) density by our procedure RevDir \CDRodeo{}.

%, even if the assumption in Proposition~\ref{propfXtilde} is not satisfied. However, we have removed $X_i$ from the sample to estimate $f_{X_j\vert X_{1:(j-1)}}(X_{i,1:j})$ as in leave-one-out procedures (even if for dimension $j>3$, due to the input $\lbrace \tilde{f}_{X_{j-1}\vert X_{1:(j-2)}}(X_{k,1:(j-1)})\rbrace_{k=1}^n$, some dependency on $X_i$ will remain). But, notice that in our theoretical results, we do not need any independence assumption between $X$ and~$\tilde{X}$.

Observe also that the previous calibration of $a$, namely $a=\log(d-1)$, does not extend for univariate densities. Based on preliminary numerical experiments, we set $a=-1$ for the univariate case.

Implemented in \textsc{R}, with a $3.1$ GHz Intel Core i7 processor, the running times for $\hfX$ in Example~(a) in dimension  $d_1=2$ and  in Examples~(b) and ~(c) in dimension  $d_1=3$ is summarized in the following table:\\

\indent\indent \begin{tabular}{|c||c|c|c||c|}
\hline 
 & \multicolumn{3}{c||}{Mean time per run (seconds)} & Total time for $100\ 000$ runs \\ 
\hline 
 \rule[-1ex]{0pt}{2.5ex} & $\tilde{\mathrm{f}}_{X_1}$ & $\tilde{f}_{X_2\vert X_1}$  &  $\tilde{f}_{X_3\vert X_{1:2}}$ & $\hfX$ \\ 
\hline 
Model (a) & 0.734 & 0.654 & N.A. & 138 780s (around 1d 15h)\\ 
\hline 
Model (b) & 1.31 & 1.61 & 1.72 & 463 559s (around 5d 9h) \\ 
\hline 
Model (c) & 0.675 & 1.17 & 1.05 & 289 695s (around 3d 8h)\\ 
\hline 
\end{tabular}\\

The running times strongly depend on the distance between the initialization bandwidth and the selected one, which explains non increasing running times when the dimension grows for Models~(a) and~(c).

One may object that several days of computation for the preliminary estimator is quite long.  But, note that it is done without parallelization. Given a powerful enough cluster,  the running time can be divided by $n$ using parallelization over the evaluation points.
%Fortunately, the execution for each evaluation point can be performed in parallel (given a sufficiently powerful computing cluster).
%Therefore, the executions can be performed in parallel only dimension after dimension, multiplying the running time by a factor $d_1$. 
}

\begin{figure}[h]
\includegraphics[width=\linewidth, trim=0cm 0cm 0cm 1cm,clip]{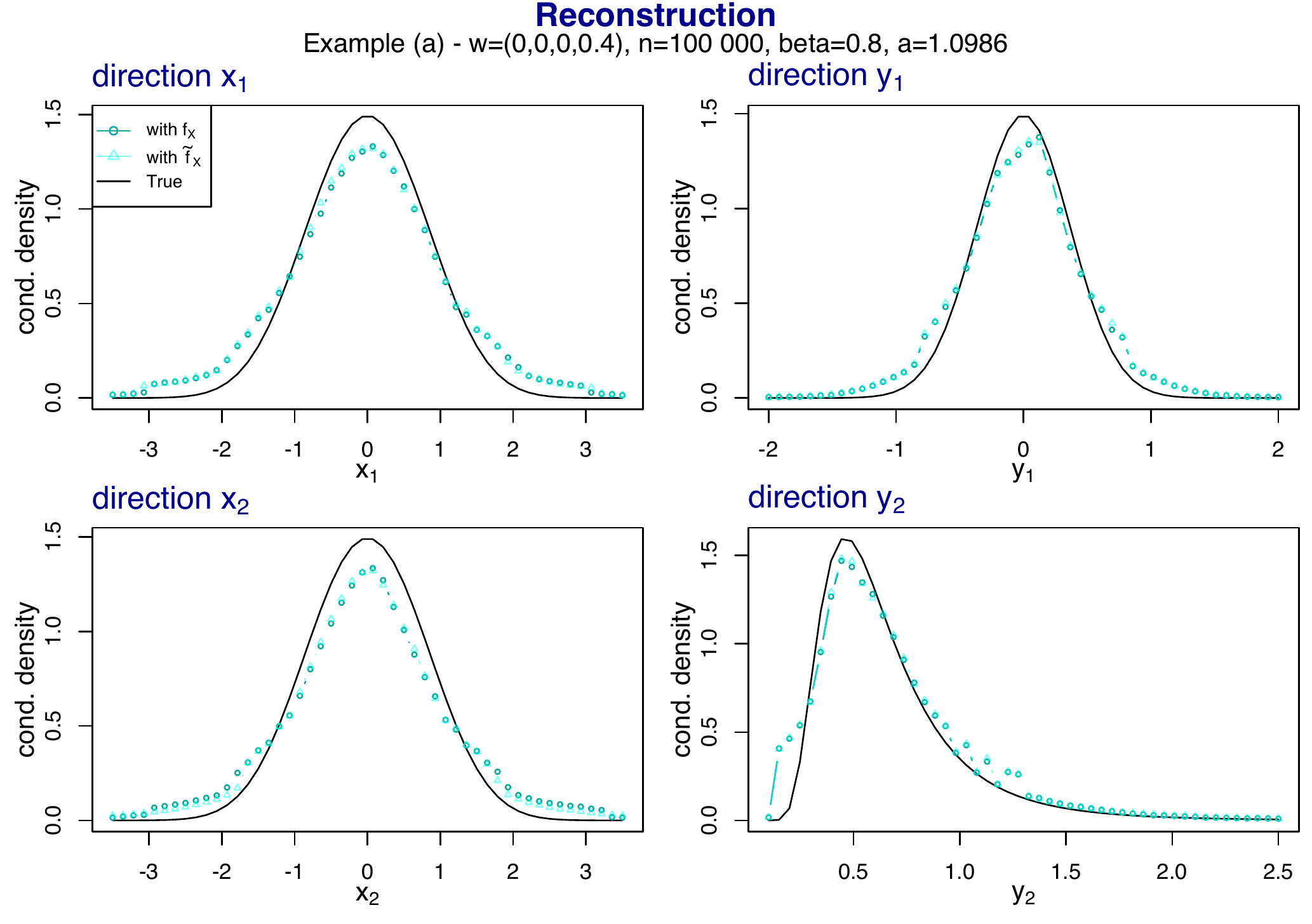}\\
\vspace{-1cm}
\caption{\textbf{Reconstruction for Example~(a).} Compared to the true conditional density (in black line), our pointwise estimates with the true marginal density $\fX$ (circles in darker shade) or with the estimate $\hfX$ (triangles in lighter shade). Only one direction is varying specified at the top left of each graph, the others being fixed to a set point: $w=(0,0,0,0.4)$. The sample size is $n=100\ 000$.}\label{FigReconst(a)} 
\end{figure}

\begin{figure}[p]
\includegraphics[width=\linewidth, trim=0cm 0cm 0cm 1cm,clip]{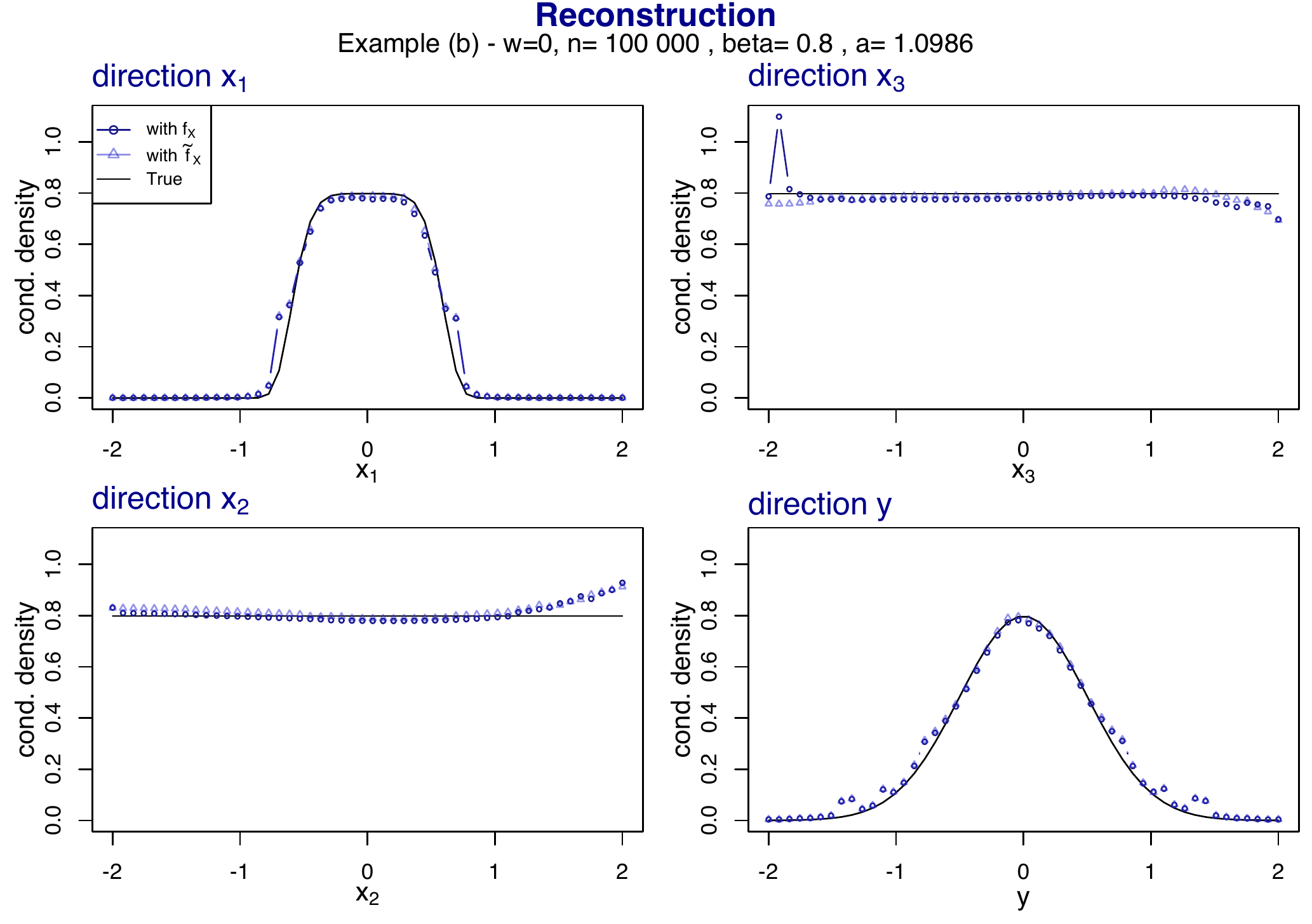} 
\vspace{-1cm}
\caption{\textbf{Reconstruction for Example~(b).} See the description in Figure~\ref{FigReconst(a)}, except: $w=(0,0,0,0)$.}\label{FigReconst(b)}
\vspace*{0.5cm}
\includegraphics[width=\linewidth, trim=0cm 0cm 0cm 1cm,clip]{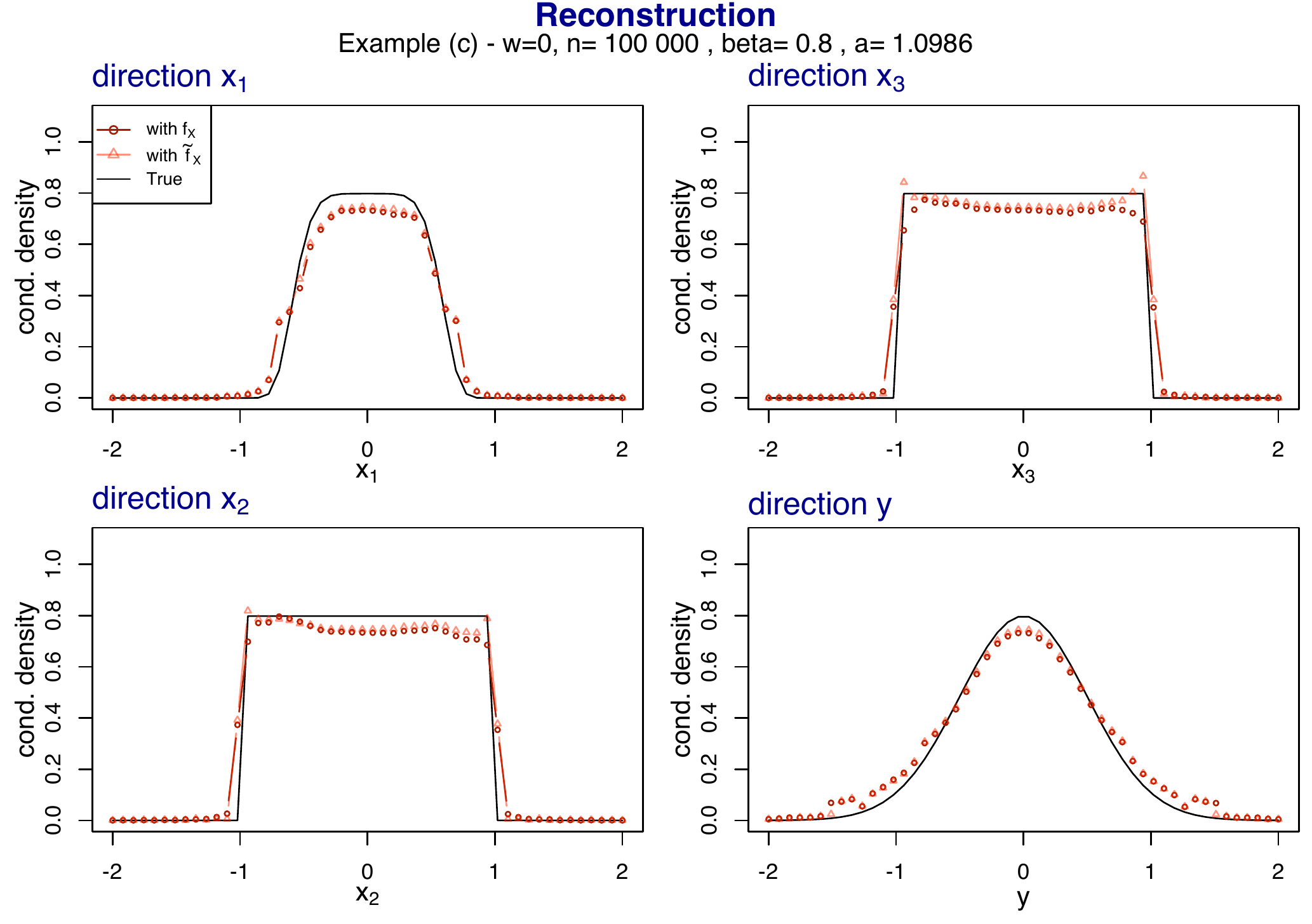} \vspace{-1cm}
\caption{\textbf{Reconstruction for Example~(c).} See the description in Figure~\ref{FigReconst(a)}, except: $w=(0,0,0,0)$.}\label{FigReconst(c)}
\end{figure}

\paragraph{Visualization.}  In Figures~\ref{FigReconst(a)}, \ref{FigReconst(b)}, and \ref{FigReconst(c)}, 
%we compare, for respectively Examples~(a), (b) and (c), 
the two kinds of estimates 
%with $\fX$ or withthe estimate $\hfX$ following the above methodology   
are built from a sample of dimension $d=4$ and size $n=100\ 000$
for respectively 
%simulated according to 
Examples~(a), (b) and (c). 
Limited to two-dimensional visualizations, we vary only one component at a time, the others being fixed to a set point: $w=(0,0,0,0.4)$ for Example~(a) and $w=(0,0,0,0)$ for Examples~(b) and (c).


The overall signal is nicely recovered.
Comparing the different examples,  Example~(a) is the least accurately estimated: the estimates are oversmoothed near the modes. It was expected since it is the example without sparsity and even in dimension as small as $4$, the curse deteriorates the convergence rate.

Thanks to the strong similarity between Example~(b) and Example~(c), the impact of the discontinuity can be properly visualized: (b) is clearly more accurately estimated than (c), even though the focus point $w=0$ is not really close to the discontinuity points $\pm 1$ (in the directions $x_j$). The loss of accuracy is once again due to the curse, as the directions $x_j$, $j\geq2$, in Example~(c) are not completely irrelevant.  The \CDRodeo{} procedure does not consider the relevance of a variable as a binary answer: {in fact, when a variable is relevant, it can be more or less relevant.} See the analysis of the selected bandwidths in the next section for more details.

Besides, in all examples, the estimation is less accurate at the specific points where Assumption \Monoton{} is not satisfied. %Note that every (standard or conditional) density presents at least one point where Assumption \Monoton{} is not satisfied. 
Taking account that the chosen kernel is Gaussian, thus of order $2$, it especially happens around the zeros of the second derivative. % (which induces a change of dominant term in $\E{[}Z_{hj}{]}$ when $h\rightarrow 0$, leading to non monotony). \red{\\ Revoir en fonction de la section 2}. 

 Then, note that RevDir \CDRodeo{} may stop during the direct Step (as $\vert Z_{hj}\vert$ may become smaller than $\seuil{hj}$) but has no impact on the increasing step (Reverse Step). That is why the initialization $h_0$ is set as the lower bound of its range (see Equation~\eqref{ineq-initialisation}) to minimize the undesirable impacts. For illustration of the improvement made by the RevDir algorithm, see Figures IV.3, IV.4 and IV.5  in Ph.D. thesis \citep{Jeanne-these} which compares the Direct and the RevDir procedures.

Note lastly that the estimates with either $\fX$ or $\hfX$ are very close to each other. More precisely, the estimates with $\hfX$ is slightly better (in particular, near the modes and near the discontinuity in Example~(c)): \cite{DP16} actually prove that dividing by an estimator of the density produces better results than if the density itself was used.
That is the reason why the reliability of our results is maintained in the following part even if the true $\fX$ is used in order to save the running times of computing the $\hfX(X_i)$'s for several samples and dimensions.

%%%%%%%%%%%%%%%%%
\subsubsection{Impact of the dimension and sparsity detection\label{sssect SimuSparse}}

\begin{figure}[t]
%\includegraphics[width=\linewidth, trim=0cm 0.5cm 0cm 1cm,clip]{plotAdapt/SparseRobust_Model1_RD_w0_est.pdf} 
%\vspace{-1.1cm}
%\includegraphics[width=\linewidth, trim=0cm 0cm 0cm 1.5cm,clip]{plotAdapt/SparseRobust_Example(a)_RelErr_y2_04.pdf} 
%\vspace{-1cm}
%\caption{\label{FigNonSparseEst}\textcolor{forest}{\textbf{Estimation with increasing dimension.} Boxplots of the estimates $\frod(w)$ in function of $d$, given $50$ samples of size $n=100 \  000$ for Example~(a). The dashed horizontal segment is the true value $f(w)$ (evalated for $x=0$ and $y=(0,0.5)$).}}
\includegraphics[width=\linewidth, trim=0cm 0cm 0cm 1.5cm,clip]{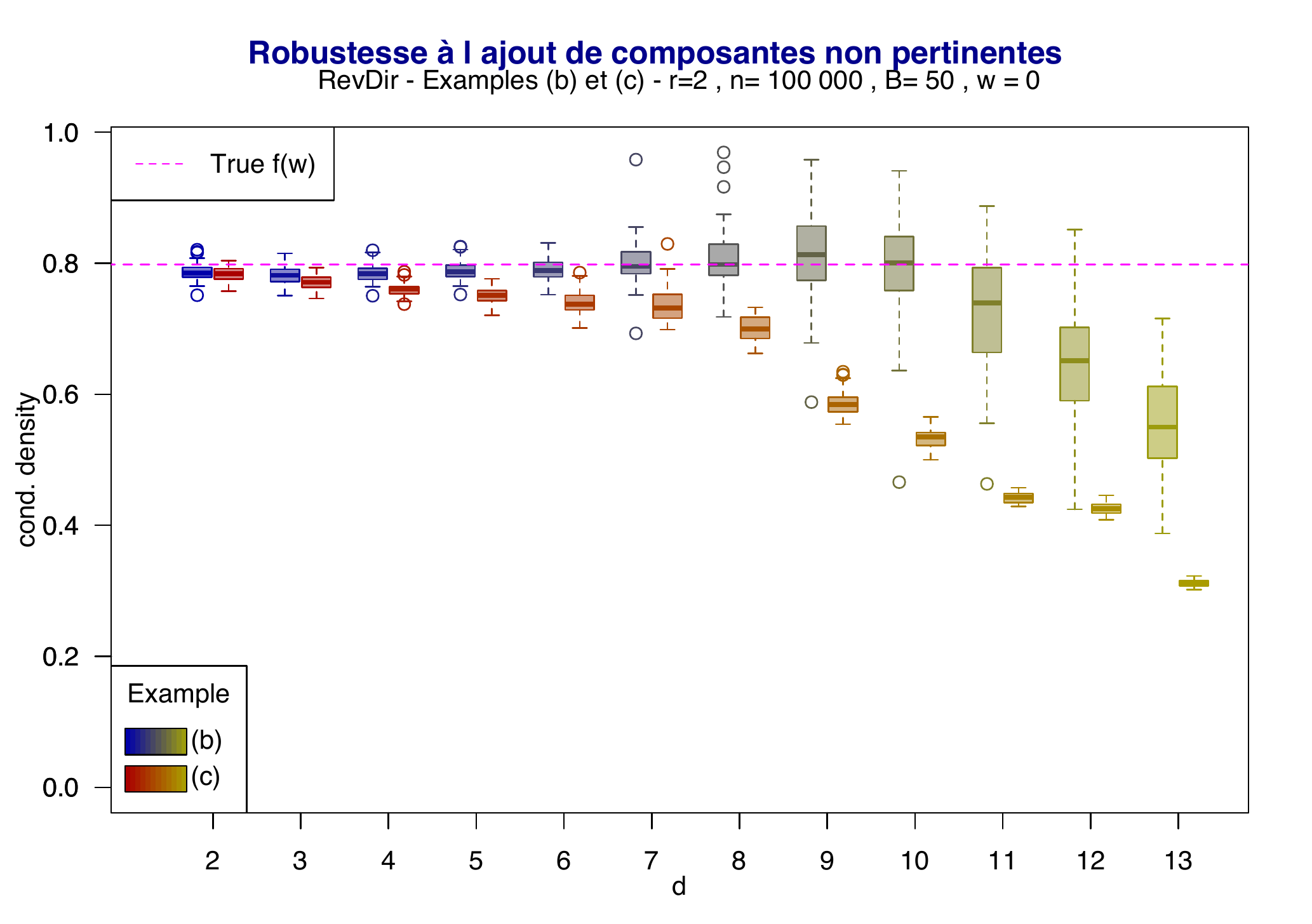}
\vspace{-1cm}
\caption{\label{FigSparseRobustEst}\textbf{Robustness to addition of irrelevant variables.} Boxplots of the estimates $\frod(w)$ in function of $d$, given $50$ samples of size $n=100 \  000$ for Example~(b) (in red shades) and Example~(c) (in green shades). The dashed horizontal line is the true value $f(w)$ (at the evaluation point $w=0$).
}
\end{figure}
\begin{figure}[h!]
\includegraphics[width=\linewidth, trim=0cm 0cm 0cm 1.5cm,clip]{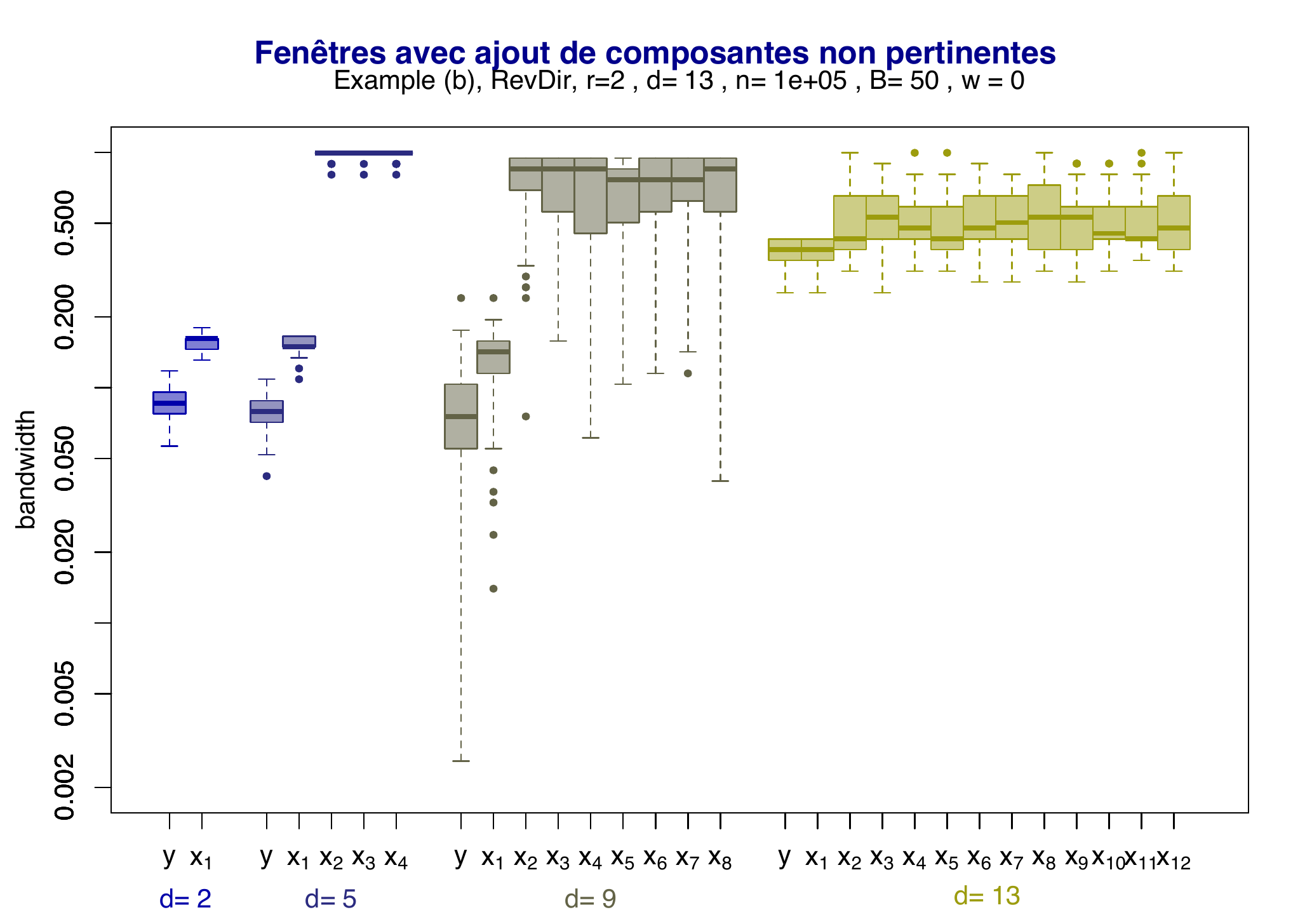} \\
\includegraphics[width=\linewidth, trim=0cm 0cm 0cm 1.5cm,clip]{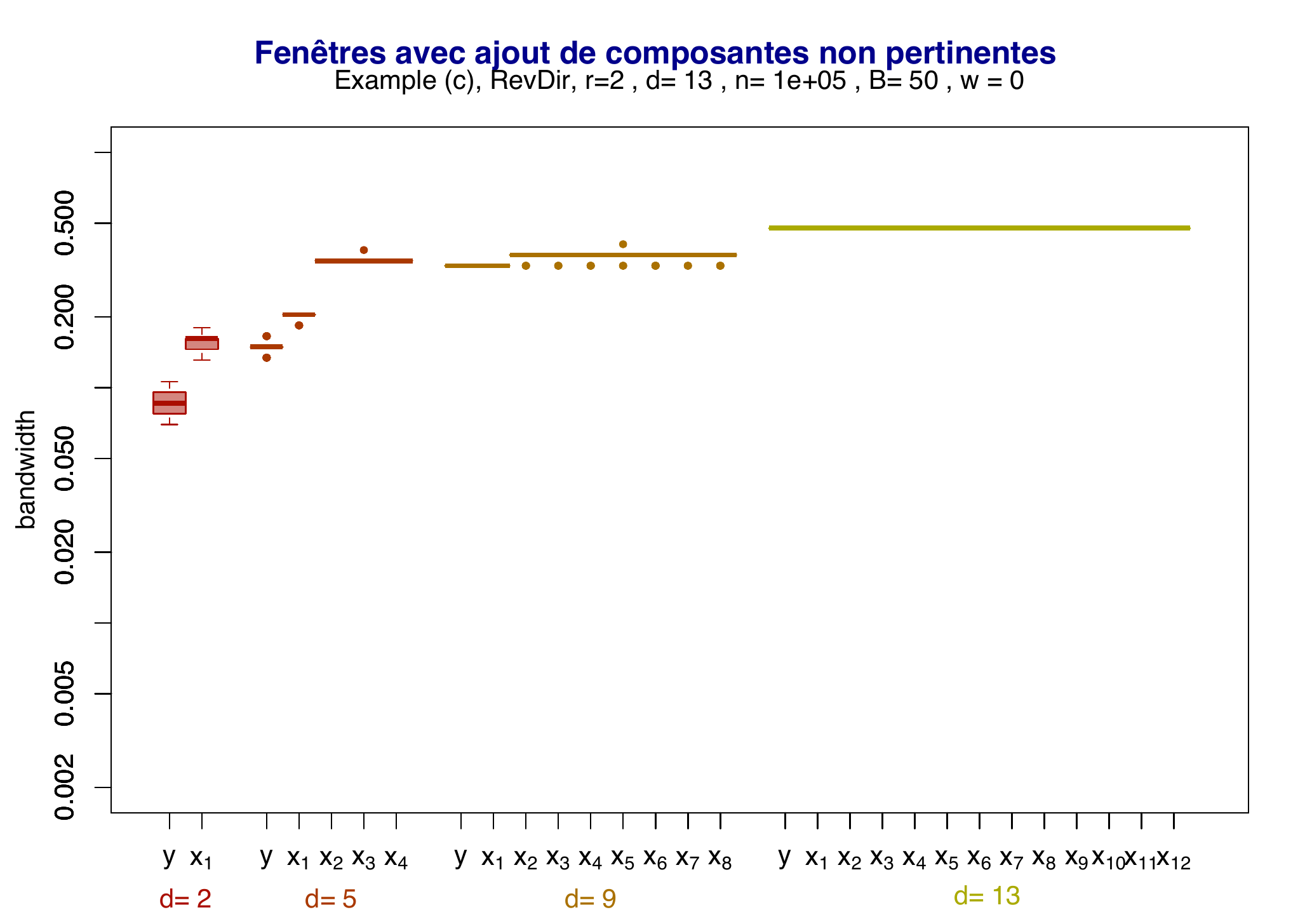} 
\vspace{-1cm}
\caption{Bandwidths associated to the estimates of Figure~\ref{FigSparseRobustEst} for the dimensions $d\in\{2,5, 9,13\}$. Top: Example~(b). Bottom: Example~(c).}
\label{FigSparseRobustBw}
\end{figure}
 Let us now consider how the RevDir \CDRodeo{}  procedure detects the sparsity structure. For examples with sparsity structure -- namely Examples (b) and (c) --, we check the robustness {to irrelevant explanatory variables: starting with the fully relevant example at dimension $d=2$, we gradually add irrelevant variables until dimension $d=13$. }
 %
 % we provides the boxplots of the \CDRodeo{} estimates for increasing dimensions $d_1$ (from $1$ to $7$). Observe as each time a variable is added, the bias increases until the signal is completely lost. 
 In Figure~\ref{FigSparseRobustEst}, the boxplots are built from 
 $50$ simulated samples of size $n=100\ 000$ with varying dimension $d_1$ from $1$ to $12$: in bluish shades, the estimates of Example~(b) and in reddish shades, the ones of Example~(c). We also provide in Figure~\ref{FigSparseRobustBw} the boxplots of the selected bandwidths for the dimensions $d\in\{2,5,9,13\}$. Notice that our fully nonparametric procedure actually ends within reasonable times for dimensions as large as $13$ (e.g. $40$ minutes for the whole $600$ estimates of Figure~\ref{FigSparseRobustEst} on samples of size $100\ 000$), while most nonparametric methods struggle to handle dataset of dimension higher than $4$.

Usually, without sparsity, each added variable worsens the estimation: see for instance  Example~(a) with increasing (relevant) dimension in Figure~\ref{FigNonSparseEst} of the supplementary material \cite{JCV-supp}, in which our method struggles providing good estimates as soon as the dimension $5$. For Example~(b) (where the relevant dimension is $r=2$), until the dimension $6$, our method has the same behavior as for dimension $2$. For larger dimensions, the estimation is progressively noised by the too many irrelevant variables, to finally lose the signal beyond the dimension $10$. 
The bandwidths in Figure~\ref{FigSparseRobustBw} give a good understanding of how the procedure handles the extra variables. 
Comparing the dimensions $d=2$ and $d=5$, the relevant bandwidth components (namely the directions $y$ and $x_1$) are selected at very similar values (called hereafter their "expected values"), while the irrelevant components (in dimension $d=5$) are taken as high as possible, around the value $1$ (the upper limit of the bandwidth grid): thus, the bias-variance trade-off is unchanged, ensuring a quality of estimation as good as in dimension $d=2$. 
In dimension $d=9$, the larger dimension makes the detection of irrelevant variables more difficult, producing variance in the bandwidth selection. 
Nevertheless, the relevant components are still selected at their expected value (but with more variance), producing rather good estimates. 
In dimension $d=13$, the sparsity is less accurately detected: the irrelevant bandwidths decrease to $0.5$. Their product $0.5^{11}$ reaches numerically the emergency stop $\approx \frac1n$. Therefore, there is not enough room left for the relevant components to decrease until their expected value, which explains the loss of signal observed in Figure~\ref{FigSparseRobustEst}. Note that in this last setting $d>\log n$, and that is the reason why the emergency stop is reached. More generally, this framework seemed to be out of reach for \Rodeo{}-type procedures: in particular, in \citep{LW08} where growing dimensions with $n$ are considered, the framework is also restricted to dimensions $d\ll \log n$. 
%In our simulations, $\log(n)\approx 11.5$ (for $n=100\ 000$), which explains the deterioration of performances beyond the dimension $10$.

Let now consider Example~(c). 
The same phenomenon occurs, but complicated by the discontinuity of $f$ in the directions $x_j$: away from $\pm1$, the relevant dimension is $r=2$, but in the neighborhood of $\pm1$, these components are highly relevant. 
In fact, these neighborhoods depend on the bandwidth: the larger the bandwidth, the larger the support of $K_h$ until reaching the points $\pm1$, and 
once $\pm1$ belongs to  the support of $K_h$, the components $x_j, j>1$, are detected as relevant.
% for bandwidth values smaller than 
This is the reason why these bandwidth components are much smaller in Figure~\ref{FigSparseRobustBw} (bottom) (around the value $0.48$ instead of $1$ in Example~(c)). 
These smaller components amplify the phenomenon described for Example~(b): as soon as dimension $d=5$, the relevant components can no longer decrease to their expected value; in dimension $d=9$, there is almost no room left for the relevant components, and in dimension $d=13$, the relevant components are completely lost.

{
All in all, the overall behavior of our procedure is very satisfying: the RevDir \CDRodeo{} procedure nicely detects relevant variables and is robust to extra irrelevant in moderate dimensions ($d\leq\log n$).
The difficulties described in the last paragraphs are inherent to the curse of dimensionality and is bound to occur with any nonparametric procedure. }
%%%%%%%%%%%%%%%%%%%%%%%%%%
%\subsection{Ce qui suit est du vieux}
%\textcolor{blue}{Commentaire sur l'estimation de $\fX$ a cout faible}
%\textcolor{blue}{We show that Assumption~\Monoton{} is satisfied for all studied functions in our numerical setting. :-)}
%We study the estimation procedure by generating samples $(X_i,Y_i)$ following 2 models:
%\begin{itemize}
%\item $X_i\sim\mathcal{U}_{[0,1]}$ and $f(x,y)=...,$
%$w=$
%\item
%\end{itemize}
%We use $\beta=$, $a=$, 
%%%%%%%%%%%%%%%%%%%%%%%%%%
%%%%%%%%%%%%%%%%%%%%%%%%%%
\input{preuvesrodeo.tex}\label{sec:preuves}
%Preuves 
%%%%%%%%%%%%%%%%%%%%%%%%%%
%%%%%%%%%%%%%%%%%%%%%%%%%%
\input{appendice.tex}
% Appendix

\noindent {\bf Acknowledgements:} We are very grateful to Benjamin Auder  (Universit\'e Paris-Saclay) who helped us  for parallelization of Rodeo algorithms.
%%%%%%%%%%%%%%%%%%%%%%%%%%
%%%%%%%%%%%%%%%%%%%%%%%%%%
\bibliographystyle{apalike}
\bibliography{biblio}
\end{document}

\title{Supplementary material for Adaptive greedy algorithm for moderately large dimensions in kernel conditional density estimation}

\author{
{\sc Minh-Lien Jeanne Nguyen},\\[2pt]
 Mathematical Institute,\\ University of Leiden\\
Niels Bohrweg 1, 2333 CA Leiden, Netherlands\\
m.j.nguyen@math.leidenuniv.nl\\[7pt]
{\sc Claire Lacour},\\[2pt]
 LAMA, CNRS\\ Univ Gustave Eiffel,
        Univ Paris Est Creteil\\ 
        F-77447 Marne-la-Vall\'ee, France\\
        claire.lacour@univ-eiffel.fr \\[7pt]
{\sc Vincent Rivoirard} \\[2pt]
CEREMADE, CNRS, UMR 7534\\
       Universit\'e Paris-Dauphine, PSL University\\
       75016 Paris, France\\
       Vincent.Rivoirard@dauphine.fr\\[7pt]
}

\maketitle

\begin{abstract}
{
This supplementary material of \cite{JCV} provides both theoretical and numerical supplementary results.
First,  we relax the assumption~\Monoton{} of bias convexity in an assumption of bias monotony and provide a new theorem specifying the rate of convergence under this new assumption.
The risk is now bounded by $\left(\frac{(\log n)^{(2+a)}}{n}\right)^{\frac{s}{2s+r}}$,  meaning that the price to pay is an extra logarithmic term $(\log n)^{\frac{2s}{2s+r}}$.
In the second part,  we first provide several graphs to complete the calibration of the parameter $a$, then
we give an illustration of the curse of dimensionality without sparsity structure by studying Example (a) in growing dimensions.}
\end{abstract}
%%%%%%%%%%%%%%%%%%%%%%%%%%
%%%%%%%%%%%%%%%%%%%%%%%%%%%
\section{Minimax rates under relaxation of Assumption~\Monoton{}}
Subsequent Theorem~\ref{thm-bis}
is a variation of Theorem~\ref{thm} of \cite{JCV}. Except for Theorem~\ref{thm-bis}, all number references are from the paper \cite{JCV}. We relax the assumption \Monoton{} of bias convexity by an assumption of \emph{bias monotony}:
\begin{assumption}{\Monotonbis\label{Monotonbis}}
{\\
For all $j\in\rel$, for all $h$ and $h'\in(\mathbb{R}^*_+)^d$ such that $h\preceq h'$, 
$h_j\vert\mathbb{E}{[}\bar Z_{h,j}{]}\vert 
\leq h'_j\vert \mathbb{E}{[}\bar Z_{h',j}{]} \vert $,
where $\bar Z_{h,j}$ is defined as $Z_{h,j}$ in \eqref{formulaZhj} but with true $\fX$ replacing $\hfX$.
}
\end{assumption}
The previous assumption consists in replacing the condition $\vert\mathbb{E}{[}\bar Z_{h,j}{]}\vert 
\leq \vert \mathbb{E}{[}\bar Z_{h',j}{]} \vert $ of Assumption~\Monoton{}  by the weaker condition $h_j\vert\mathbb{E}{[}\bar Z_{h,j}{]}\vert 
\leq h'_j\vert \mathbb{E}{[}\bar Z_{h',j}{]} \vert $, which allows more room in the monotony assumption, in particular when components of $h$ are much smaller than components of $h'$. {
Notice that $h_j \mathbb{E}{[}\bar Z_{h,j}{]}$ is proportional to the difference of bias in the direction $j$ (cf Equation\eqref{diffBiais}).} The new result is as follows:
\begin{theorem}\label{thm-bis} For any $r\in 0:d$, $1<s\leq p$ and $L>0$, if $f$ has only $r$ relevant components  and belongs to $\mathcal{H}_d(s, L)$, then
under Assumptions \AfXmin, \AestfX and \Monotonbis, the pointwise risk of the RevDir \CDRodeo{}  estimator $\frod(w)$ is bounded as follows: for any $l\geq 1$, for $n$ large enough, 
\begin{equation}\label{eq:ratebis}
\E\left[
\left\vert
\frod(w)-f(w)
\right\vert^l
\right]^{1/l}
\leq  \text{C} 
 	\left(\frac{(\log n)^{(2+a)}}{n}\right)^{\frac{s}{2s+r}} 
\end{equation}
where $\text{C}$ only depends on $d,r,K, \beta,\delta, L, s,\|f\|_{\infty,\UUn}$. 
%{\color{violet} Sûre que ça ne dépend pas de $l$}
\end{theorem}
\begin{remark}
%The logarithmic term of the upper bound of Inequality~\eqref{eq:ratebis} is $(\log n)^{\frac{s(2+a)}{2s+r}},$ instead of $(\log n)^{\frac{sa}{2s+r}}$ obtained in Theorem~\ref{thm} of \cite{JCV}.
Compared to the logarithmic term $(\log n)^{\frac{sa}{2s+r}}$ in the original result (Theorem~\ref{thm} of \citep{JCV}), the price to pay is the extra factor $(\log n)^{\frac{2s}{2s+r}}$.
\end{remark}
\begin{proof}[Proof of Theorem~\ref{thm-bis}] In the following, the notations of \Cref{ssectNotations} and \ref{sec:eMajBias} of \citep{JCV} are re-used.

Analyzing the proof of Theorem~\ref{thm} of  \citep{JCV}, observe that Assumption \Monoton{} is only used to obtain the bound \eqref{eMajBias} of the bias term $\left\vert\biasbar{h}\right\vert$ on $\lbrace\hat{h}=h\rbrace\cap\Ehp$ for all $h\in\Hhp$. So, to obtain Theorem~\ref{thm-bis}, we only have to bound $\left\vert\biasbar{h}\right\vert$ with the relaxed assumption, then adjust the new bias-variance tradeoff.
In the sequel, we just specify the modifications from the original proof.

 Following the proof of the bound \eqref{eMajBias} (see \Cref{sec:eMajBias}), the first change occurs in 
Subcases (C.a) and (C.b).
%To bound $\left\vert\biasbar{h}\right\vert$, we follow the steps of the proof of Inequality (5.6) (see Section 5.5.1 of \citep{JCV}). Bounds obtained for (Case A) and (Case B) remain unchanged. 
%We fix $j\in \rel$ such that $t_j\leq t_{j_0}$. 
%For (Case C), we follow the proof until the introduction of the subcases.  I
\begin{itemize}
\item[Subcase (C.a):] Since $\hu{j_0}\preccurlyeq h^{(t_j)}$, we apply Assumption~\Monotonbis{} instead of using Assumption~\Monoton{}:
$$\hu{j_0}_j\left\vert \E\left[\bar{Z}_{\hu{j_0},j}\right]\right\vert 
\leq h^{(t_j)}_j\left\vert \E\left[\bar{Z}_{h^{(t_j)},j}\right]\right\vert.$$
\item[Subcase (C.b):] Since $\huvar{j_0}\preccurlyeq h^{(t_j)}$ and $j\in\rel$, we apply Assumption~\Monotonbis{}:
$$\hu{j_0}_j\left\vert \E\left[\bar{Z}_{\hu{j_0},j}\right]\right\vert 
=\huvar{j_0}_j\left\vert \E\left[\bar{Z}_{\huvar{j_0},j}\right]\right\vert
\leq h^{(t_j)}_j\left\vert \E\left[\bar{Z}_{h^{(t_j)},j}\right]\right\vert.$$
\end{itemize}
Then, in the right hand side of \eqref{poursup1}, we artificially factorize (and divide)
by $\hu{j_0}_j$, then apply the previous inequality to obtain instead of Inequality~\eqref{poursup2}:
\begin{align}
\left\vert \biasbar{\hint{t_{j_0}}} -\biasbar{\hint{t_{j_0-1}}}\right\vert
&\leq \sum_{\substack{j\in\rel\\t_j\leq t_{j_0}}} \int_{u=0}^1  \hu{j_0}_j \left\vert\E\left[\bar{Z}_{\hu{j_0},j}\right]\right\vert
\times \left( \hint{t_{j_0}}_j-\hint{t_{j_0-1}}_j\right)\frac{du}{\hu{j_0}_j}\nonumber\\
&\leq \sum_{\substack{j\in\rel\\t_j\leq t_{j_0}}}  h^{(t_j)}_j \left\vert \E\left[\bar{Z}_{h^{(t_j)},j}\right]\right\vert \times \left( \hint{t_{j_0}}_j-\hint{t_{j_0-1}}_j\right)\int_{u=0}^1 \frac{du}{\hu{j_0}_j}\nonumber\\
&= \sum_{\substack{j\in\rel\\t_j\leq t_{j_0}}} h^{(t_j)}_j\left\vert \E\left[\bar{Z}_{h^{(t_j)},j}\right]\right\vert\times\left(\log\Big(\hint{t_{j_0}}_j\Big)-\log\Big(\hint{t_{j_0-1}}_j\Big)\right).\label{diffBiais}
\end{align}
Then remember we have to replace $\hint{t_{\jm-1}}$ by $\hint{\tminimax}$ given our slight abuse of notation, so:
\begin{align*}
\sum\limits_{j_0=\jm}^r \left\vert \biasbar{\hint{t_{j_0}}} -\biasbar{\hint{t_{j_0-1}}} \right\vert
&\leq \sum\limits_{j_0=\jm}^r\sum_{\substack{j\in\rel\\t_j\leq t_{j_0}}} h^{(t_j)}_j\left\vert \E\left[\bar{Z}_{h^{(t_j)},j}\right]\right\vert\times\left(\log\Big(\hint{t_{j_0}}_j\Big)-\log\Big(\hint{t_{j_0-1}}_j\Big)\right)\\
&\leq\sum\limits_{j=\jm}^r h^{(t_j)}_j\left\vert \E\left[\bar{Z}_{h^{(t_j)},j}\right]\right\vert\sum\limits_{j_0=\jm}^j\left(\log\Big(\hint{t_{j_0}}_j\Big)-\log\Big(\hint{t_{j_0-1}}_j\Big)\right)\\
&\leq\sum\limits_{j=\jm}^r h^{(t_j)}_j\left\vert \E\left[\bar{Z}_{h^{(t_j)},j}\right]\right\vert\left(\log\Big(h^{(t_j)}_j\Big)-\log\Big(\hint{\tminimax}_j\Big)\right)\\
&\leq \sum\limits_{j=\jm}^r h^{(t_j)}_j\left\vert \E\left[\bar{Z}_{h^{(t_j)},j}\right]\right\vert \Big(-\log\hint{\tminimax}_j\Big).
\end{align*}
From line 2 to line 3, the sum is telescoping and notice that $\hint{t_{j}}_j=h^{(t_j)}_j$. 
For the last line, note that: 
$$\CAm (\log n)^{\Am} n^{-\frac1{2s+r}}\leq \hint{\tminimax}_j \leq h^{(t_j)}_j \leq 1.$$
Indeed, for the first and second inequalities, notice that $t_j<\tminimax$ (Case C); for the last inequality, the procedure only explores bandwidth components no larger than $1$.
So, for $n$ large enough:
$$0\leq \log\Big(h^{(t_j)}_j\Big)-\log\Big(\hint{\tminimax}_j\Big)\leq 0-\log\Big(\CAm (\log n)^{\Am} n^{-\frac1{2s+r}}\Big)\leq \log n.$$
Therefore, the control \eqref{biaispoursupp} of $\left\vert\biasbar{h}\right\vert$ becomes:
$$\left\vert\biasbar{h}\right\vert\leq r\CBbar  {\CAm}^s \left(\log n\right)^{\Am s} n^{-\frac{s}{2s+r}}
+ (\log n)\sum\limits_{j=\jm}^r  \left\vert \E\left[\bar{Z}_{h^{(t_j)},j}\right]\right\vert 
 h^{(t_j)}_j.$$
%since $\hint{\tminimax}_j\geq n^{-1}$ by using the emergency stop and Condition \eqref{ineq-initialisation}. 
Using final arguments of 
Section~\ref{sec:eMajBias}, we obtain
 \begin{align*}
\1_{\Ehp\cap \lbrace\hat{h}=h\rbrace} \left\vert \biasbar{h} \right\vert
&\leq  r \CBbar{\CAm}^s \left(\log n\right)^{\Am s} n^{-\frac{s}{2s+r}}
\\
&\hspace{2cm}+ r\log n\times \max\left( \tfrac{7\Cl }{4\beta^{\frac{d-r}2}{\CAm}^{\frac{r}2}} \frac{(\log n)^{\frac{a-\Am r}2}}{n^{\frac{s}{2s+r}}}
, {\rcp} \left(\frac{(\log n )^{a}}{n}\right)^{\frac{p}{2p+1}}\right).
\end{align*}
To obtain the final result, we now follow the proof of Theorem~\ref{thm} until \Cref{eqPourSuppThm}, then we apply our new bound of $\vert\biasbar{h}\vert$, which adds a factor $\log n$ in the last term. 
Thus, we have to modify the optimization of $\Am$: 
the optimal value is now $\Am=\frac{2+a}{2s+r}$. This leads to the new logarithmic exponent $\frac{s(2+a)}{2s+r}$.
\end{proof}

\section{Numerical results}
In this section, we extend the numerical results of \citep{JCV} by providing supplementary graphs. In Section~\ref{sec:a:supp}, we complete the calibration of the parameter $a$ (Section~\ref{sec:a} of \citep{JCV}) for different sample sizes. In Section~\ref{sec:curse:supp}, we provide an illustration of the curse of dimensionality without sparsity structure.
%%%%%%%%%%%%%%%%%%%%%%%%%%
\subsection{Graphs for the tuning of $a$.}\label{sec:a:supp}

This section contains the supplementary graphs of the calibration of the parameter $a$. Like in Figure~\ref{fig Calib a} of \cite{JCV}, we display  for Examples (a), (b) and (c), the absolute error of our estimates (abbreviated AE) in function of $a$
 for the extra sample sizes $n\in\{10 \   000 ; 50 \  000 ; 200 \  000\}$.

\begin{figure}[p]
\includegraphics[width=\linewidth, trim=0cm 0.5cm 0cm 1cm,clip]{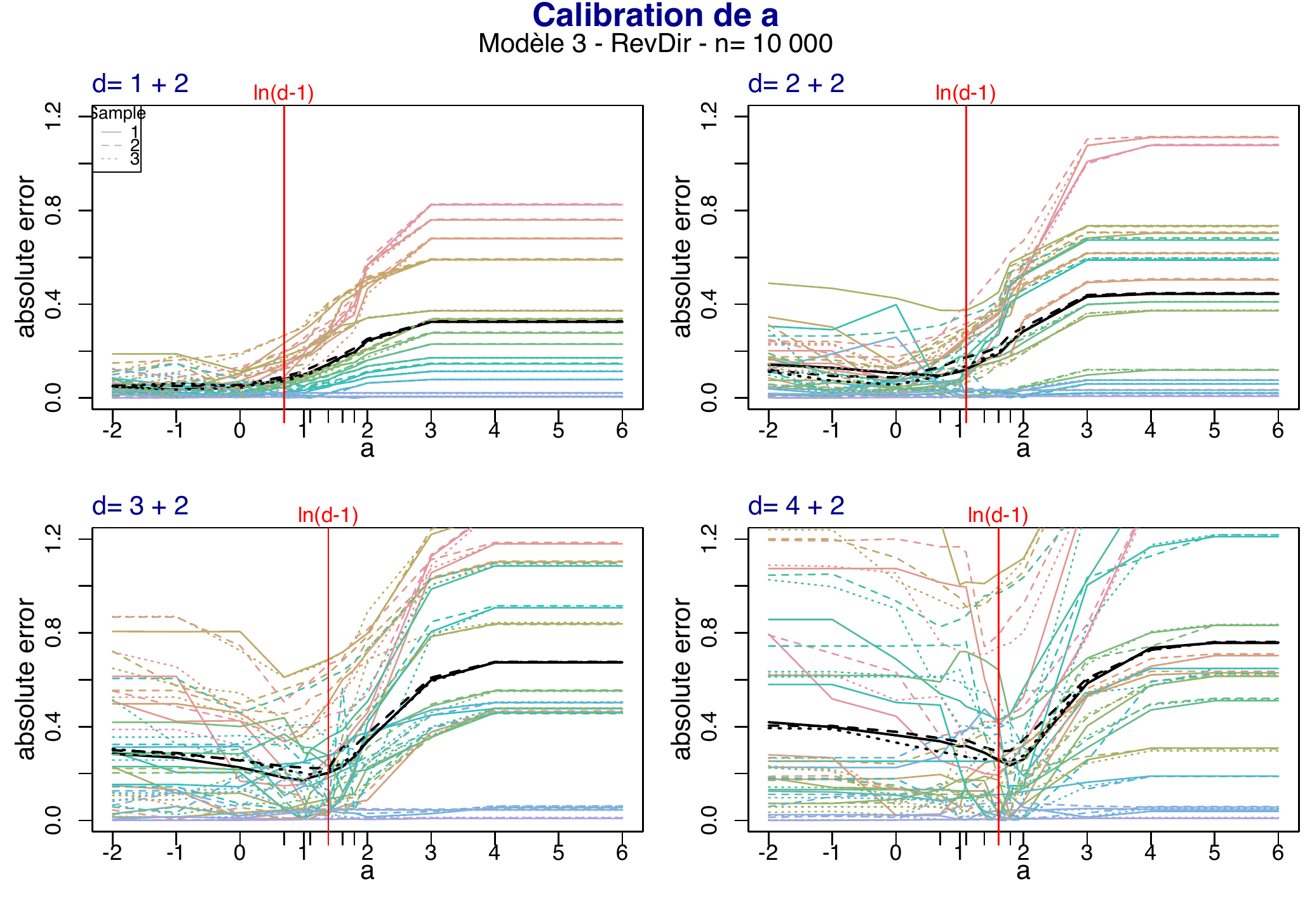} 
\vspace{-0.8cm}
\caption{\textbf{Tuning of $a$ for Example (a) with $n=10 \  000$.} In each subgraph: AE curves in function of $a$ for $16$ evaluation points $w^k$ (the warmer the pastel color, the larger $f(w^k)$) given $B=3$ samples (differentiated by line type) at fixed dimension (specified top left in the form $d=d_1+d_2$). In black lines: the average per sample of the $16$ pastel curves. The vertical straight red line: our final choice.}
\label{fig Calib a a 10000}
\vspace{-0.0cm}
\includegraphics[width=\linewidth, trim=0cm 0.5cm 0cm 1cm,clip]{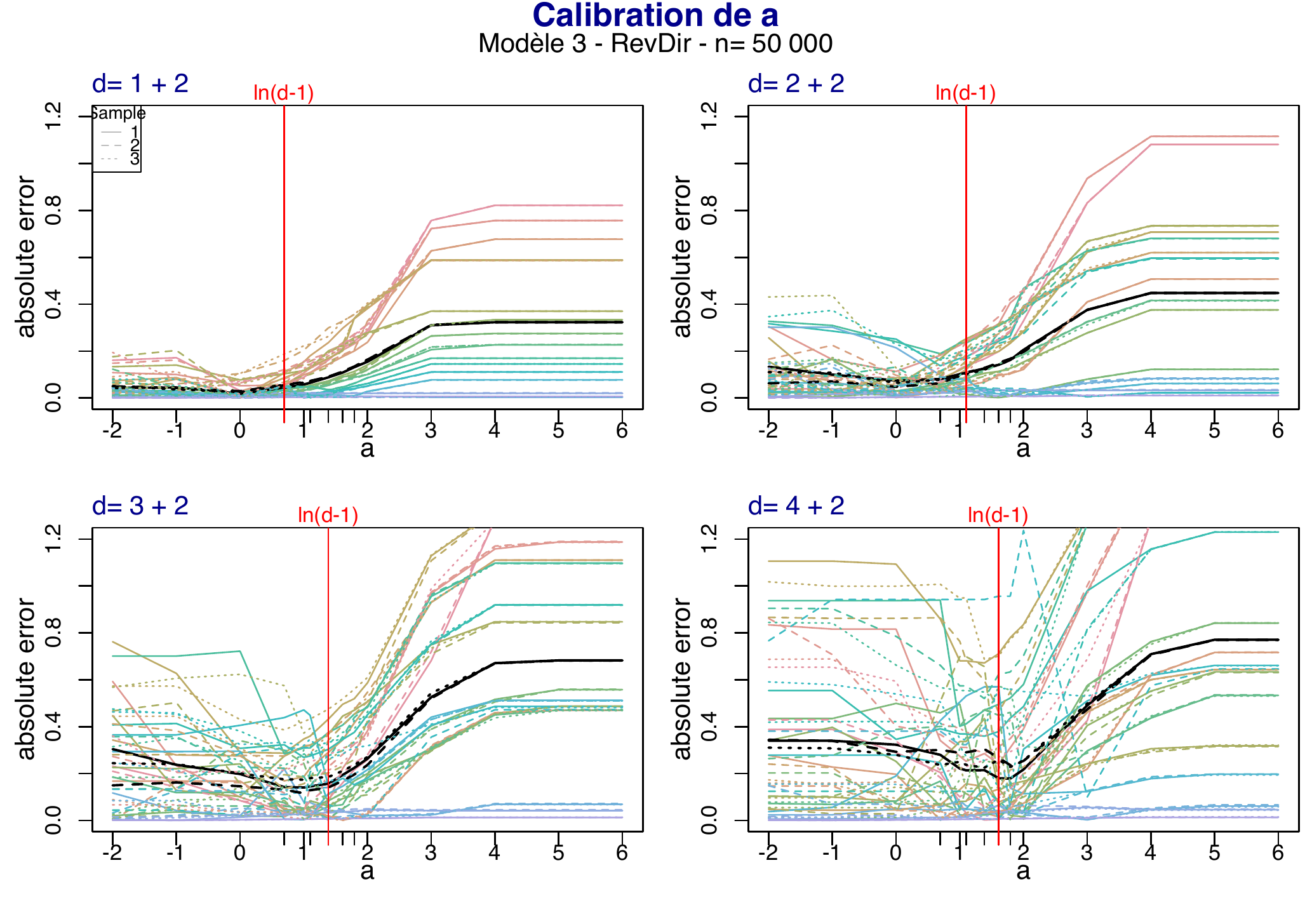} 
\vspace{-0.5cm}
\caption{\textbf{Tuning of $a$ for Example (a) with $n=50 \  000$.} Same description as for Figure~\ref{fig Calib a a 10000}.}\label{fig Calib a a 50000}
\end{figure}
\begin{figure}[p]
\includegraphics[width=\linewidth, trim=0cm 0.5cm 0cm 1cm,clip]{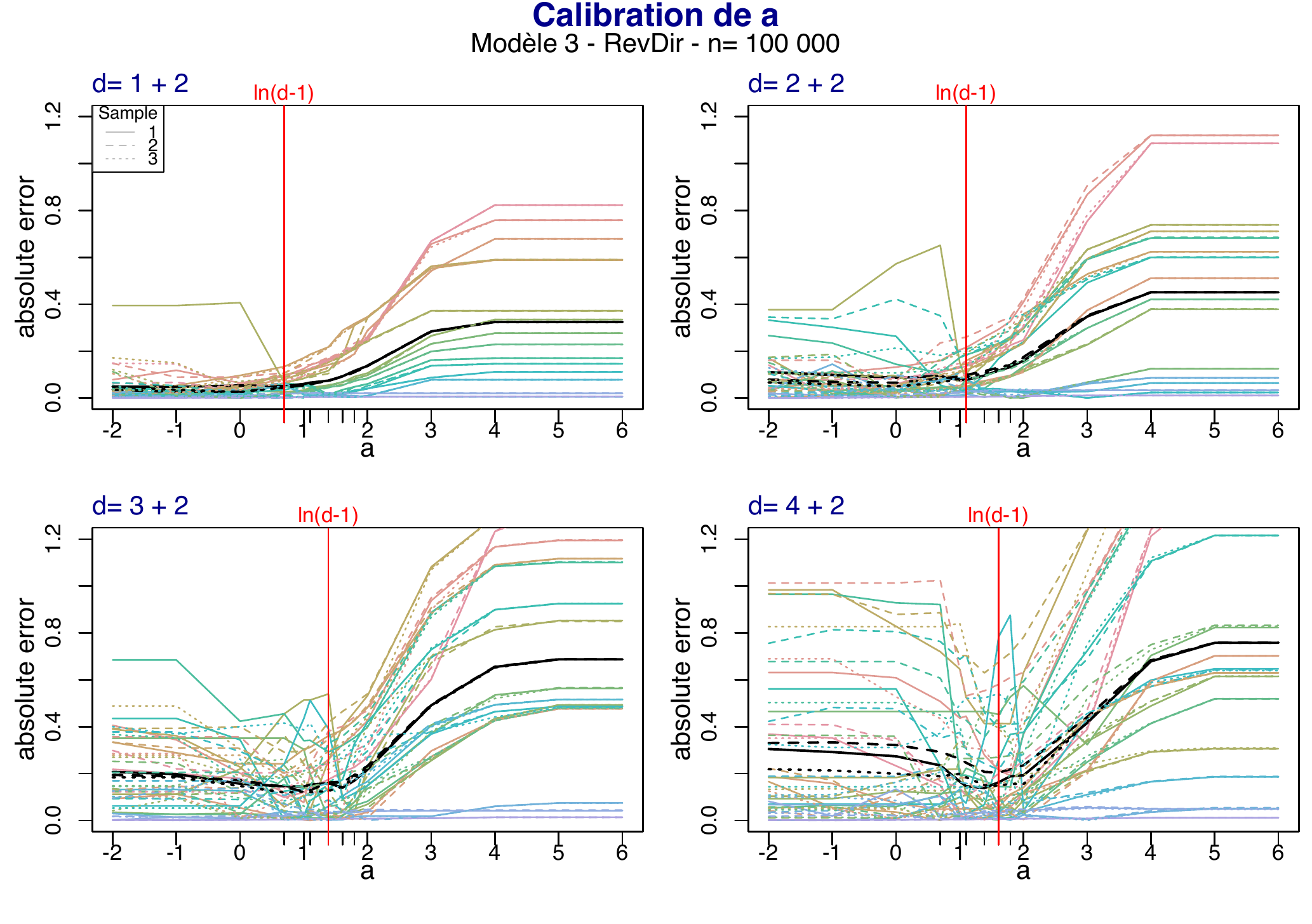} 
\vspace{-0.5cm}
\caption{\textbf{Tuning of $a$ for Example (a) with $n=100 \  000$.} Same description as for Figure~\ref{fig Calib a a 10000}.}\label{fig Calib a a 100000}
\vspace{0.5cm}
\includegraphics[width=\linewidth, trim=0cm 0.5cm 0cm 1cm,clip]{CalibrationaModel3RevDirn200000.pdf} 
\vspace{-0.5cm}
\caption{\textbf{Tuning of $a$ for Example (a) with $n=200 \  000$.} Same description as for Figure~\ref{fig Calib a a 10000}.}
\label{fig Calib a a 200000}
\end{figure}
%\begin{figure}[t]
%\includegraphics[width=\linewidth, trim=0cm 0.5cm 0cm 1cm,clip]{Calibration_a-Model3-RevDir-n200000.pdf} 
%\vspace{-0.5cm}
%\caption{\textbf{Illustration of the tuning of $a$ for Example (a) with $n=200 \  000$ for growing dimensions.} Same description as for Figure~\ref{fig Calib a a 10000}.}
%\label{fig Calib a a 200000}
%\end{figure}
\begin{figure}[p]
\includegraphics[width=\linewidth, trim=0cm 0.5cm 0cm 0.7cm,clip]{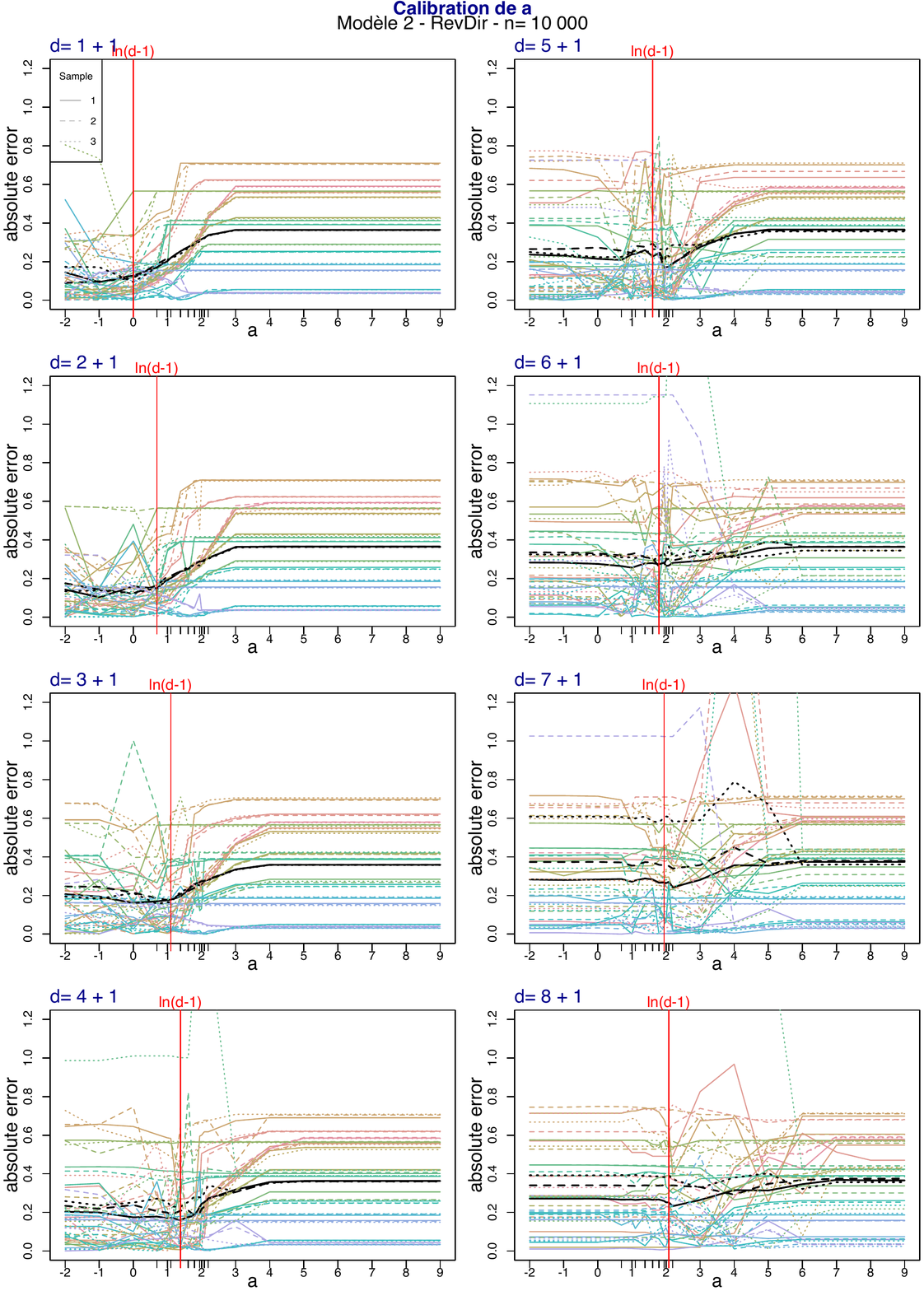} 
\vspace{-0.5cm}
\caption{\textbf{Tuning of $a$ for Example (b) with $n=10 \  000$.} Same description as for Figure~\ref{fig Calib a a 10000}.}
\label{fig Calib a b 10000}
\end{figure}
\begin{figure}[p]
\includegraphics[width=\linewidth, trim=0cm 0.5cm 0cm 0.7cm,clip]{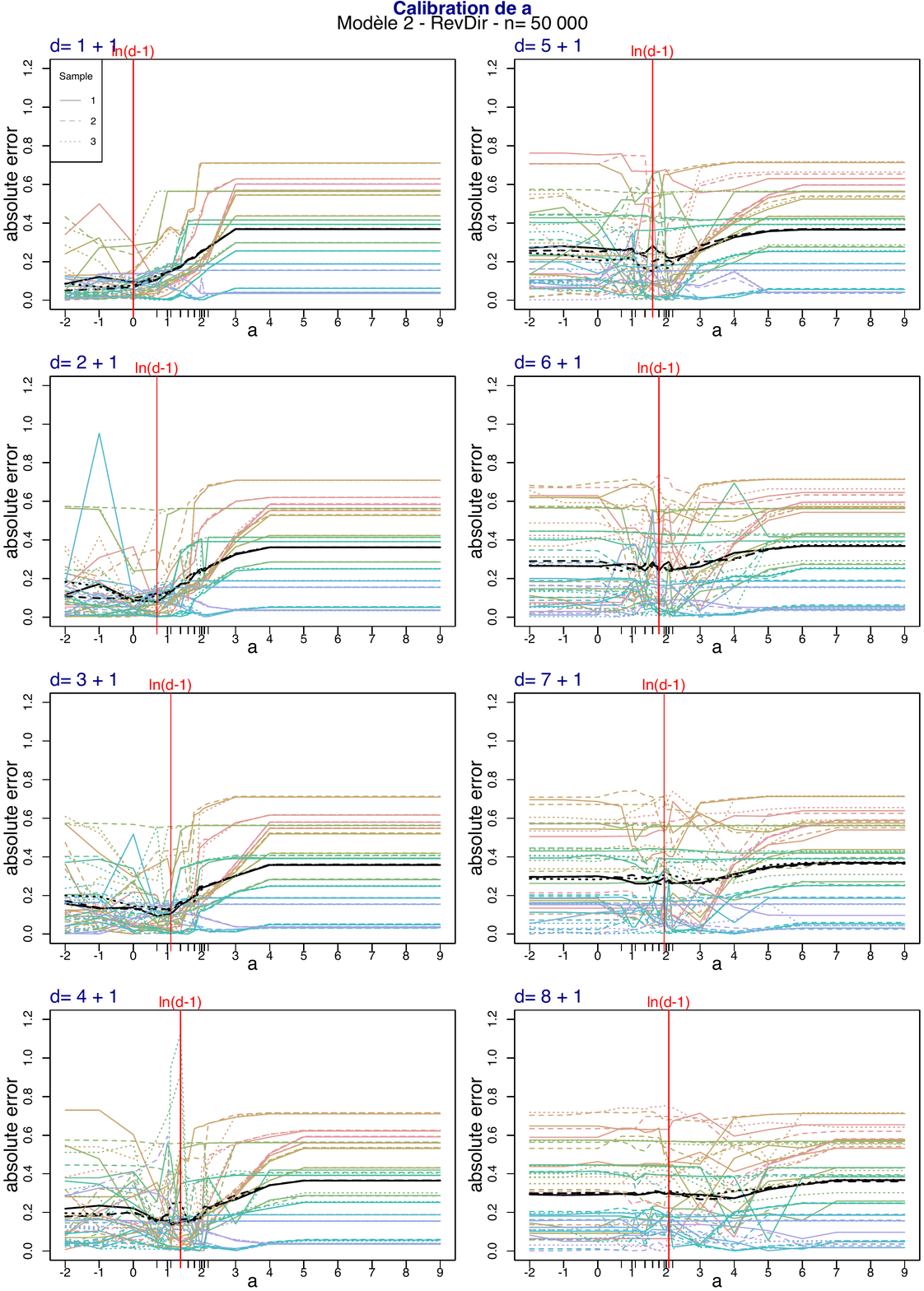} 
\vspace{-0.5cm}
\caption{\textbf{Tuning of $a$ for Example (b) with $n=50 \  000$.} Same description as for Figure~\ref{fig Calib a a 10000}.}
\label{fig Calib a b 50000}
\end{figure}
\begin{figure}[p]
\includegraphics[width=\linewidth, trim=0cm 0.5cm 0cm 0.7cm,clip]{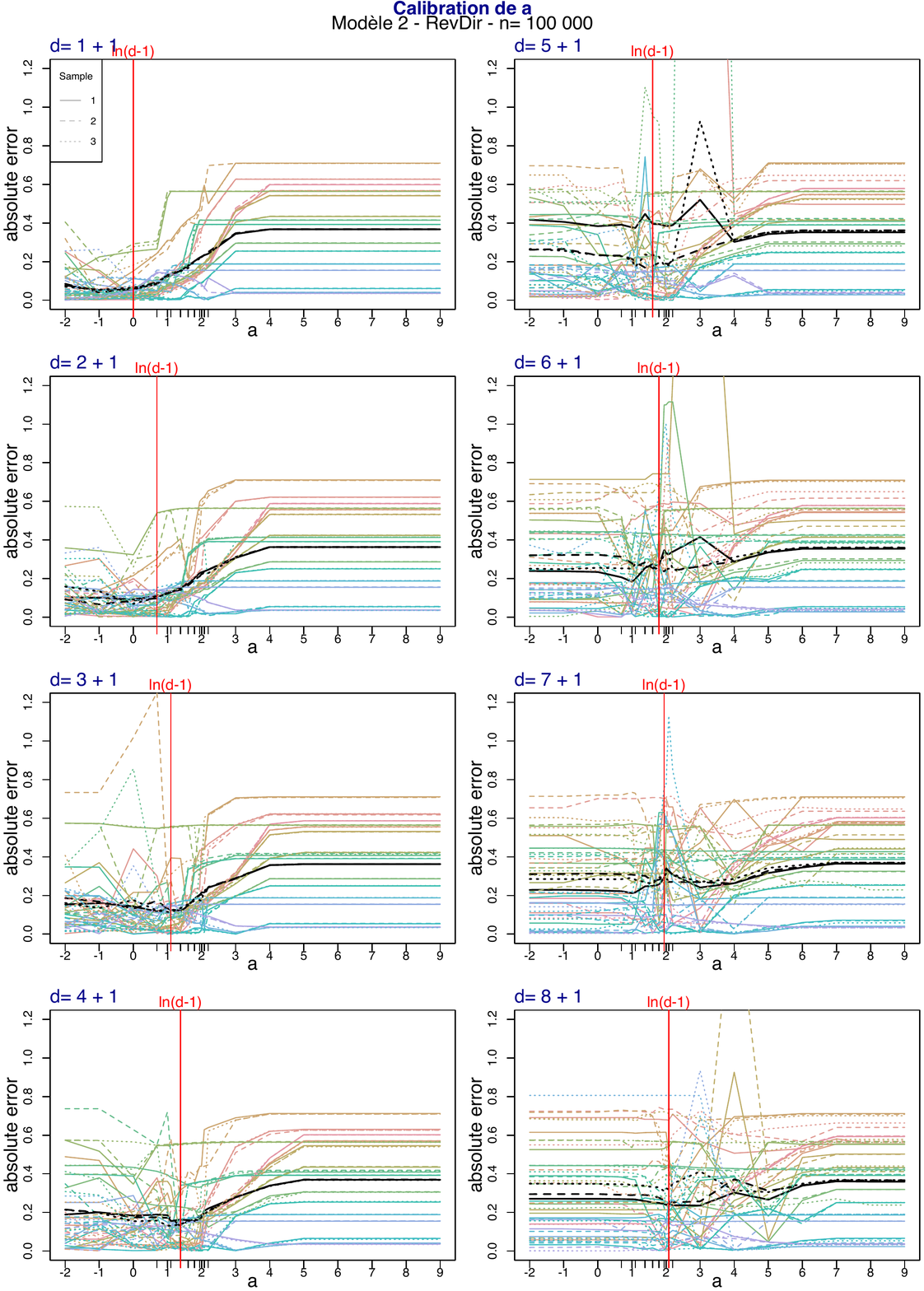} 
\vspace{-0.5cm}
\caption{\textbf{Tuning of $a$ for Example (b) with $n=100 \  000$.} Same description as for~Figure~\ref{fig Calib a a 10000}.}
\label{fig Calib a b 100000}
\end{figure}
\begin{figure}[p]
\includegraphics[width=\linewidth, trim=0cm 0.5cm 0cm 0.7cm,clip]{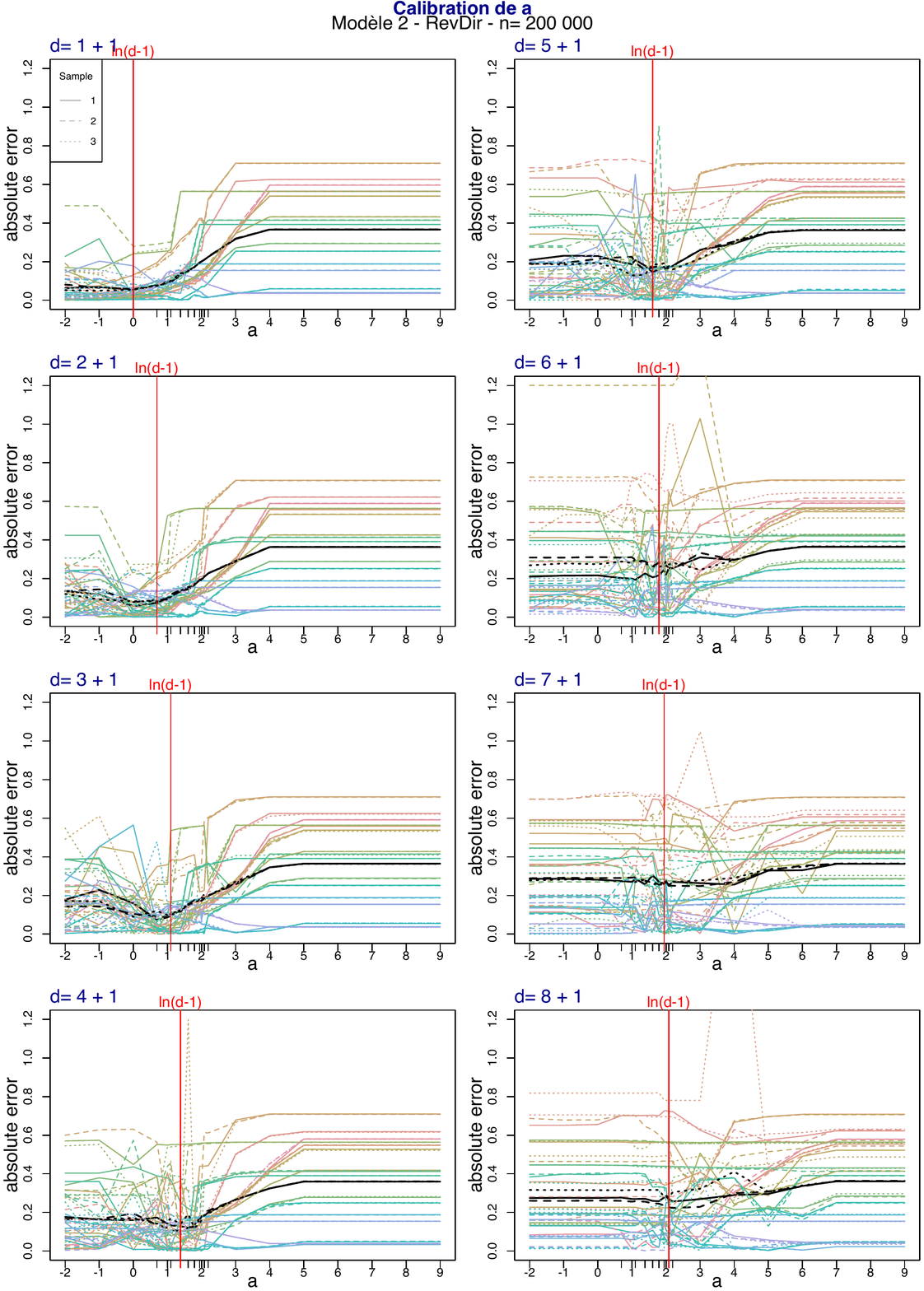} 
\vspace{-0.5cm}
\caption{\textbf{Tuning of $a$ for Example~(b) with $n=200 \  000$.} Same description as for~Figure~\ref{fig Calib a a 10000}.}
\label{fig Calib a b 200000}
\end{figure}
\begin{figure}[p]
\includegraphics[width=\linewidth, trim=0cm 0.5cm 0cm 0.7cm,clip]{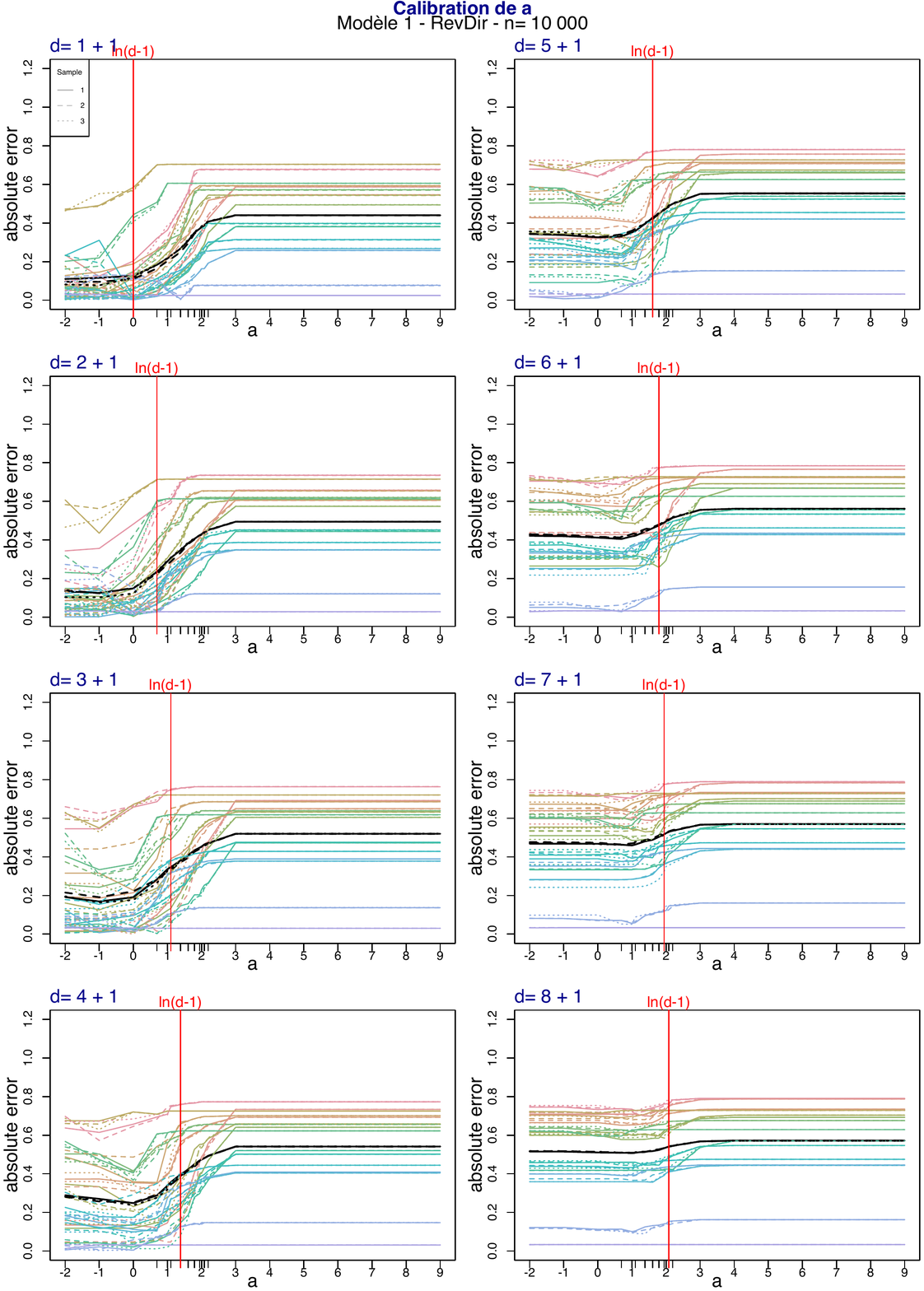} 
\vspace{-0.5cm}
\caption{\textbf{Tuning of $a$ for Example (c) with $n=10 \  000$.} Same description as for Figure~\ref{fig Calib a a 10000}.}
\label{fig Calib a c 10000}
\end{figure}
\begin{figure}[p]
\includegraphics[width=\linewidth, trim=0cm 0.5cm 0cm 0.7cm,clip]{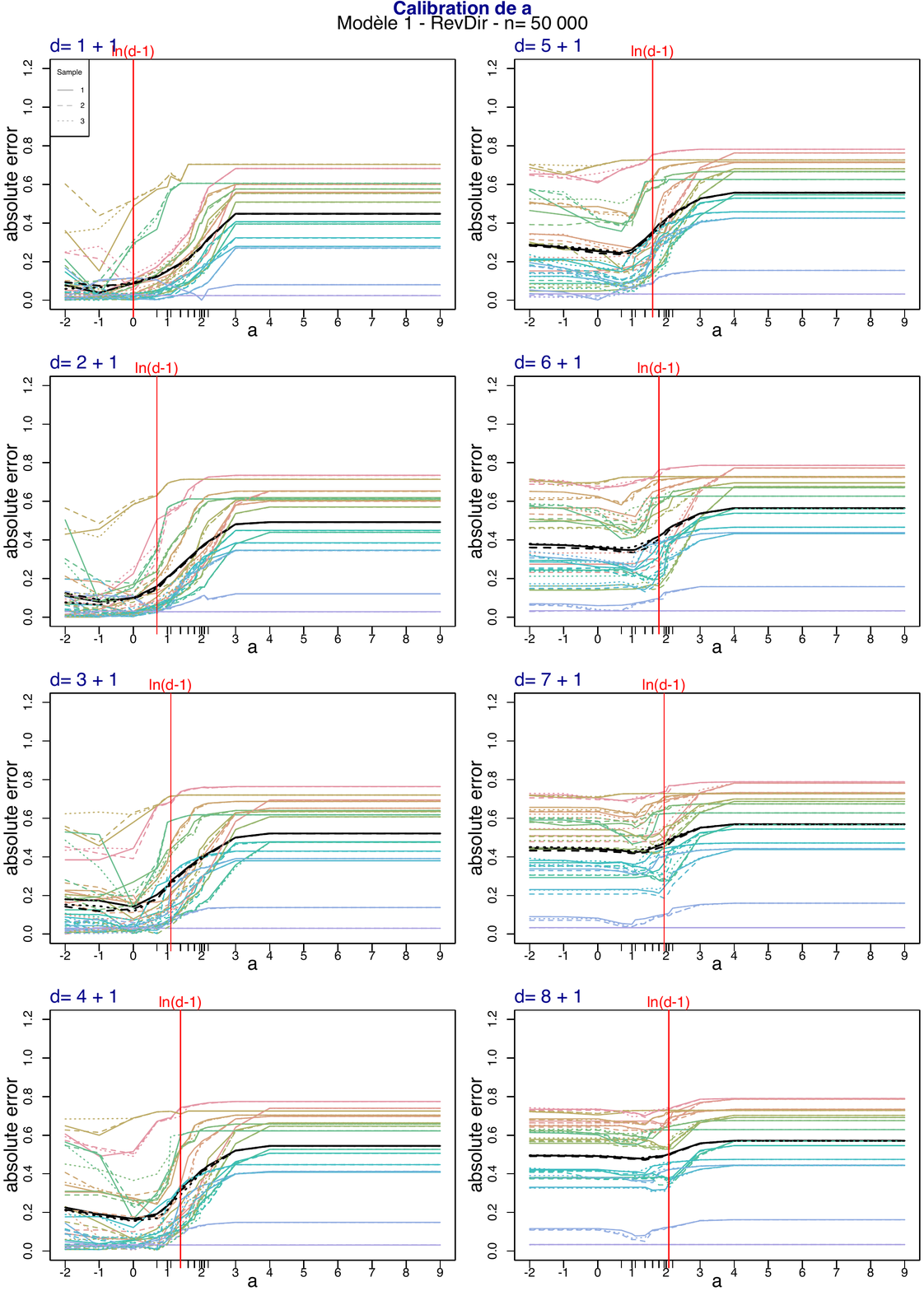} 
\vspace{-0.5cm}
\caption{\textbf{Tuning of $a$ for Example (c) with $n=50 \  000$.} Same description as for Figure~\ref{fig Calib a a 10000}.}
\label{fig Calib a c 50000}
\end{figure}
\begin{figure}[p]
\includegraphics[width=\linewidth, trim=0cm 0.5cm 0cm 0.7cm,clip]{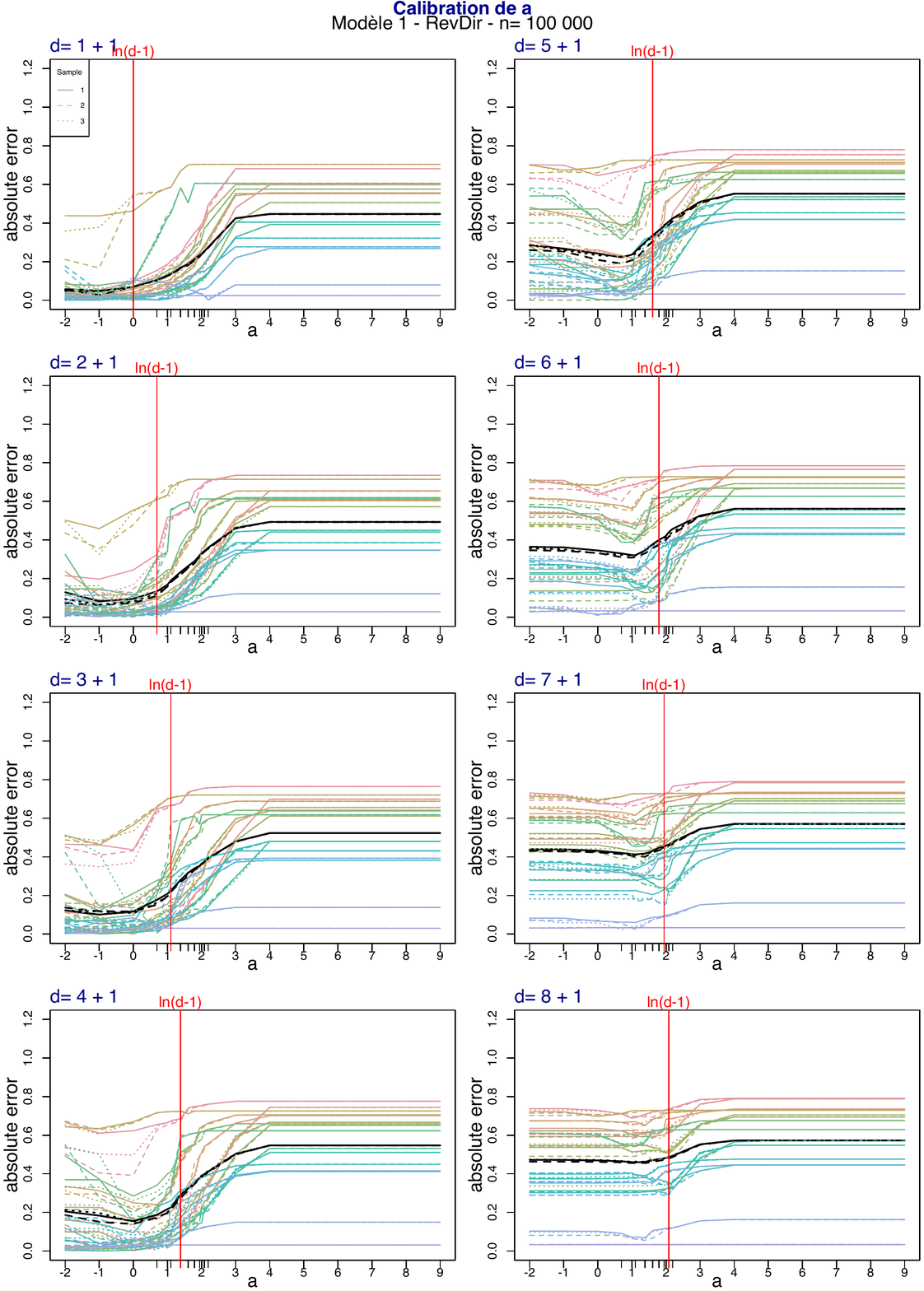} 
\vspace{-0.5cm}
\caption{\textbf{Tuning of $a$ for Example (c) with $n=100 \  000$.} Same description as for~Figure~\ref{fig Calib a a 10000}.}
\label{fig Calib a c 100000}
\end{figure}
\begin{figure}[p]
\includegraphics[width=\linewidth, trim=0cm 0.5cm 0cm 0.7cm,clip]{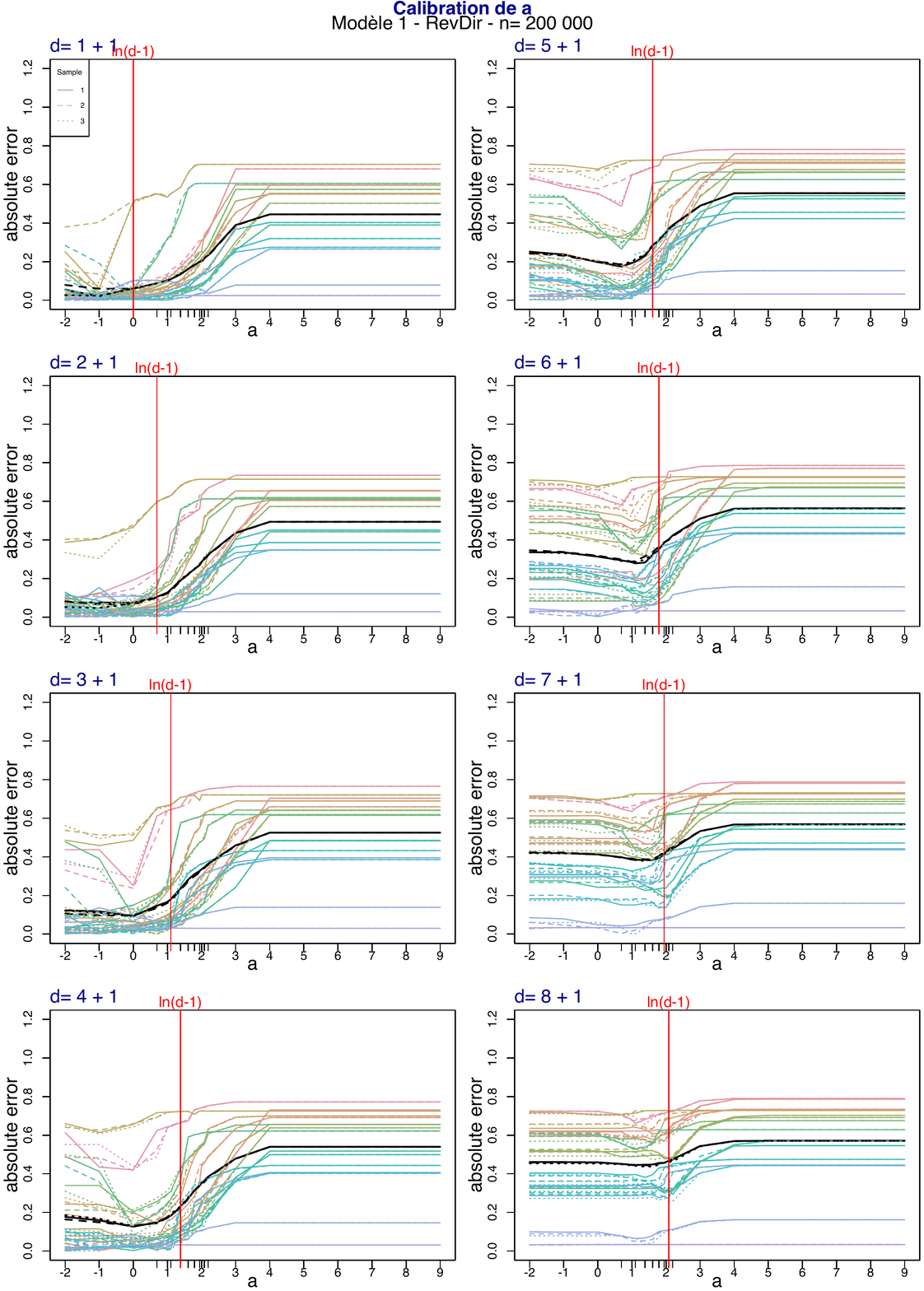} 
\vspace{-0.5cm}
\caption{\textbf{Tuning of $a$ for Example (c) with $n=200 \  000$.} Same description as for~Figure~\ref{fig Calib a a 10000}.}
\label{fig Calib a c 200000}
\end{figure}
\clearpage
%%%%%%%%%%%%%%%%%%%%%%%%%%
\subsection{Curse of dimensionality for Example~(a)}\label{sec:curse:supp}
\begin{figure}[h!]
\includegraphics[width=\linewidth, trim=0cm 0cm 0cm 1.5cm,clip]{%SparseRobust_Examples(b)et(c)_BarEst.pdf}
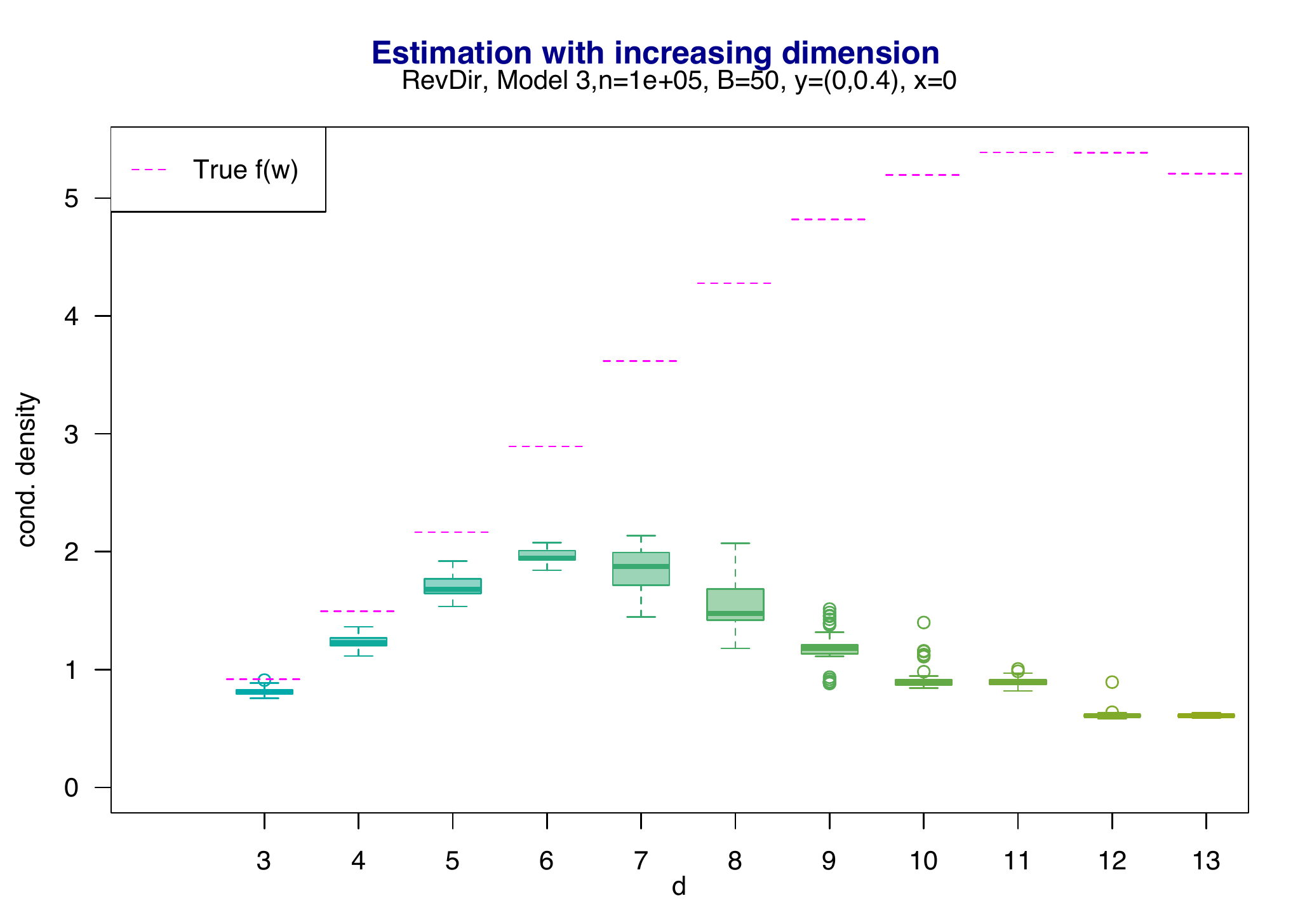}
\vspace{-1cm}
\caption{\label{FigNonSparseEst}\textbf{Estimation with increasing dimension.} Boxplots of the estimates $\frod(w)$ in function of $d$, given $50$ samples of size $n=100 \  000$ for Example~(a).  In greenish shades: the estimates of Example (a). The purple dashed horizontal segment is the true value $f(w)$ (evaluated for $x=0$ and $y=(0,0.5)$).}
\end{figure}

The following example illustrates the difficulties met by our procedure to face with the curse of dimensionality when data have no sparsity structure.
Figure~\ref{FigNonSparseEst}  provides the boxplots of 50 simulated samples of size $n = 100 \ 000$ with varying dimension $d_1$ from $1$ to $11$.

Without sparsity our method struggles providing good estimates as soon as the dimension is larger than $5$. To the best of our knowledge, all classical kernel procedures cope with same difficulties.
%%%

%%%%%%%%%%%%%%%%%%%%%%%%%%
%%%%%%%%%%%%%%%%%%%%%%%%%%

\bibliographystyle{apalike}
\bibliography{biblio}

%% file: preuvesrodeo.tex
% !TEX root =Rodeo-Adap-submitted.tex

%%%%%%%%%%%%%%%%%%%%%%%%
%%%%%%%%%%%%%%%%%%%%%%%%
\section{Proofs}\label{sec:proofs}
\subsection{Notations}\label{ssectNotations}
In order to prove the theorem, some intermediate lemmas are needed. See Appendix for their statements.
%The lemmas are mainly picked up from \cite{Jeanne1}, see  Nonetheless, note that some refinements have been made, see in particular the improved proofs of Lemma \ref{lmDelta}, and the points 2. of \Cref{lmfbar,lmZbar} in Appendix.\\
First, we define some general notations: We denote
\begin{itemize}
\item $\partial_j g$ the partial derivative of a function $g$ with respect to its $j$-th component;
\item $v\cdot v'$ the multiplication term by term of two vectors $v$ and $v'$;
%\item $l:m$ the set of consecutive integers from $l$ to $m$;
\item $v_{\mathcal{I}}$ the vector $v$ restricted to its components indexed in $\mathcal{I}$;
\item $b\vee c=\max(b,c)$ the maximum value of two reals $b$ and $c$.
\end{itemize}
 Let us now introduce the key quantities of the proofs. For any bandwidth $h\in(\R_+^*)^d$ and any component $k\in 1:d$, we consider the estimator $\fbar{h}$ that we would have used if the density $\fX$ were known:
$$\fbar{h}:=\frac1{n} \sum\limits_{i=1}^n \fbar{hi},\quad \fbar{hi}:= \frac{\K_h(w-W_i)}{\fX(X_i)}$$
and we denote $\Delta_h$ its difference with the real estimator:
$$\Delta_h:= \f{h}-\fbar{h}.$$
We denote $\biasbar{h}:= \E\left[\fbar{h}\right]-f(w)$ the bias of $\fbar{h}$.
We also consider its partial derivative~$\bar{Z}_{hk}$:
$$ \bar{Z}_{hk} :=  \frac{\partial}{\partial h_k}\fbar{h}.$$
We can write $$\bar{Z}_{hk}:=\frac{1}{n}\sum_{i=1}^n \bar{Z}_{hik}, \quad \bar{Z}_{hik}:= \frac{1}{\fX\left(X_i\right)} \dfrac{\partial}{\partial h_k}\left(\prod\limits_{k=1}^{d}h_k^{-1}K\left(\frac{{w}_k-W_{ik}}{h_k}\right)\right).$$
We shall consider $\DZ{hk}$ the difference between $Z_{hk}$ and $\bar{Z}_{hk}$: 
$$\DZ{hk}:=Z_{hk}-\bar{Z}_{hk}.$$
Note that the value of the final bandwidth of our procedure provides the value of the bandwidth at each iteration. More precisely, if a bandwidth $h$ is the output of the RevDir procedure, we denote $(h^{(t)})_{t\in\ZZ}$,  the different values of the bandwidth for all iterations $t$.\\ 
- On the one hand, if  $h_k>h_0$, it means that  at Initialization, the component $k$ was in $\Act^{(-1)} $ and then the bandwidth path of this component has increased during the Reverse Step according to the following path $h_0\beta^{-1}, h_0\beta^{-2},...$ until $h_k:=h_0\beta^{-|t_k|}$, and remains fixed during the whole Direct Step ($t\geq 0$).\\
- On the other hand, if $h_k<h_0$, the component $k$ was in $\Act^{(0)}$ at Initialization. Thus the value of the bandwidth component was fixed and equals to $h_0$ during the Reverse Step (i.e for every $t<0$ ). Then, it  decreases during the Direct step: $h_0\beta,h_0\beta^2,...$ until $h_k:=h_0\beta^{t_k}$ is achieved (see Figure~\ref{figrodeo}). This gives the following formula: for any $k \in 1:d$, during the Reverse Step (when $t<0$),
$$\quad h^{(t)}_k:=\max(h_0,\min\left(h_k, \beta^{t}h_0\right))
=\begin{cases}
h_0 &\text{ if } k \text{ is active during the Direct Step,}\\
\beta^{t}h_0&\text{ if } k \text{ is active during the Reverse Step and not deactivated yet,}\\
h_k &\text{ if } k \text{ has already been deactivated during the Reverse Step,}
\end{cases}$$
and during Direct Step (when $t\geq 0$),  
$$h^{(t)}_k:=\max\left(h_k, \beta^{t}h_0\right)
=\begin{cases}
\beta^{t}h_0 &\text{ if } k \text{ is active during the Direct Step and not deactivated yet,}\\
h_k &\text{ if } k \text{ has already been deactivated (during the Reverse or the Direct Step).}
\end{cases}$$
Now we can define the set of bandwidths $\Hhp$ which contains with high probability the bandwidth selected by the RevDir procedure:
\begin{align*}
\Hhp := \lbrace
h\in\left(\R_+^*\right)^d :\,  &\forall k\in 1:d, h_k=\beta^{t_k}h_0\leq 1 \text{ with } t_k\in\ZZ,  \\
&\, \text{ and } \prod\limits_{k=1}^d h_k\geq \beta^r \tfrac{(\log n)^{a+1}}{n},\\
&\text{ and } \forall k\in\rel^c, h_k=\hirr
\rbrace,
\end{align*}
where $\hirr$ is uniquely defined by $\tirr\in\mathbb Z$ such that
$\beta<\hirr:=\beta^{\tirr} h_0\leq 1$.
We also denote $\HhpRev$ (respectively $\HhpDir$) the set which contains the different states of the bandwidth during the Reverse Step (respectively the Direct Step) provided that the selected bandwidth is in $\Hhp$:
\begin{equation}
\HhpRev := \lbrace h^{(t)} : h\in\Hhp, t<0\rbrace
\end{equation}
\begin{equation}
\HhpDir := \lbrace h^{(t)} : h\in\Hhp, t\geq 0\rbrace.
\end{equation}
Finally, we introduce the high probability event $\Ehp$ on which $\hat{h}$ systematically belongs to $\Hhp$:
\begin{equation}
\Ehp:= \Atilde 
\cap  \bigcap\limits_{h\in\Hhp} \left( \Bernf{h} \cap \Bernfabs{h} \right)
\cap  \bigcap\limits_{h\in (\HhpRev\cup\HhpDir)} \bigcap\limits_{k=1}^d\left( \BernZ{h,k} \cap \BernZabs{h,k} \right),
\end{equation}
where $\Atilde$ is the high probability event of \ref{fXtildeAccuracy} in Assumption \AestfX:
$$\Atilde=\left\lbrace\sup\limits_{u\in\Un}\left\vert \fX(u)-\hfX(u)\right\vert \leq M_{X} \frac{(\log n)^{\frac{a}2}}{\sqrt{n}}\right\rbrace,$$
and $\mathcal{B}\text{ern}_{\dagger}(\ddagger)$ is the high probability event resulting of Bernstein's Inequality applied on the random variable $\dagger$ with parameter(s) $\ddagger$. More formally:
$$\Bernf{h}:=\lbrace \vert \fbar{h}-\mathbb{E}[\fbar{h}]\vert \leq \sigma_h\rbrace,$$
$$\Bernfabs{h}:=\left\{ \left\vert \frac1n\sum\limits_{i=1}^n \vert\fbar{hi}\vert-\mathbb{E}[\vert\fbar{h}\vert]\right\vert \leq \CEbar \right\},$$
$$\BernZ{h,k}:=\left\{ \vert \bar{Z}_{hk}-\E{\bar{Z}_{hk}}\vert\leq\frac12\seuil{hk} \right\},$$
$$\BernZabs{h,k}:=\left\{\left\vert \frac1n \sum\limits_{i=1}^n \vert\bar{Z}_{hik}\vert -\E{\vert\bar{Z}_{h1k}\vert}\right\vert\leq \CEZA h_k^{-1} \right\},$$
where $$\sigma_h=\Cs\sqrt{ \frac{(\log n)^a}{ n \prod\limits_{k=1}^d h_k}}$$ with
$ \Cs=\frac{2\Vert K\Vert^d_2 \Vert f\Vert_{\infty,\UUn}^{\frac12}}{\delta^{\frac12}}$.
See \Cref{lmfbar,lmZbar} in Appendix for the details and definitions of 
constants $\CEbar, \CEZA$. 
%%%%%%%%%%%%%%%%%%%
\subsection{Main steps of the proof}
\Cref{propHhp} describes the form of the bandwidth selected by the RevDir procedure with high probability. Given this selection, \Cref{propBiasDeviation} gives upper bounds on the bias and the deviation of the estimator $\bar{f}_{\hat{h}}(w)$.
\begin{proposition} \label{propHhp}The selected bandwidth belongs to $\Hhp$ with high probability. More precisely:
\begin{equation}\label{eProp2hinHhp}
\Ehp\subset \lbrace \hat{h}\in\Hhp \rbrace 
\end{equation}
and for $n$ large enough:
\begin{equation}\label{eProp2MajEhpc}
\P\left( \Ehp^c\right)\leq 2 e^{-(\log n)^{1+\frac{a-1}2}}.
\end{equation}
\end{proposition}
Note in particular that with high probability the irrelevant components of the  selected bandwidth are equal to $\hirr$.\\

Recall that $\biasbar{h}:= \E\left[\fbar{h}\right]-f(w)$ is the bias of $\fbar{h}$.
\begin{proposition}\label{propBiasDeviation}
The following upper bounds are satisfied for all $h\in\Hhp$, and any constants $\Am\in\R$ and $\CAm>0$:

\begin{align}
\1_{\lbrace\hat{h}=h\rbrace\cap\Ehp} 
&\left\vert  
\biasbar{h}
\right\vert
\leq  r \CBbar{\CAm}^s \frac{\left(\log n\right)^{\Am s}}{n^{\frac{s}{2s+r}}}
+ r \max\left( \tfrac{7\Cl }{4\beta%{\hirr}
^{\frac{d-r}2}{\CAm}^{\frac{r}2}} \frac{(\log n)^{\frac{a -Ar}2}}{n^{\frac{s}{2s+r}}}
, {\rcp} \left(\frac{(\log n)^{a}}{n}\right)^{\frac{p}{2p+1}}\right),\label{eMajBias}\\
\1_{\lbrace\hat{h}=h\rbrace\cap\Ehp} 
&\left\vert  
\fbar{h}-\E\left[ \fbar{h}\right]
\right\vert
\leq \1_{\lbrace\hat{h}=h\rbrace\cap\Ehp} \sigma_h  \nonumber
\\
&\leq \max\left( \tfrac{\Cs}{\beta%{\hirr}
^{(d-r)/2}{\CAm}^{r/2}} (\log n)^{(a-\Am r)/2},
 \tfrac{4{\CAm}^s\CEZ\Cs\beta^{-\frac{r}2-s}}{\Cl}  (\log n)^{s\Am} \right)
n^{-\frac{s}{2s+r}},\label{eMajDeviation}
\end{align}
where $\Cl$ is the constant defined in \eqref{threshold} and $\CBbar,\Cs,\CEZ$ are constants defined  in \Cref{lmfbar,lmZbar}  in Appendix.
\end{proposition}
\subsection{Proof of Theorem \ref{thm}}\label{sec:proofmainthm}
Let us fix $l>1$. From \Cref{propHhp}: $\Ehp\subset \lbrace \hat{h}\in\Hhp \rbrace$, thus:
\begin{equation}\label{eDecompRiskl}
\E\left[ \left\vert \frod(w)-f(w)\right\vert^l \right]
= \E\left[\1_{\Ehp^c}\left\vert \frod(w)-f(w)\right\vert^l\right]
+\sum\limits_{h\in\Hhp} 
\E\left[\1_{\lbrace \hat{h}=h\rbrace \cap \Ehp}\left\vert \hat{f}_h(w)-f(w)\right\vert^l\right].
\end{equation}
We first control the terms $\E\left[\1_{\lbrace \hat{h}=h\rbrace \cap \Ehp}\left\vert \hat{f}_h(w)-f(w)\right\vert^l\right]$. 
We fix $h\in\Hhp$. Then, we decompose the difference $\hat{f}_h(w)-f(w)$ as follows:
\begin{equation}
\hat{f}_h(w)-f(w)
=\Delta_h + \left(\fbar{h} - \E\left[\fbar{h}\right]\right) + \biasbar{h},
\end{equation}
where we recall the notations $\Delta_h:=\hat{f}_h(w) - \fbar{h}$ and $\biasbar{h}:=\E\left[\fbar{h}\right] - f(w)$. Remark that $\prod_{k=1}^d h_k\leq 1$, since $h\in\Hhp$.
We apply 2. of \Cref{lmDelta} and 3. of \Cref{lmfbar}. Since $\Ehp\subset \left( \Atilde \cap \Bernfabs{h}\right)  \cap \Bernf{h}$:
$$\1_{\Ehp} \left\vert \Delta_h \right\vert
\leq  \CMD\sigma_h
$$
and 
$$ \1_{\Ehp} \left\vert \fbar{h}- \E\left[\fbar{h}\right] \right\vert 
\leq  \sigma_h.$$
Therefore:
\begin{equation}\label{eqPourSuppThm}
\1_{\lbrace \hat{h}=h\rbrace \cap \Ehp}\left\vert \hat{f}_h(w)-f(w)\right\vert
\leq 
\1_{\lbrace \hat{h}=h\rbrace \cap \Ehp} 
\left( ( \CMD + 1) \sigma_h 
+ \left\vert \biasbar{h}\right\vert\right).
\end{equation}
From \Cref{propBiasDeviation} which controls both $\sigma_h$ and $\left\vert \biasbar{h}\right\vert$, we deduce:
\begin{align*}
\1_{\lbrace \hat{h}=h\rbrace \cap \Ehp}
&\left\vert \hat{f}_h(w)-f(w)\right\vert\\
&\leq 
	\left( \CMD + 1\right)
	\max\left( 
		\tfrac{\Cs}{{\beta}^{(d-r)/2}{\CAm}^{r/2}} (\log n)^{\frac{a-\Am r}2},
 		\tfrac{4{\CAm}^s\CEZ\Cs\beta^{-\frac{r}2-s}}{\Cl}  (\log n)^{s\Am} 
 	\right) 
 	n^{-\frac{s}{2s+r}}
\\
&\quad+ r \CBbar{\CAm}^s \left(\log n\right)^{\Am s} n^{-\frac{s}{2s+r}}
+ r \max\left( 
	\tfrac{7\Cl }{4\beta^{\frac{d-r}2}{\CAm}^{\frac{r}2}} \frac{(\log n)^{\frac{a-Ar}2}}{n^{\frac{s}{2s+r}}},
	{\rcp} \left(
		\frac{(\log n)^{a}}{n}
	\right)^{\frac{p}{2p+1}}
\right).
\end{align*}
%%%%%%%%%%%%%%%%%%%%%%%%%%%%%%%%%%%%%%%%%%%%%%%
%We optimize in $\Am$ and $\CAm$: in particular, denoting
%$C_{opt,1}:=\max\left(r\CBbar,
% \tfrac{4\CEZ\Cs\beta^{-\frac{r}2-s}\left( 1+\tfrac{\CMD}{(\log n)^{\frac{a}2}} \right)}{\Cl}\right)$ 
%and 
%$C_{opt,2}:=\max \left( \tfrac{\left( 1+\tfrac{\CMD}{(\log n)^{\frac{a}2}} \right)\Cs}{{\hirr}^{d-r}},
%\tfrac{3r\Cl }{4\hirr^{\frac{d-r}2}} \right)$,
% for $\Am=\frac{a}{2s+r}$ and $\CAm=\left(\frac{C_{opt,2}}{C_{opt,1}}\right)^{\frac{2}{2s+r}}$, we obtain:
%\begin{align*}
%\1_{\lbrace \hat{h}=h\rbrace \cap \Ehp}
%\left\vert \hat{f}_h(w)-f(w)\right\vert
%&\leq 3
%\max\left(
%	\left(\frac{C_{opt,2}}{C_{opt,1}}\right)^{\frac{2s}{2s+r}}
% 	\left(\frac{(\log n)^a}{n}\right)^{\frac{s}{2s+r}},
%	r\rcp \left(
%		\frac{(\log n)^{a}}{n}
%	\right)^{\frac{p}{2p+1}}
%\right).
%\end{align*}
%The first term is dominant: $C_p$ are defined such that when $r=1$ and $s=p$: 
%$$r\rcp(\log n)^{\frac{a.p}{2p+1}} \leq \left(\frac{C_{opt,2}}{C_{opt,2}}\right)^{\frac{2s}{2s+r}}(\log n)^{\frac{as}{2s+r}},$$ 
% and otherwise $  n^{-\frac{p}{2p+1}}\ll n^{-\frac{s}{2s+r}}$. Therefore, for $n$ large enough:
%\begin{align}
%\1_{\lbrace \hat{h}=h\rbrace \cap \Ehp}
%\left\vert \hat{f}_h(w)-f(w)\right\vert
%&\leq 3
%	\left(\frac{C_{opt,2}}{C_{opt,1}}\right)^{\frac{2s}{2s+r}}
% 	\left(\frac{(\log n)^a}{n}\right)^{\frac{s}{2s+r}}. \label{eMajErreurEhp}
%\end{align}
%%%%%%%%%%%%%%%%%%%%%%%%%%%%%%%%%%%%%%%%%%%%%
We optimize in $\Am$ and $\CAm$:
With $\Am=\frac{a}{2s+r}$, we obtain
\begin{align*}
\1_{\lbrace \hat{h}=h\rbrace \cap \Ehp}
\left\vert \hat{f}_h(w)-f(w)\right\vert
&\leq 
\max\left(
	C_1
 	\left(\frac{(\log n)^a}{n}\right)^{\frac{s}{2s+r}},
	\rcp r\left(
		\frac{(\log n)^{a}}{n}
	\right)^{\frac{p}{2p+1}}
\right).
\end{align*}
where $C_1$ depends on $\beta,d, r,s,\CBbar, \CEZ,\Cs,\CMD,\Cl.$
If $r=0$, the last term in the right hand side vanishes, otherwise $p/(2p+1)\geq s/(2s+r)$ (since $p\geq s$).
 Therefore, for $n$ large enough:
\begin{align}
\1_{\lbrace \hat{h}=h\rbrace \cap \Ehp}
\left\vert \hat{f}_h(w)-f(w)\right\vert
&\leq 
	\text{C}'
 	\left(\frac{(\log n)^a}{n}\right)^{\frac{s}{2s+r}}. \label{eMajErreurEhp}
\end{align}
To prove the theorem, it then remains to control $\left\vert \hat{f}_h(w)-f(w)\right\vert$ on $\Ehp^c$.
Recall that: $$\prod\limits_{k=1}^d\hat{h}_k\geq \beta^r \frac{(\log n)^{1+a}}{n},$$ and \ref{fXtildemin}: $$\dX:=\inf\limits_{u\in \Un} \hfX(u)>n^{-1/2},$$
then we can roughly bound $\frod(w)$ by: 
\begin{align*}
\left\vert\frod(w) \right\vert
\leq \frac{\Vert K\Vert_{\infty}^d n}{\dX \beta^r(\log n)^{1+a} } 
=o( n^2).
\end{align*}
So:
$$\left\vert \hat{f}_h(w)-f(w)\right\vert^l
=o( n^{2l})=o(e^{2l \log n}).$$
Besides, from \Cref{propHhp}: 
\begin{align*}
\P\left( \Ehp^c\right)
&\leq 2 e^{-(\log n)^{1+\frac{a-1}2}}.
\end{align*}
Note that, since $a>1$,
%$$(\log n)^{1+\frac{a-1}2}\gg (2+\frac12)l \log n,$$
\begin{equation}\label{lntau}
2l \log n+  l \log( n^{\frac12})=o((\log n)^{1+\frac{a-1}2}),
\end{equation}
therefore:
\begin{align*}
\E\left[\1_{\Ehp^c}\left\vert \frod(w)-f(w)\right\vert^l\right]^{1/l}
\leq  \left(\P\left( \Ehp^c\right) e^{2l \log n}\right)^{1/l}
=o( n^{-\frac12}).
\end{align*}
To conclude, we combine Equation \eqref{eDecompRiskl} with the above upper bound and Inequality \eqref{eMajErreurEhp}:
%%%%%%%%%%%%%%%%%%%%%%%%%%%%%%%%%%%%%%%%%%%%%
%\begin{align*}
%\E\left[ \left\vert \frod(w)-f(w)\right\vert^l \right]^{1/l}
%&\leq o( n^{-\frac12}) + 
%\left\lbrace
%\left(3
%	\left(\frac{C_{opt,2}}{C_{opt,1}}\right)^{\frac{2s}{2s+r}}
% 	\left(\frac{(\log n)^a}{n}\right)^{\frac{s}{2s+r}}\right)^l \sum\limits_{h\in\Hhp} \E[\1_{\hat{h}=h}]
% 	\right\rbrace^{1/l}\\
%&\leq \text{C}
% 	\left(\frac{(\log n)^a}{n}\right)^{\frac{s}{2s+r}} + o(n^{-1/2}),
%\end{align*}
%with 
%$$\text{C}:=3\left(\frac{C_{opt,2}}{C_{opt,1}}\right)^{\frac{2s}{2s+r}}.$$
%%%%%%%%%%%%%%%%%%%%%%%%%%%%%%%%%%%%%%%%%%%%%%%
\begin{align*}
\E\left[ \left\vert \frod(w)-f(w)\right\vert^l \right]^{1/l}
&\leq o( n^{-\frac12}) + 
\left\lbrace
\left(\text{C}'
 	\left(\frac{(\log n)^a}{n}\right)^{\frac{s}{2s+r}}\right)^l \sum\limits_{h\in\Hhp} \E[\1_{\hat{h}=h}]
 	\right\rbrace^{1/l}\\
&\leq \text{C}
 	\left(\frac{(\log n)^a}{n}\right)^{\frac{s}{2s+r}},
\end{align*}
with 
$\text{C}$ depending on $d,r,\|f\|_{\infty,\UUn}, \delta, L, s, K, \beta$.
%%%%%%%%%%%%%%%%%%%%%%%
\subsection{Proof of Proposition \ref{propHhp}}
By definition of the procedure, any selected bandwidth $\hat{h}$ satisfies
$$\exists (t_1,\dots,t_d)\in\ZZ^d, \forall k \in 1:d, \hat{h}_k=\beta^{t_k} h_0$$
The loop condition in the Reverse Step imposes for any active component $k$ that at the beginning of an iteration $t\in\ZZ_-$ :
$$\hat{h}^{(t)}_k\leq \beta.$$
At most, $\hat{h}^{(t)}_k$ is multiplied by $\ibeta$. Then after the last update of the component $\hat{h}_k$:
$$\hat{h}_k\leq 1=\ibeta \beta.$$

Now let us prove that on $\Ehp$, the irrelevant components are deactivated at value $\hirr$. It suffices to show that during the initialization, the irrelevant components activate for  Reverse Step, \textit{i.e.}:
$$\rel^c\subset\Act^{(-1)},$$ 
and in the case where $h_0\leq \beta$, it suffices to prove that they remain active at all iterations $t \in -1:\tirr$. Remember that $\tirr\in\ZZ$ is defined such that:
 $\hirr= \beta^{\tirr} h_0.$
 
 Note that %the state of the selected bandwidth determines exactly the states of the bandwidth at all iterations. 
 %In particular, if $\hat{h}=h$, then for any iteration $t$, $\hat{h}^{(t)}=h^{(t)} $ (defined in  \Cref{ssectNotations}). 
 if the irrelevant components remain active at all iteration $t \in -1:\tirr$, then for $k\in\rel^c$, $\hat{h}^{(t)}_k=H^{(t)}_k=\beta^t h_0$. 
It corresponds to the definition of $\Hhp$, since for all $h\in\Hhp$, $t \in  -1:\tirr$ and $k\in\rel^c$, $$h_k^{(t)}=\beta^t h_0.$$
Therefore, there exists $h\in\Hhp$ such that $\hat{h}^{(t)}=h^{(t)} $ for all iterations $t \in -1:\tirr$.% since during the Reverse step, the bandwidth is increasing, and since $\hirr\leq\beta$.
We will then  prove that for any $h\in\Hhp$, $t\in -1:\tirr$ and $k\in\rel^c$, 
$$\1_{\Ehp} \vert Z_{h^{(t)}k }\vert \leq \lambda_{h^{(t)}k}.$$
Let us fix $h\in\Hhp$, $t\in -1:\tirr$ and $k\in\rel^c$.
We decompose $Z_{h^{(t)}k}$ as follows:
\begin{equation}
Z_{h^{(t)}k}= \left(Z_{h^{(t)}k}-\bar{Z}_{h^{(t)}k}\right) + \left(\bar{Z}_{h^{(t)}k} - \E\bar{Z}_{h^{(t)}k}\right) + \E\bar{Z}_{h^{(t)}k}.
\end{equation}
We use:\begin{itemize}
\item 1. of \Cref{lmDelta}: Recall the notation $\DZ{h^{(t)}k}:=Z_{h^{(t)}k}-\bar{Z}_{h^{(t)}k}$, then remark that $\forall h'\in\HhpRev\cup\HhpDir$, $\prod_{k=1}^d h'_k\leq 1$, and 
%$$\hirr^{d/2}=\frac{\delta\Cl}{8M_X \CEZA}.$$ 
$\Ehp\subset \BZA{h^{(t)},k}\cap\Atilde$, therefore:
$$\1_{\Ehp}\left(Z_{h^{(t)}k}-\bar{Z}_{h^{(t)}k}\right)\leq \frac14 \seuil{h^{(t)}k} ,$$
\item the definition of $\BZ{h^{(t)},k}$: since $\Ehp\subset \BZ{h^{(t)},k}$,
$$\1_{\Ehp} \left\vert\bar{Z}_{h^{(t)}k} - \E\bar{Z}_{h^{(t)}k}\right\vert \leq \frac12 \seuil{h^{(t)}k},$$
\item 2. of \Cref{lmZbar}: since $k\in\rel^c$,
$$\E\bar{Z}_{h^{(t)}k}=0.$$
\end{itemize}
Therefore:
$$\1_{\Ehp} \vert Z_{h^{(t)}k}\vert \leq \frac34 \seuil{h^{(t)}k}\leq \seuil{h^{(t)}k},$$
and so, every irrelevant component is active during Reverse Step until Iteration $\tirr$. In particular, we have proved that: 
$$\Ehp\subset \lbrace\forall k\in\rel^c : \hat{h}_k = \hirr\rbrace.$$

Let us now prove that on $\Ehp$, 
$$\prod\limits_{k=1}^d\hat{h}_k\geq \beta^r \frac{(\log n)^{1+a}}{n}.$$
The loop condition in the Direct Step imposes that at the beginning of any iteration $t\geq 0$:
$$\prod\limits_{k=1}^d\hat{h}^{(t)}_k\geq \frac{(\log n)^{1+a}}{n}.$$
For our algorithm, the bandwidth can only decrease during the Direct Step. Since on $\Ehp$, the irrelevant components are active the during Reverse Step, they are inactive during the Direct Step. This is the reason why during the last iteration, only relevant components could decrease and be multiplied by $\beta$. Therefore:
$$\prod\limits_{k=1}^d\hat{h}_k\geq \beta^r \frac{(\log n)^{1+a}}{n},$$
which ends the proof of the inclusion \eqref{eProp2hinHhp} of \Cref{propHhp}.

Finally, we control $\P\left( \Ehp^c\right)$.
We first control the cardinal of $\Hhp$ by enumerating the possible values for a component of a bandwidth in $\Hhp$. For $h\in\Hhp$ and $k\in\rel$, 
$$\beta (\log n)^{1+a} n^{-1}  \leq h_k \leq 1,$$ 
 thus:
\begin{equation*}
\left\vert \left\lbrace  h_k: h\in\Hhp \right\rbrace\right\vert
=\left\vert \left\lbrace  \beta^t h_0\in[\beta (\log n)^{1+a} n^{-1},1], t\in\mathbb{Z}\right\rbrace\right\vert
\leq 1+\logb\left(\frac{1}{\beta (\log n)^{1+a} n^{-1} }\right)
\leq \logb n
\end{equation*} 
(for $n$ large enough). For $k\in\rel^c$, 
$$ h_k=\hirr,$$ thus, we have
\begin{equation*}
\left\vert \left\lbrace  h_k: h\in\Hhp \right\rbrace\right\vert
= 1.
\end{equation*} 
Therefore:
\begin{equation}
\left\vert \Hhp\right\vert 
\leq  \left(\logb n\right)^r.
\end{equation}
Let us also control the cardinal of $\HhpRev\cup\HhpDir$.
The only supplementary bandwidths are the ones whose irrelevant components are smaller than $\hirr$. We consider the irrelevant components as the relevant ones, and we obtain the rough bound
\begin{equation}
\left\vert \HhpRev\cup\HhpDir\right\vert 
\leq  \left(\logb n\right)^d.
\end{equation}
By Assumption \AestfX, \ref{fXtildeAccuracy}: 
$$\P\left(\Atilde^c \right) \leq \exp(-(\log n)^{1+\frac{a-1}2}).$$
\\ 
We bound the events $\Bernf{h}^c$'s and $\Bernfabs{h}^c$'s using \Cref{lmfbar}. Since for all $h\in\Hhp$, 
$$\prod\limits_{k=1}^d h_k\geq \beta^r \tfrac{(\log n)^{a+1}}{n},$$ 
note that:
\begin{itemize}
\item Cond($h$): $\prod\limits_{k=1}^d h_k\geq \frac{4\Vert K \Vert^{2d}_{\infty}}{9\delta^2\Cs^2} \frac{(\log n)^a}{n}$ 
is satisfied for any $h\in\Hhp$ for $n$ large enough (when $\log n\geq \frac{4\Vert K \Vert^{2d}_{\infty}}{9\beta^r\delta^2\Cs^2}$). So, we have
$$\mathbb{P}\left( \Bernf{h}^c \right) 
\leq 2 e^{-(\log n)^a}.$$
\item Moreover, 
$$\mathbb{P}\left( \Bernfabs{h}^c \right) 
\leq  2 e^{-\CgfA n\prod_{k=1}^d h_k} \leq 2 e^{-\CgfA \beta^r (\log n)^{a+1}}.$$
\end{itemize}
Similarly, we bound the probability of events $\BernZ{h}^c$'s and $\BernZabs{h}^c$'s using \Cref{lmZbar}.
Note that for all $h\in\HhpRev\cup\HhpDir$:
\begin{itemize}
\item $\text{Cond}_{\bar{Z}}(h)$: $ \prod\limits_{k=1}^d h_k
\geq   \condZbar\frac{(\log n)^a}{n}$ is satisfied for $n$ large enough (when $\log n\geq  \frac{\condZbar}{\beta^r}$). So, we have
$$\PP{\BZ{h,j}^c}
\leq 2 e^{-\frac{\delta }{\Vert f\Vert_{\infty,\UUn}} (\log n)^a}.$$
\item Moreover, 
$$
\PP{\BZA{h,j}^c}
\leq 2 e^{-\CgZA n\prod_{k=1}^d h_k}
\leq 2 e^{-\CgZA \beta^r (\log n)^{a+1}}.
$$
\end{itemize}
Therefore,
\begin{align*}
\P\left( \Ehp^c\right)
&\leq \P\left( \Atilde ^c\right)
+ \sum\limits_{h\in\Hhp} \left( \P\left(\Bernf{h}^c\right) + \P\left(\Bernfabs{h}^c \right) \right)\\
&\quad +  \sum\limits_{h\in (\HhpRev\cup\HhpDir)} \sum\limits_{k=1}^d \left( \P\left(\BernZ{h,k}^c \right) + \P\left(\BernZabs{h,k}^c \right) \right)
\\
&\leq e^{-(\log n)^{1+\frac{a-1}2}}
+ \sum\limits_{h\in\Hhp} \left( 2 e^{-(\log n)^a} + 2 e^{-\CgfA \beta^r (\log n)^{a+1}}
\right)\\
&\quad +  \sum\limits_{h\in (\HhpRev\cup\HhpDir)} \sum\limits_{k=1}^d \left( 2 e^{-\frac{\delta }{\Vert f\Vert_{\infty,\UUn}} (\log n)^a} + 2 e^{-\CgZA \beta^r (\log n)^{a+1}}
\right)\\
&\leq  e^{-(\log n)^{1+\frac{a-1}2}} 
\left(1
+  4 \left(\logb n\right)^r e^{-(\log n)^{\frac{a-1}2}} 
+ 4 d \left(\logb n\right)^d e^{- \frac{\delta }{\Vert f\Vert_{\infty,\UUn}} (\log n)^{\frac{a-1}2}}
\right)\\
&\leq 2 e^{-(\log n)^{1+\frac{a-1}2}},
\end{align*}
for $n$ large enough.
%%%%%%%%%%%%%%%%%%%%%%%%%%%%%%
\subsection{Proof of \Cref{propBiasDeviation}}
We fix $h\in\Hhp$ and consider the event $\lbrace\hat{h}=h\rbrace\cap \Ehp$.
Let $(t_1,\dots,t_d)\in\mathbb{Z}^d$ such that for all $k \in 1:d$,
$$h_k=\beta^{t_k} h_0.$$
Given positive constants $\Am$ and $\CAm$ (to be opimized), we call $ \CAm \left(\log n\right)^A n^{-\frac{1}{2s+r}}$ the minimax bandwidth level and we define $\tminimax\in \mathbb{R}$ such that
$$\beta^{\tminimax} h_0 = \CAm \left(\log n\right)^A n^{-\frac{1}{2s+r}}.$$
Using the definition \eqref{ineq-initialisation} of $h_0$, observe that $\tminimax>0$ (for $n$ large enough). 
To simplify the notation (permutation of the labels), we consider: 
$$\rel=1:r$$
and  \begin{equation}\label{ordretk}
t_1\geq t_2\geq \dots \geq t_r.
\end{equation}

\subsubsection{Proof of Inequality~\eqref{eMajBias}}\label{sec:eMajBias}
The bias of $\fbar{h}$ is denoted $\biasbar{h}$. Note that it does not depend on $\lbrace h_k\rbrace_{k\in\rel^c}$. Indeed, we have
\begin{align}
\biasbar{h}:&= \E\left[\fbar{h}\right] - f(w) \nonumber\\
&=\int_{u\in\R^d}\K_h(w-u) \frac{{\rm f}_W(u)}{\fX(u_{1:d_1})} du- f(w)\nonumber\\
&=\int_{u\in\R^d}\K_h(w-u)f(u) du- f(w)\nonumber\\
%&=\int_{z\in\R^d}\left(\prod\limits_{k=1}^d K(z_k)\right) f(w-h\cdot z)dz- f(w)\int_{z\in\R^d}\left(\prod\limits_{k=1}^d K(z_k)\right)dz\nonumber\\
&=\int_{z\in\R^d}\left(\prod\limits_{k=1}^d K(z_k)\right) \left[f(w-h\cdot z)- f(w)\right]dz\label{eBiaisNul}\\
&=\int_{z'\in\R^r} \left(\prod\limits_{k=1}^r K(z'_{k})\right) \left[f_{\rel}\left(w_{1:r}-h_{1:r}\cdot z' \right)-f_\rel(w_{1:r})\right] dz'. \nonumber
\end{align}
We consider the following disjunction of cases:
\begin{enumerate}[label=(Case \Alph*)]
\item \label{casA} without relevant component: $\rel=\emptyset$ 
\item \label{casB} with small relevant bandwidth components: $\min\limits_{j\in\rel} t_j\geq \tminimax$
\item \label{casC} with at least one large relevant bandwidth component: $\exists j\in\rel, t_j < \tminimax$.
\end{enumerate}
Then we control the bias in each case.
\begin{enumerate}[label=(Case \Alph*)]
\item {Assume $\rel=\emptyset$. In particular, $f$ is constant on the neighborhood $\mathcal{U}$. Note that for any $z\in\text{supp}\left(K\right)^d$, $w-h\cdot z\in\mathcal{U}$.
We then derive from Equation \eqref{eBiaisNul}:
$$\biasbar{h}=0.$$
}
\item {Assume $\min\limits_{j\in\rel} t_j\geq \tminimax$. We apply 2. of \Cref{lmfbar}
\begin{align*}
\left\vert\biasbar{h} \right\vert
&\leq \CBbar \sum\limits_{j\in\rel} h_j^s
=\CBbar \sum\limits_{j\in\rel} \left(\beta^{t_j}h_0\right)^s\\
&\leq \CBbar\times r \left(\beta^{\tminimax}h_0\right)^s
= r \CBbar  {\CAm}^s \left(\log n\right)^{\Am s} n^{-\frac{s}{2s+r}}
\end{align*}
}
\item {Assume $\exists j\in\rel, t_j < \tminimax$. Then we consider 
$$\jm=\min \left(j\in\rel : t_j <\tminimax\right).$$
In particular, for all $j\geq \jm$, the bandwidth components are larger than the minimax level:
\begin{equation}
h_j\geq \CAm (\log n)^{\Am} n^{-\frac1{2s+r}}.
\end{equation}
For the previously fixed bandwidth $h$ (and its relevant deactivation times $(t_1,\dots,t_r)$), we define the following intermediate bandwidths $\hint{t}$, $t\in\mathbb{R}$:
$$\hint{t}_k=\left\lbrace 
\begin{tabular}{@{}ll}
		$\beta^{t\vee t_k} h_0$ &$\text{ if }k\in\rel$  \\ 
		$h_k$ &$\text{ else.}$\\
		\end{tabular}\right. $$

Then we decompose the bias by splitting $f(w-h\cdot z)-f(w)$ (note that $\hint{t_r}=h$):
\begin{align}
\biasbar{h}
= \int_{z\in\R^d}\left(\prod\limits_{k=1}^d K(z_k)\right) 
&{\large[}f(w-\hint{\tminimax}\cdot z)-f(w) \nonumber\\
&+ f(w-\hint{t_{\jm}}\cdot z)-f(w-\hint{\tminimax}\cdot z)\nonumber\\
&+\sum\limits_{j_0=\jm+1}^r  f(w-\hint{t_{j_0}}\cdot z)-f(w-\hint{t_{j_0-1}}\cdot z)
{\large]}dz\nonumber\\
=\biasbar{\hint{\tminimax}} 
+ (\biasbar{\hint{t_{\jm}}} &-\biasbar{\hint{\tminimax}})
+\sum\limits_{j_0=\jm+1}^r \left(\biasbar{\hint{t_{j_0}}} -\biasbar{\hint{t_{j_0-1}}} \right).
\label{eDecompBiashint}
\end{align}
For the first term, note that $\hint{\tminimax}$ satisfies the condition of \ref{casB}, thus:
\begin{equation}\label{eMajBhminimax}
\left\vert\biasbar{\hint{\tminimax}} \right\vert
\leq r \CBbar {\CAm}^s \left(\log n\right)^{\Am s} n^{-\frac{s}{2s+r}}.
\end{equation}
Let us now control the other terms. The same arguments are used to control the terms in the sum $\biasbar{\hint{t_{j_0}}} -\biasbar{\hint{t_{j_0-1}}}$ (for $j_0 \in (\jm+1):r$) and the second term $\biasbar{\hint{t_{\jm}}} -\biasbar{\hint{\tminimax}}$. To shorten the proof, the followings lines are also applied to control the second term: for the added case $j_0=j_\Am$, one just has to replace $\hint{t_{{j_0}-1}}_{j}$ by $\hint{\tminimax}$.\\ 
Let us now fix $j_0\in \jm : r$ and consider the path between $\hint{t_{{j_0}-1}}_{j}$ and $\hint{t_{j_0}}_{j}$. Namely for $u\in[0,1]$, we denote $\hu{j_0}:=\hint{t_{j_0-1}}+ u\left(\hint{t_{j_0}}-\hint{t_{j_0-1}}\right)$.
Remark that, for any $j \in 1:d$,
\begin{align*}
\hint{t_{j_0}}_{j}-\hint{t_{j_0-1}}_{j}\neq 0&\Rightarrow \left(j\in\rel \text{ and }\beta^{t_{j_0}\vee t_j}\neq \beta^{(t_{j_0-1})\vee t_j}\right)\\
&\Rightarrow \left(j\in\rel \text{ and }{t_j< t_{j_0} \mbox{ or }t_j< t_{j_0-1}}\right)\\
&\Rightarrow \left(j\in\rel \text{ and }t_j\leq t_{j_0}\right).
%&\Rightarrow \left(j\in\rel \text{ and }j\geq j_0\right)
\end{align*}
The last implication is due to the fact that a component could not be deactived between the consecutive deactivation times  $t_{j_0}$ and $t_{j_0-1}$.
%Indeed, given the definition of $\hint{t}$ for all $t$, each irrelevant component $j$ keeps the value $h_j$.  For $j\in\rel$, note that $\beta^{t_j\vee t_{j_0}}\neq \beta^{t_j\vee t_{j_0-1}} \textcolor{magenta}{\Rightarrow t_j< t_{j_0} \mbox{ or }t_j< t_{j_0-1}\Rightarrow t_j\leq t_{j_0}.}$\\

Then, we introduce the function $g:u\in[0,1]\mapsto f(w-\hu{j_0}\cdot z)$ (for a fixed $z\in\R^d$). In particular, using the above remark:
\begin{align*}
g'(u)
&= \sum\limits_{\substack{j\in\rel\\t_{j}\leq t_{j_0}}}  \left( \hint{t_{j_0}}_j-\hint{t_{j_0-1}}_j\right)\times z_{j} \partial_{j}f(w-\hu{j_0}\cdot z).
\end{align*}
%\red{ici on utilise $s>1$}\\
Then we write:
\begin{align*}
f(w-\hint{t_{j_0}}\cdot z)-&f(w-\hint{t_{j_0-1}}\cdot z)\\
&= g(1)- g(0)
=\int_{u=0}^1 g'(u) du\\
&=\sum_{\substack{j\in\rel\\t_j\leq t_{j_0}}} \int_{u=0}^1 \left( \hint{t_{j_0}}_j-\hint{t_{j_0-1}}_j\right)\times z_{j}  \partial_{j}f(w-\hu{j_0}\cdot z)du.
\end{align*}
Hence, we obtain
\begin{align}
\biasbar{\hint{t_{j_0}}} -\biasbar{\hint{t_{j_0-1}}}
&=\int_{z\in\R^d}\left(\prod\limits_{k=1}^d K(z_k)\right) 
 [f(w-\hint{t_{j_0}}\cdot z)-f(w-\hint{t_{j_0-1}}\cdot z)]
dz\nonumber\\
&= \sum_{\substack{j\in\rel\\t_j\leq t_{j_0}}} \int_{u=0}^1 \left( \hint{t_{j_0}}_j-\hint{t_{j_0-1}}_j\right)\int_{z\in\R^d} \left(\prod\limits_{k=1}^d K(z_k)\right) z_{j}  \partial_{j}f(w-\hu{j_0}\cdot z)dz\ du
\nonumber\\
&=\sum_{\substack{j\in\rel\\t_j\leq t_{j_0}}} \int_{u=0}^1 \left( \hint{t_{j_0}}_j-\hint{t_{j_0-1}}_j\right) \E\left[\bar{Z}_{\hu{j_0},j}\right] du,\label{eDiffBiais}
\end{align}
using Equation \eqref{eEZ}: 
\begin{align*}
\E\left[\bar{Z}_{\hu{j_0},j}\right]
=\int_{\R^d}  \left(\prod\limits_{k=1}^d K(z_k)\right) z_j \partial_j f(w-\hu{j_0}\cdot z)dz.
\end{align*}
Now the idea is to control $\left\vert \E\left[\bar{Z}_{\hu{j_0},j}\right]\right\vert$ with the test at the iteration $t_j$ on $\vert Z_{h^{(t_j)},j}\vert$. More precisely, we will first apply Assumption~\Monoton{} to move from $\left\vert \E\left[\bar{Z}_{\hu{j_0},j}\right]\right\vert$ to $\left\vert \E\left[\bar{Z}_{h^{(t_j)},j}\right]\right\vert$. Then, we will apply Bernstein's inequality to convert the control on $\left\vert Z_{h^{(t_j)},j}\right\vert$ to a control on $\left\vert \E\left[\bar{Z}_{h^{(t_j)},j}\right]\right\vert$.\\
Let us fix $j\in\rel$ such that $t_j\leq t_{j_0}$. We distinguish the cases where the component $j$ is deactivated during the Reverse Step or when it happens during the Direct Step.
\begin{enumerate}[label=Subcase (C.\alph*)]
\item{$t_j\geq 0$, \textit{i.e.}: $j$ is deactivated during the Direct Step.\\
Let us show $\hu{j_0}\preccurlyeq h^{(t_j)}$:  
\begin{itemize}
\item for $k\in\rel^c$, since $\hint{t_{j_0-1}}_k=h_k=\hint{t_{j_0}}_k$,
\begin{equation*}
\hu{j_0}_k=h_k.
\end{equation*}
Remember that the irrelevant components deactivate during the Reverse Step, therefore they already have their final value during the Direct Step. Formally, since $t_k<0\leq t_j$, we have
\begin{align*}
\hu{j_0}_k=h_k=\beta^{t_k} h_0=\beta^{t_j\wedge t_k}h_0=h^{(t_j)}_k.
\end{align*}
\item for $k\in\rel$, notice $\hint{t_{j_0-1}}\preccurlyeq \hint{t_{j_0}}$. Therefore:
\begin{align*}
\hu{j_0}_k
&\leq  \hint{t_{j_0}}_k 
= \beta^{t_{j_0}\vee t_{k}} h_0\\
&\leq \beta^{t_j\wedge t_{k}} h_0=h^{(t_j)}_k.
\end{align*}
\end{itemize}
Then, we have proved $\hu{j_0}\preccurlyeq h^{(t_j)}$. Using Assumption~\Monoton{}:
$$\left\vert \E\left[\bar{Z}_{\hu{j_0},j}\right]\right\vert 
\leq \left\vert \E\left[\bar{Z}_{h^{(t_j)},j}\right]\right\vert.$$
 }
 \item{ $t_j<0$, \textit{i.e.}: $j$ is deactivated during Reverse Step.\\
 As well as $h'\mapsto\biasbar{h'}$, $h'\mapsto\E\left[\bar{Z}_{h',j}\right]$ is independent of the irrelevant components of the bandwidth (see for instance Equation \eqref{eEZ}).\\
 Then we modify the irrelevant components of $\hu{j_0}$ and use the value of the irrelevant components of $h^{(t_j)}$. Formally, we introduce the notation $\huvar{j_0}$ such that
 $$\huvar{j_0}_k=\left\lbrace 
\begin{tabular}{@{}ll}
		$\hu{j_0}_{k}$ &$\text{ if }k\in\rel$  \\ 
		$h^{(t_j)}_k$ &$\text{ else,}$\\
		\end{tabular}\right. $$
so that:
$$\E\left[\bar{Z}_{\hu{j_0},j}\right] =\E\left[\bar{Z}_{\huvar{j_0},j}\right].$$
Now we just have to verify $\huvar{j_0}\preccurlyeq h^{(t_j)}$:
\begin{itemize}
\item for $k\in\rel^c$, by definition of $\huvar{j_0}$:
$$\huvar{j_0}_k=h^{(t_j)}_k$$
\item for $k\in\rel$,
\begin{align*}
\huvar{j_0}_k
&=\hu{j_0}_k\\
&\leq \hint{t_{j_0}}_k
=\beta^{t_{j_0}\vee t_k} h_0\\
&\leq \beta^{t_{j}\vee t_k} h_0, \text{ since }t_j\leq t_{j_0},\\
&\leq\max\left(h_k, \beta^{t_j}h_0\right)=:h^{(t_j)}_k. 
\end{align*}
\end{itemize}
Then we have proved $\huvar{j_0}\preccurlyeq h^{(t_j)}$. Using Assumption~\Monoton{}:
$$\left\vert \E\left[\bar{Z}_{\hu{j_0},j}\right]\right\vert 
=\left\vert \E\left[\bar{Z}_{\huvar{j_0},j}\right]\right\vert
\leq \left\vert \E\left[\bar{Z}_{h^{(t_j)},j}\right]\right\vert. $$
 }
\end{enumerate}
In each case (C.a and C.b), we have proved $\left\vert \E\left[\bar{Z}_{\hu{j_0},j}\right]\right\vert 
\leq \left\vert \E\left[\bar{Z}_{h^{(t_j)},j}\right]\right\vert$, then we apply this inequality in Equation \eqref{eDiffBiais}:
\begin{align}
\left\vert \biasbar{\hint{t_{j_0}}} -\biasbar{\hint{t_{j_0-1}}}\right\vert
&\leq \sum_{\substack{j\in\rel\\t_j\leq t_{j_0}}} \int_{u=0}^1 \left( \hint{t_{j_0}}_j-\hint{t_{j_0-1}}_j\right) \left\vert\E\left[\bar{Z}_{\hu{j_0},j}\right]\right\vert du\label{poursup1}\\
&\leq \sum_{\substack{j\in\rel\\t_j\leq t_{j_0}}}  \int_{u=0}^1\left( \hint{t_{j_0}}_j-\hint{t_{j_0-1}}_j\right) \left\vert \E\left[\bar{Z}_{h^{(t_j)},j}\right]\right\vert du\nonumber\\
&\leq \sum_{\substack{j\in\rel\\t_j\leq t_{j_0}}} \left( \hint{t_{j_0}}_j-\hint{t_{j_0-1}}_j\right) \left\vert \E\left[\bar{Z}_{h^{(t_j)},j}\right]\right\vert.\label{poursup2}
\end{align}
Then, the previous decomposition of the bias \eqref{eDecompBiashint} leads to:
\begin{align}
\left\vert \biasbar{h} \right\vert
&\leq \left\vert \biasbar{\hint{\tminimax}} \right\vert
+\sum\limits_{j_0=\jm}^r \left\vert \biasbar{\hint{t_{j_0}}} -\biasbar{\hint{t_{j_0-1}}} \right\vert\nonumber \\
&\leq 
r\CBbar  {\CAm}^s \left(\log n\right)^{\Am s} n^{-\frac{s}{2s+r}}
+ \sum\limits_{j_0=\jm}^r  \sum_{\substack{j\in\rel\\t_j\leq t_{j_0}}}  \left( \hint{t_{j_0}}_j-\hint{t_{j_0-1}}_j\right) \left\vert \E\left[\bar{Z}_{h^{(t_j)},j}\right]\right\vert \nonumber \\
&\leq 
r\CBbar{\CAm}^s \left(\log n\right)^{\Am s} n^{-\frac{s}{2s+r}}
+ \sum\limits_{j=\jm}^r  \left\vert \E\left[\bar{Z}_{h^{(t_j)},j}\right]\right\vert 
\sum\limits_{j_0=\jm}^{j}  \left( \hint{t_{j_0}}_j-\hint{t_{j_0-1}}_j\right)\nonumber \\
&\leq 
r\CBbar  {\CAm}^s \left(\log n\right)^{\Am s} n^{-\frac{s}{2s+r}}
+ \sum\limits_{j=\jm}^r  \left\vert \E\left[\bar{Z}_{h^{(t_j)},j}\right]\right\vert 
 h^{(t_j)}_j,\label{biaispoursupp}
\end{align}
 since the sum is telescoping, and by noticing that:  $\hint{t_{j}}_j=h^{(t_j)}_j$.\\
 
 Now, it remains to control $\left\vert \E\left[\bar{Z}_{h^{(t_j)},j}\right]\right\vert$ for $j \in \jm:r$ using the test at the iteration $t_j$ on $Z_{h^{(t_j)},j}$:
 \begin{align*}
 \1_{\Ehp\cap \lbrace\hat{h}=h\rbrace} \left\vert \E\left[\bar{Z}_{h^{(t_j)},j}\right]\right\vert 
 &\leq 
 \1_{\hat{h}=h}\left\vert Z_{h^{(t_j)},j}\right\vert 
 + \1_{\Atilde \cap \BernZabs{h^{(t_j)},j}} \left\vert Z_{h^{(t_j)},j}-\bar{Z}_{h^{(t_j)},j}\right\vert \\
 &+\1_{\BernZ{h^{(t_j)},j}} \left\vert \bar{Z}_{h^{(t_j)},j}-\E\left[\bar{Z}_{h^{(t_j)},j}\right]\right\vert
 \end{align*}
By construction of the \CDRodeo{} procedure, if $\hat{h}=h$, then $j$ is deactivated at iteration $t_j$, in other words:
 $$ \1_{\Ehp\cap \lbrace\hat{h}=h\rbrace} \left\vert Z_{h^{(t_j)},j}\right\vert \leq \seuil{h^{(t_j)},j}.$$
We also apply:
 \begin{itemize}
 \item the definition of $\BernZ{h^{(t_j)},j}$:
 $$\1_{\BernZ{h^{(t_j)},j}} \left\vert \bar{Z}_{h^{(t_j)},j}-\E\left[\bar{Z}_{h^{(t_j)},j}\right]\right\vert 
 \leq \frac12 \seuil{h^{(t_j)},j},$$
 %Note that the condition $\text{Cond}_{\bar{Z}}(h^{(t_j)})$ is satisfied for $n$ large enough: since $h^{(t_j)}\in\Hhp$,
 %$$\prod\limits_{k=1}^d h_k^{(t_j)}\geq \beta^r \frac{(\log n)^{a+1}}n
%\gg \condZbar\frac{(\log n)^a}{n}$$
 \item 1. of \Cref{lmDelta} (note in particular $\prod_{k=1}^d h^{(t_j)}_k\leq 1$):
 $$\1_{\Atilde \cap \BernZabs{h^{(t_j)},j}} \left\vert Z_{h^{(t_j)},j}-\bar{Z}_{h^{(t_j)},j}\right\vert 
 =\1_{\Atilde \cap \BernZabs{h^{(t_j)},j}} \left\vert \DZ{h^{(t_j)}j} \right\vert 
 \leq \frac14 \seuil{h^{(t_j)},j}.$$
 \end{itemize}
Therefore:
\begin{align*}
 \1_{\Ehp\cap \lbrace\hat{h}=h\rbrace} \left\vert \E\left[\bar{Z}_{h^{(t_j)},j}\right]\right\vert 
 \leq \1_{\Ehp\cap \lbrace\hat{h}=h\rbrace} \frac74 \seuil{h^{(t_j)},j}.
 \end{align*}
 Hence:
 \begin{align}
\1_{\Ehp\cap \lbrace\hat{h}=h\rbrace} \left\vert \biasbar{h} \right\vert
&\leq \1_{ \lbrace\hat{h}=h\rbrace} \left (r \CBbar {\CAm}^s \left(\log n\right)^{\Am s} n^{-\frac{s}{2s+r}}
+ \sum\limits_{j=\jm}^r  \frac74 \seuil{h^{(t_j)},j}\times
 h^{(t_j)}_j\right),
 \nonumber\\
 &\leq \1_{ \lbrace\hat{h}=h\rbrace} \left( r\CBbar {\CAm}^s \left(\log n\right)^{\Am s} n^{-\frac{s}{2s+r}}
+ \sum\limits_{j=\jm}^r  \frac{7\Cl (\log n)^{a/2}}{4\left(n\prod_{k=1}^d h_k^{(t_j)}\right)^{1/2}} 
\right).\label{eMajBiaisAvantDisjonctionDeCas}
\end{align}
Then we control $\prod_{k=1}^d h_k^{(t_j)}$ using the same disjunction of subcases as above:
\begin{enumerate}[label=Subcase (C.\alph*)]
\item{$t_j\geq 0$. At the iteration $t_j\geq 0$, the Direct Step has begun, thus the Reverse Step is over. Since $h\in\Hhp$, the irrelevant components have already their final value: for all $k\in\rel^c$,
$$1
\geq h_k^{(t_j)}
=h_k
=\hirr
>\beta.$$
Moreover, during the Direct Step, at iteration $t_j$, all components are lower bounded by the current active bandwidth value $\beta^{t_j} h_0$, \textit{i.e.}: for any $k\in\rel$,
$$h_k^{(t_j)}\geq \beta^{t_j} h_0.$$
Recall that $j\geq \jm$, thus: $$t_{j}\leq t_{\jm}\leq \tminimax.$$
It follows:
 $$h_k^{(t_j)}\geq \beta^{\tminimax} h_0=\CAm \left(\log n\right)^A n^{-\frac{1}{2s+r}}.$$ 
 Therefore:
 $$ \prod_{k=1}^d h_k^{(t_j)} \geq \beta^{d-r} \left( \CAm \left(\log n\right)^A n^{-\frac{1}{2s+r}} \right)^r.$$
 Then the upper bound in Equation \eqref{eMajBiaisAvantDisjonctionDeCas} becomes:
 \begin{align*}
 \frac{7\Cl (\log n)^{a/2}}{4\left(n\prod_{k=1}^d h_k^{(t_j)}\right)^{1/2}} 
&\leq \tfrac{7\Cl }{4\beta^{\frac{d-r}2}{\CAm}^{\frac{r}2}} (\log n)^{\frac{a-Ar}2}n^{-\frac12\left(1-\frac{r}{2s+r}\right)}\\
&=\tfrac{7\Cl }{4\beta^{\frac{d-r}2}{\CAm}^{\frac{r}2}} (\log n)^{\frac{a-Ar}2}n^{-\frac{s}{2s+r}}.
 \end{align*}
}
\item{$t_j< 0$. At iteration $t_j$, only iterations of the Reverse Step have been performed. Thus, the current bandwidth has only been increased. Therefore:
$$ 
\frac{7\Cl (\log n)^{a/2}}{4\left(n\prod_{k=1}^d h_k^{(t_j)}\right)^{1/2}} 
\leq \frac{7\Cl (\log n)^{a/2}}{4\left(n h_0^d\right)^{1/2}} .
$$
Remark that the lower bound on $h_0$ \eqref{ineq-initialisation} is exactly defined so, we have
$$
\frac{7\Cl (\log n)^{a/2}}{4\left(n h_0^d\right)^{1/2}} 
\leq \frac{7}{4}\left(\frac{(\log n) ^{a}}{n}\right)^{\frac{p}{2p+1}}.
$$
Note that $n^{-\frac{p}{2p+1}}$ is smaller than the minimax optimal rate for any regularity and any sparsity structure (except for the degenerate case where $r=0$ and which is solved separately: cf (Case A)):
$$n^{-\frac{p}{2p+1}}=\min\limits_{\substack{1\leq r'\leq d\\ 1\leq s'\leq p}}\left( n^{-\frac{s'}{2s'+r'}}\right).$$
%\red{Prendre $C_p=1$(plus simple) ou $\CAm$ optimal pour $s=p$ et $r=1$ (si le ratio $\frac{C_{opt,2}}{C_{opt,1}}$ se simplifie en constantes ne dépendant pas de $f$ (à $r=1$ et $s=p$ connu)\\
%Claire : prendre $C_p=2$ car simple et $\geq 7/4$\\
%Vincent : au log pres, prendre $h_0^d\geq n^{-1/(2p+1)}$}
}
\end{enumerate}
When we reunite the two subcases, Inequality \eqref{eMajBiaisAvantDisjonctionDeCas} becomes:
\begin{align*}
\1_{\Ehp\cap \lbrace\hat{h}=h\rbrace} \left\vert \biasbar{h} \right\vert
&\leq  r \CBbar{\CAm}^s \left(\log n\right)^{\Am s} n^{-\frac{s}{2s+r}}
\\
&+ r\times \max\left( \tfrac{7\Cl }{4\beta^{\frac{d-r}2}{\CAm}^{\frac{r}2}} \frac{(\log n)^{\frac{a-Ar}2}}{n^{\frac{s}{2s+r}}}
, {\rcp} \left(\frac{(\log n )^{a}}{n}\right)^{\frac{p}{2p+1}}\right),
\end{align*}
}
\end{enumerate}
which concludes the proof of Inequality \eqref{eMajBias} %(in particular, note that the upper bounds in (Case A) and (Case B) can be rewrite into this last inequality).\\

\subsubsection{Proof of Inequality~\eqref{eMajDeviation}}
Let us now prove the second inequality \eqref{eMajDeviation}.
By definition: $\Ehp\subset\Bernf{h}$. Thus, we have
$$ \1_{\lbrace\hat{h}=h\rbrace\cap \Ehp} \left\vert \fbar{h}- \E\left[\fbar{h}\right]\right\vert 
\leq \sigma_h
:=\Cs \sqrt{\frac{(\log n)^a}{n\prod_{k=1}^d h_k}}.$$
Two cases occur: in the first case, the deviation is controlled  by a concentration inequality; in the second case, we control the deviation by $\E Z_{hj}$ thanks to the tests on the $Z_{hj}$'s.
\begin{enumerate}
\item {$\max\limits_{k\in\rel} t_k 
\leq \tminimax$.
Then, $\forall k\in\rel$:
 $$h_k=\beta^{t_k} h_0 
 > \beta^{\tminimax} h_0 = \CAm (\log n)^{\Am} n^{-\frac{1}{2s+r}}.
$$
Besides, for $k\in\rel^c$:
$$h_k= \hirr>\beta.$$
Therefore:
\begin{align*}
\sigma_h
\leq \Cs \sqrt{\frac{(\log n)^a}{n {\beta}^{d-r} \left(\CAm (\log n)^{\Am} n^{-\frac{1}{2s+r}} \right)^r}}
=\tfrac{\Cs}{{\beta}^{(d-r)/2}{\CAm}^{r/2}} (\log n)^{(a-\Am r)/2} n^{-\frac{s}{2s+r}}.
\end{align*}

}
\item {$\max\limits_{k\in\rel} t_k > \tminimax$. First remark that for any $k \in 1:d$,
$$\sigma_h=\frac{\Cs}{\Cl} h_{k}\;\seuil{hk}.$$
Hence, it suffices to control the threshold in order to bound the deviation.
 Let us consider $j_0\in\argmax_{k\in\rel} t_k$ (actually assuming \eqref{ordretk} means that $j_0=1$).
 In particular, when $\hat{h}=h$, the component $j_0$ is deactivated during the last iteration, and during the Direct Step (recall that $\tminimax>0$). Let us consider the penultimate iteration, i.e. Iteration $t_{j_0}-1$. At this iteration, $j_0$ is not deactivated, \textit{i.e.}:
$$ \1_{\hat{h}=h} \left\vert Z_{h^{(t_{j_0}-1)}j_0}\right\vert 
>\1_{\hat{h}=h} \seuil{h^{(t_{j_0}-1)}j_0}.$$
Then we use 1. of \Cref{lmDelta}. Note that $\prod_{k=1}^d h^{(t_{j_0}-1)}_k\leq 1$, thus:
$$ \1_{\Ehp}\left\vert\DZ{h^{(t_{j_0}-1)}j_0}\right\vert 
\leq  \frac{1}{4} \seuil{h^{(t_{j_0}-1)}j_0}.$$
Remember the definition of $\BernZ{h,j}$, thus
$$ \1_{\Ehp}\left\vert \bar{Z}_{h^{(t_{j_0}-1)}j_0} - \E\left[  \bar{Z}_{h^{(t_{j_0}-1)}j_0}\right]\right\vert 
\leq  \frac12 \seuil{h^{(t_{j_0}-1)}j_0}.$$
Therefore:
\begin{equation}
\1_{\lbrace\hat{h}=h\rbrace\cap\Ehp} \left\vert \E\left[  \bar{Z}_{h^{(t_{j_0}-1)}j_0}\right]\right\vert 
>\1_{\lbrace\hat{h}=h\rbrace\cap\Ehp}\ \frac14 \seuil{h^{(t_{j_0}-1)}j_0}.
\end{equation}
Let us compare $h^{(t_{j_0}-1)} $ to $h$. Recall $h=h^{(t_{j_0})}$, since $t_{j_0}$ is the final iteration of  our algorithm. We have:
\begin{itemize}
\item for $k\in\rel^c$, $h_k^{(t_{j_0}-1)}=h_k$. Indeed, $t_k<0$, hence the components $k$ have been deactivated before Iteration $t_{j_0}-1$, and have the same value for the last two iterations. 
\item  for $k\in\rel$, $h_k\geq \beta h_k^{(t_{j_0}-1)}$. Indeed, at worst, the component $k$ was active during Iteration $t_{j_0}-1$ and have been multiplied by $\beta$.
\end{itemize}
Therefore:
$$\prod\limits_{k=1}^d h_k  \geq \beta^r \prod\limits_{k=1}^dh_k^{(t_{j_0}-1)}$$
%Note that $h_{j_0}=\beta h_{j_0}^{(t_{j_0}-1)}$, hence:
and
$$ h_{j_0} \seuil{hj_0}=\Cl \sqrt{\frac{(\log n)^a}{n\prod_{k=1}^d h_k}} 
\leq \beta^{-\frac{r}2} h_{j_0}^{(t_{j_0}-1)} \seuil{h^{(t_{j_0}-1)} j_0}.$$
To summarize, we have
\begin{align*}
\1_{\lbrace\hat{h}=h\rbrace\cap \Ehp} \left\vert \fbar{h}- \E\left[\fbar{h}\right]\right\vert 
&\leq \1_{\lbrace\hat{h}=h\rbrace\cap \Ehp} \sigma_h
=\1_{\lbrace\hat{h}=h\rbrace\cap \Ehp} \frac{\Cs}{\Cl} h_{j_0} \seuil{h j_0}\\
&\leq \1_{\lbrace\hat{h}=h\rbrace\cap \Ehp} \beta^{-\frac{r}2}\frac{\Cs}{\Cl} h_{j_0}^{(t_{j_0}-1)} \seuil{h^{(t_{j_0}-1)} j_0}\\
&\leq \1_{\lbrace\hat{h}=h\rbrace\cap \Ehp} 4\beta^{-\frac{r}2}\frac{\Cs}{\Cl} h_{j_0}^{(t_{j_0}-1)} \left\vert \E\left[ \bar{Z}_{h^{(t_{j_0}-1)} j_0}\right] \right\vert.
\end{align*}
Then we apply 2. of \Cref{lmZbar}:
$$\left\vert \E\left[ \bar{Z}_{h^{(t_{j_0}-1)} j_0}\right]\right\vert
\leq \CEZ \left(h^{(t_{j_0}-1)}_{j_0}\right)^{s-1}.$$
Therefore:
\begin{align*}
\1_{\lbrace\hat{h}=h\rbrace\cap \Ehp} \left\vert \fbar{h}- \E\left[\fbar{h}\right]\right\vert 
&\leq \1_{\lbrace\hat{h}=h\rbrace\cap \Ehp}  4\beta^{-\frac{r}2}\frac{\Cs}{\Cl} h_{j_0}^{(t_{j_0}-1)}\times \CEZ  \left(h^{(t_{j_0}-1)}_{j_0}\right)^{s-1}\\
&\leq \tfrac{4\CEZ\Cs\beta^{-\frac{r}2}}{\Cl}  \left(\beta^{t_{j_0}-1} h_{0}\right)^{s}
= \tfrac{4\CEZ\Cs\beta^{-\frac{r}2-s}}{\Cl} \left(\beta^{t_{j_0}} h_{0}\right)^{s}\\
&\leq \tfrac{4\CEZ\Cs\beta^{-\frac{r}2-s}}{\Cl} \left(\beta^{\tminimax} h_{0}\right)^{s}
=\tfrac{4\CEZ\Cs\beta^{-\frac{r}2-s}}{\Cl} \left(\CAm (\log n)^{\Am} n^{-\frac1{2s+r}}\right)^{s} \\
&= \tfrac{4{\CAm}^s\CEZ\Cs\beta^{-\frac{r}2-s}}{\Cl}  (\log n)^{s\Am} n^{-\frac{s}{2s+r}}.
\end{align*}
}
Reuniting the two cases, we obtain Inequality \eqref{eMajDeviation}:
\begin{align*}
\1_{\lbrace\hat{h}=h\rbrace\cap \Ehp} &\left\vert \fbar{h}- \E\left[\fbar{h}\right]\right\vert 
\leq \1_{\lbrace\hat{h}=h\rbrace\cap \Ehp}\sigma_h\\
&\leq \max\left( \tfrac{\Cs}{{\beta}^{(d-r)/2}{\CAm}^{r/2}} (\log n)^{(a-\Am r)/2},
 \tfrac{4{\CAm}^s\CEZ\Cs\beta^{-\frac{r}2-s}}{\Cl}  (\log n)^{s\Am} \right)
n^{-\frac{s}{2s+r}}.
\end{align*}
\end{enumerate}

\subsection{Proof of  Proposition \ref{complexity}}
\label{proofcomplexity}
Let us evaluate the number of operations of our procedure.
During the Reverse Step, each bandwidth of $\Act^{(-1)}$ can be multiplied by $\beta^{-1}$ several times until the loop condition is achieved:
$$(\Act^{(t)}\neq\emptyset ) \& (\max \hat{h}_k^{(t)}\leq \beta).$$
In particular, $\max \hat{h}_k^{(t)}\leq 1.$
Since $\hat{h}_k^{(t)}=h_0\beta^{-|t_k|}$, 
$$|t_k|
=\log\left( \frac{\hat{h}_k^{(t)}}{h_0}\right)/ \log \left(\beta^{-1}\right)
\leq \frac{\log(h_0^{-1})}{\log(\beta^{-1})}=O\left(\frac{\log(n)}{d(2p+1)}\right)$$
using the lower bound on $h_0$ \eqref{ineq-initialisation}.
Thus, during this Reverse Step,  note that only 
$|\Act^{(-1)}|$ components are updated and:
\begin{itemize}
\item the number of updates of the $Z_{hj}$'s is of order  $\frac{|\Act^{(-1)}|}{d(2p+1)}\log(n)$ given the above remark,
\item the computation of the $Z_{hj}$'s and the comparison to the threshold cost $\mathcal{O}(|\Act^{(-1)}|n)$ operations.
\end{itemize}
Therefore at worst, there are $\mathcal{O}\left(\frac{|\Act^{(-1)}|^2}{d}\log(n)n\right)$ operation during the Reverse Step.

For the Direct Step, the stopping condition is $\left(\prod\limits_{k=1}^d \hat{h}_k^{(t)} >\frac{(\log n)^{1+a}}{n}\right)$, which is satisfied for the penultimate iteration, hence:
$$\prod\limits_{k=1}^d \hat{h}_k
> \beta^d \frac{(\log n)^{1+a}}{n}.$$
We denote ${t}_k$ the deactivation times of $\hat{h}$, then
$$h_0^d\beta^{\sum_{k=1}^d t_k}
> \beta^d \frac{(\log n)^{1+a}}{n},$$
which gives
$$\sum_{k=1}^d t_k
< \frac{\log(\beta^{-d}(\log n)^{-(1+a)}nh_0^{d})}{\log(1/\beta)}.$$
Thus, during the Direct Step, note that only $|\Act^{(0)}|$ components are updated and
\begin{itemize}
\item the total number of updates of the $Z_{hj}$'s is of order  $\logb(n)$ given the above remark, 
\item the computation of the $Z_{hj}$'s and the comparison to the threshold cost $\mathcal{O}(|\Act^{(0)}|n)$ operations.
\end{itemize}
Therefore at worst, there are $\mathcal{O}({|\Act^{(-1)}|}\log(n)n)$ operations during the Direct Step.
%
%
%
%Number of updates of the bandwidth (for all components): 
%\begin{itemize}
%\item Reverse Step : $\frac{|\Act^{(-1)|}{d(2p+1)}\log(n)$ because of the stopping condition and the value of $h_0$
%\item Direct Step: $\frac{1}{\log (1/\beta)}\log(n)$ because of the stopping condition
%\end{itemize}
%
%\medskip
%Each time : 
%\begin{itemize}
%\item
%computation of $Z_{hj}$ : $O(n?)$ operations
%\item 
%computation of the threshold and comparison  : $O(d)$ operations
%\end{itemize}
%
Using $|\Act^{(-1)}|+|\Act^{(0)}|\leq d$, the sum of these two steps leads to the proposition.

%%%%%%%%%%%%%%%%%%%%%%%%%%%

%% file: appendice.tex
% !TEX root =Rodeo-Adap-submitted.tex

%%%%%%%%%%%%%%%%%%%%%%%%
%%%%%%%%%%%%%%%%%%%%%%%%
\section{Appendix }
\subsection{Lemmas}
The following lemmas are mainly proved in \cite{Jeanne1}. Note that
some adjustments have been made from their initial versions. In particular, we have refined points 2. of \Cref{lmfbar} and of \Cref{lmZbar} to take into account the extension of our results to H\"older smoothness. In the sequel, we only prove results of subsequent lemmas which were not established in \cite{Jeanne1}.
\begin{lemma}[Lemma~5 of \cite{Jeanne1}: $\fbar{h}$ behaviour]\label{lmfbar}
 %For any bandwidth $h\in(0,1]^d$, and any $i=1:n$, let us denote $\fbar{hi}:= \frac{\K_h(w-W_i)}{\fX(X_i)}$. Then, 
 Under Assumption \AfXmin, for any bandwidth $h\in(0,1]^d$, and any $i\in 1:n$, 
\begin{enumerate}
\item Let $\CEbar:=\Vert f\Vert_{\infty,\UUn}\Vert K\Vert_1^d$. Then
$$\left\vert \E{\fbar{h1}} \right\vert \leq \E{\left\vert\fbar{h1} \right\vert}
\leq \CEbar.$$
\item If $f$ has only $r$ relevant components $\rel$ and belongs to $\mathcal{H}_d(s, L)$ and
%Under Assumptions \Afsparse and \Afreg,  
if the order $p$ of the kernel $K$ is larger than or equal to $s$,
\begin{equation}\label{lmfbarMajBbar}
\left\vert\biasbar{h} \right\vert
\leq \CBbar \sum\limits_{k\in\rel} h_k^s,
\end{equation}
with $\CBbar>0$ a constant only depending on  $L$, $s$ and $K$.
\item Let $\Bernf{h}:=\lbrace \vert \fbar{h}-\mathbb{E}[\fbar{h}]\vert \leq \sigma_h\rbrace$,
 where $\sigma_h:=\Cs\sqrt{ \frac{(\log n)^a}{ n \prod\limits_{k=1}^d h_k}}$ with
$ \Cs=\frac{2\Vert K\Vert^d_2 \Vert f\Vert_{\infty,\UUn}^{\frac12}}{\delta^{\frac12}}$.
If  Cond($h$): $\prod\limits_{k=1}^d h_k\geq \frac{4\Vert K \Vert^{2d}_{\infty}}{9\delta^2\Cs^2} \frac{(\log n)^a}{n}$ is satisfied, then:
 $$\mathbb{P}\left( \Bernf{h}^c \right) \leq 2 e^{-(\log n)^a}.$$
 \item Let  $\Bernfabs{h}:=\lbrace \left\vert \frac1n\sum\limits_{i=1}^n \vert\fbar{hi}\vert-\mathbb{E}[\vert\fbar{h}\vert]\right\vert \leq \CEbar \rbrace$.
 Then 
 $$\mathbb{P}\left( \Bernfabs{h}^c \right) \leq 2 e^{-\CgfA n\prod_{k=1}^d h_k},$$
 \end{enumerate}
 with $\CgfA:= \min\left( \frac{\CEbar{}^2}{\Cs{}^2} ; \frac{3\delta \CEbar}{4\Vert K \Vert^d_{\infty}} \right)$.
\end{lemma}

%%%%%%%%
\begin{lemma}[Lemma~6 of \cite{Jeanne1}: $\bar{Z}_{hj}$ behaviour]\label{lmZbar}
  If $K$ is chosen  as in \Cref{sec:sparsesmooth}, and under Assumption \AfXmin, for any $j\in 1,\dots,d$ and any bandwidth $h\in(0,h_0]^{d}$, we have the following results.
\begin{enumerate}
\item 
Let $\CEZA
 :=\Vert f\Vert_{\infty,\UUn} \Vert J\Vert_1 \Vert K\Vert_1^{d-1}$. We have
 \begin{equation*}%\label{e2.lmZbar2}
   \E{\vert \bar{Z}_{h1j}\vert}
 \leq \CEZA h_j^{-1}.
\end{equation*}  
\item 
%Under Assumption \Afsparse, 
 If $f$ has only $r$ relevant components $\rel$,
for $j\notin\mathcal{R}$:
\begin{equation*}%\label{eEnp}
\E{ \bar{Z}_{hj}}=0,
\end{equation*}
and if in addition $f$ belongs to $\mathcal{H}_d(s, L)$, 
%under Assumption \Afreg, 
for $j \in \rel$: 
 \begin{equation}\label{lmZbarMajEZ}
\vert\mathbb{E}[\bar Z_{h,j}]\vert
\leq \CEZ h_j^{s-1},
\end{equation}
where $\CEZ:=\left(\int |z^s K(z)|dz\right)\frac{\Vert K\Vert_1^{r-1}L }{(s-1)!}$ denoting $(s-1)!:=(s-q+1)(s-q+2)\dots(s-1).$
\item Let $\BernZ{h,j}:=\lbrace \vert \bar{Z}_{hj}-\E{\bar{Z}_{hj}}\vert\leq\frac12\seuil{hj} \rbrace$.
%Under \Cref{AfXmin,Afreg,Afsparse},
%for some $\rho>0$,
If the bandwidth satisfies:
\begin{enumerate}
\item[$\text{Cond}_{\bar{Z}}(h)$:] $ \prod\limits_{k=1}^d h_k
\geq   \condZbar\frac{(\log n)^a}{n}$, with $\condZbar:=\frac{4 \Vert J\Vert_\infty^2\Vert K\Vert_\infty^{2(d-1)}}{ 3^2 \Vert f \Vert_{\infty,\UUn} \Vert J\Vert_2^2 \Vert K\Vert_2^{2(d-1)}} $,
\end{enumerate} 
then: 
$$\PP{\BernZ{h,j}^c}\leq 2 e^{-\frac{\delta }{\Vert f\Vert_{\infty,\UUn}} (\log n)^a}.$$

\item Let $\BernZabs{h,j}:=\lbrace\vert \frac1n \sum\limits_{i=1}^n \vert\bar{Z}_{hij}\vert -\E{\vert\bar{Z}_{h1j}\vert}\vert\leq \CEZA h_j^{-1} \rbrace$. Then, %under \Cref{AfXmin,Afreg,Afsparse}:
$$
\PP{\BernZabs{h,j}^c}\leq 2 e^{-\CgZA n\prod_{k=1}^d h_k},
$$
with $\CgZA:= \min\left(\frac{\delta\CEZA^2}{4\Vert f\Vert_{\infty,\UUn} \Vert J\Vert_2^2 \Vert K\Vert_2^{2(d-1)}} ; \frac{3\delta\CEZA}{4\Vert K\Vert_\infty^{d-1}\Vert J\Vert_\infty}\right)$.
\end{enumerate}
\end{lemma}
%\red{Incorporer le lemme \ref{lmEZ} en 2. du \Cref{lmZbar} }\\
%\begin{lemma}\label{lmEZ}
%If $K$ is chosen as in Section \ref{KC}, under assumptions \ref{Afreg,Afsparse}
%$$\forall j \in \rel\qquad\vert\mathbb{E}[\bar Z_{h,j}]\vert
%\leq \CEZ h_j^{s-1},$$
%where $\CEZ:=\left(\int |z^s K(z)|dz\right)\frac{\Vert K\Vert_1^{r-1}L }{(s-1)!}$ denoting $(s-1)!:=(s-q+1)(s-q+2)\dots(s-1).$
%\end{lemma}
%
% We denote $\DZ{hk}:=Z_{hk}-\bar{Z}_{hk}$ and $\Delta_h:= \f{h}-\fbar{h}$ for any bandwidth $h\in(\R_+^*)^d$ and any component $k\in\lbrace1:d\rbrace$.
\begin{lemma}\label{lmDelta}
For any $h\in\HhpRev\cup\HhpDir$ and any component $j\in 1:d$,
under Assumptions \AfXmin and \AestfX,  if $\sqrt{\prod\limits_{k=1}^dh_{k}}\leq 1$, then
\begin{enumerate}
\item we have:
\begin{equation*}
\1_{\BZA{hj}\cap\Atilde}\left\vert\DZ{hj}\right\vert 
\leq  \frac{1}{4} \seuil{hj}
\end{equation*}
\item 
for $\CMD:=\frac{4 M_X\CEbar}{\delta\Cs} $:
 \begin{equation*}
\1_{\Atilde\cap\Bernfabs{h}} \left\vert\Delta_{h}\right\vert
\leq  \CMD\sigma_h.
\end{equation*}
\end{enumerate}
\end{lemma}

%\begin{lemma}[Bernstein's inequality]{\color{red}à enlever ? }\label{lmBernstein}
%Let $U_1,\dots,U_n$  be independent random variables almost surely uniformly bounded by a positive constant $M>0$ and such that for $i=1,\dots,n$, $\mathbb{E}[U_i^2]\leq v$ . Then for any $\lambda>0$,
%$$ \mathbb{P}\left(\left|\frac{1}{n}\sum_{i=1}^n U_i-\mathbb{E}[U_i]\right|\geq \lambda\right)\leq
%2\exp\left(-\min\left(\frac{n\lambda^2}{4v},\frac{3n\lambda}{4M}\right)\right).$$
%\end{lemma}
%\noindent
%Note that this version is a simple consequence of Birg\'e and Massart (p.366 of \citep{BM98}).

\begin{lemma}[Taylor's theorem] \label{lmTaylor}
Let $g:[0,1]\rightarrow\mathbb{R}$ be a function of class $C^q$. Then we have:
\begin{equation*}
g(1)-g(0)=\sum\limits_{l=1}^q \frac{g^{(l)}(0)}{l!} +\int_{t_1=0}^1 \int_{t_2=0}^{t_1}\dots \int_{t_q=0}^{t_{q-1}} (g^{(q)}(t_q)-g^{(q)}(0)) dt_q dt_{q-1}\dots dt_1.
\end{equation*}
\end{lemma}

%%%%%
\subsection{Proof of Inequality \eqref{lmfbarMajBbar} in \Cref{lmfbar}}
We recall that the notation $\cdot$ means  the multiplication term by term of two vectors, then we have:
\begin{align*}
\bar{B}_h
=\E{ \fbar{h}} -f(w)
&=\int_{u\in\mathbb{R}^d} \left( \prod\limits_{k=1}^d \frac{K(h_k^{-1}(w_k-u_k))}{h_k} \right)f(u) du - f(w)\\
&=\int_{z\in\mathbb{R}^d} \left(\prod\limits_{k=1}^d K(z_k)\right) (f(w-h\cdot z)-f(w)) dz.%\label{eBias}
\end{align*}
For any $z\in\mathbb{R}^d$, let us introduce the notations $\overline{z}_0:=w$ and  for $k\in 1,\dots,d$, $\overline{z}_k:=w-\sum_{j=1}^k h_jz_je_j$, where $\lbrace e_j\rbrace_{j=1}^d$ is the canonical basis of $\R^d$. Then, we write:
\begin{align*}
f(w-h.z)-f(w)
=\sum\limits_{k=1}^d f(\overline{z}_k)-f(\overline{z}_{k-1})
=\sum\limits_{k\in\mathcal{R}} f(\overline{z}_k)-f(\overline{z}_{k-1}),%\label{eDecompBias}
\end{align*}
since for $k\notin\rel$, $f(\overline{z}_k)-f(\overline{z}_{k-1})=0$.
We apply Taylor's theorem (cf \Cref{lmTaylor}) to the functions $g_k:t\in[0,1]\mapsto f(\overline{z}_{k-1}-th_kz_ke_k)$, $k\in\rel$:
\begin{align*}
f(\overline{z}_k)-f(\overline{z}_{k-1})
=g_k(1)-g_k(0)
=\sum\limits_{l=1}^{q}\frac{(\text{-}z_kh_k)^l}{l!}\partial_k^l f(\overline{z}_{k-1}) + J_k,
\end{align*}
where we recall that $q$ is the largest integer smaller than $s$ and with
\begin{align*}
J_k
&:=\int\limits_{0\leq t_q\leq \dots\leq t_1\leq 1} \left(g_k^{(q)}(t_q)-g_k^{(q)}(0)\right)dt_{1:q}\\
&=(\text{-}h_kz_k)^q\int\limits_{0\leq t_q\leq \dots\leq t_1\leq 1} \left(\partial_k^qf(\overline{z}_{k-1}-t_qh_kz_ke_k)-\partial_k^qf(\overline{z}_{k-1})\right)dt_{1:q}.
\end{align*}
We denote $\text{I}_k:=\int_{z\in\mathbb{R}^d} \left(\prod\limits_{k'=1}^d K(z_{k'})\right) J_k dz$ and for any $z\in\mathbb{R}^d$, we denote $z_{-k}\in\R^{d-1}$ the vector $z$ without its $k^{th}$ variable, then we obtain:
\begin{align*}
\bar{B}_h
&=\sum\limits_{k\in\mathcal{R}} 
\int_{z\in\R^{d}} \left(\prod\limits_{k'=1}^d K(z_{k'})\right) \left(
J_k
+ \sum\limits_{l=1}^{q}\frac{(\text{-}h_k)^l}{l!}\partial_k^l f(\overline{z}_{k-1}) 
 z_k^l 
  \right)dz\\
&=\sum\limits_{k\in\mathcal{R}}\left( \text{I}_k 
+\sum\limits_{l=1}^{q}
\text{II}_{k,l}
\right),
\end{align*}
where 
\begin{align*}
\text{II}_{k,l}
:&= \int_{z_{-k}\in\R^{d-1}} \left(\prod\limits_{k'\neq k} K(z_{k'})\right)
\frac{(\text{-}h_k)^l}{l!}\partial_k^l f(\overline{z}_{k-1}) 
\int_{z_k\in\mathbb{R}} z_k^l K(z_k) dz_kdz_{-k} \\
&=\frac{(\text{-}h_k)^l }{l!}
\int_{z_{-k}\in\R^{d-1}}\partial_k^l f(\overline{z}_{k-1}) \left(\prod\limits_{k'\neq k} K(z_{k'})\right)
dz_{-k} 
\times \int_{t\in\mathbb{R}} t^l K(t) dt
=0,
\end{align*}
since $K$ is of order $p\geq s> q$. 
So,
\begin{align*}
\bar{B}_h
&=\sum\limits_{k\in\mathcal{R}}\text{I}_k.
\label{eBias2}
\end{align*}
%Using \Cref{Afreg}, for any $k\in\rel$:
%First, 
%\begin{align}
%\left\vert\partial_k^p f(\overline{z}_{k-1})-\partial_k^p f(w)\right\vert
%\leq \sum\limits_{k'\in\rel} \left\vert\partial_k^p f(\overline{z}_{k'})-\partial_k^p f(\overline{z}_{k'-1})\right\vert
%\leq L\sum\limits_{k'\in\rel} (h_{k'}\vert z_{k'}\vert)^{s-p}
%\end{align}
%Then we can write:
%$$\text{II}_k=\frac{(\text{-}h_k)^p \int_{t\in\mathbb{R}} t^p K(t) dt}{p!}
%\partial_k^p f(w) + R_k,
% $$ with $R_k:=\frac{(\text{-}h_k)^p \int_{t\in\mathbb{R}} t^p K(t) dt}{p!}
%\int_{z_{-k}\in\R^{d-1}}(\partial_k^p f(\overline{z}_{k-1})-\partial_k^p f(w) ) \left(\prod\limits_{k''\neq k} K(z_{k''})\right) dz_{-k}$ such as:
%\begin{align}
%\vert R_k\vert 
%\leq \frac{ h_k^p \vert\int_{t\in\mathbb{R}} t^p K(t) dt\vert}{p!}
%L\sum\limits_{k'\in\rel} h_{k'}^{s-p}\int_{z_{-k}\in\R^{d-1}} \vert z_{k'}\vert^{s-p}\left\vert\prod\limits_{k''\neq k} K(z_{k''})\right\vert dz_{-k}
%=\mathcal{O}(h_k^p\sum\limits_{k'\in\mathcal{R}} h_{k'}^{s-p})
%\end{align} 
Now  we control $\vert J_k\vert $:
\begin{align*}
\vert J_k\vert 
&\leq \vert h_kz_k\vert^{q}
\left\vert\int_{0\leq t_q\leq \dots\leq t_1\leq 1} \left[\partial_k^q f(\overline{z}_{k-1}-t_qh_kz_ke_k) -  \partial_k^q f(\overline{z}_{k-1})\right] dt_{1:q}\right\vert \\
&\leq  \vert h_kz_k\vert^{q}
\int_{0\leq t_q\leq \dots\leq t_1\leq 1} L \vert t_qh_kz_k\vert^{s-q} dt_{1:q}
=\frac{L(h_k\vert z_k\vert)^{s}}{s(s-1)\dots(s-q)}.
\end{align*}
So: 
\begin{align*}
\left\vert\text{I}_k\right\vert
=\left\vert\int_{z\in\mathbb{R}^d} \left(\prod\limits_{k'=1}^d K(z_{k'})\right) J_k dz\right\vert
\leq \frac{L \Vert K\Vert_1^{d-1} \Vert(\cdot)^sK(\cdot)\Vert_1}{s(s-1)\dots(s-q)} {h_k}^{s}.
\end{align*}
%\end{itemize}
%Therefore, 
%\begin{align*}
% \overline{B}_h
% &=\sum\limits_{k\in\mathcal{R}}(\text{II}_k+\text{I}_k)\\
%&= \sum\limits_{k\in\mathcal{R}} \frac{(\text{-}h_k)^p \int_{t\in\mathbb{R}} t^p K(t) dt\  \partial_k^p f(w)}{p!}
% +\mathcal{O}(\sum\limits_{k\in\mathcal{R}} h_{k}^{s}),
% \end{align*}
% (notice that 
% $\mathcal{O}(\sum\limits_{k\in\mathcal{R}}h_k^p\sum\limits_{k'\in\mathcal{R}} h_{k'}^{s-p})
%=\mathcal{O}(\sum\limits_{k\in\mathcal{R}}h_k^s)$).\\
Finally, 
\begin{equation}
\left\vert\biasbar{h} \right\vert
\leq \CBbar \sum\limits_{k\in\rel} h_k^s,
\end{equation}
with $\CBbar:=  \frac{L \Vert K\Vert_1^{d-1} \Vert(\cdot)^sK(\cdot)\Vert_1}{s(s-1)\dots(s-q)}$.
%%%%%%%%%%%%%%%%%%%%%%%%%%%%%%%
\subsection{Proof of Inequality \eqref{lmZbarMajEZ} in \Cref{lmZbar}}
Let $j\in \rel$. Denoting $J:\R\rightarrow\R$ the function $t\mapsto tK'(t)+K(t)$, we can write
\begin{align*}
\bar{Z}_{h,j}
&= \frac1n\sum_{i=1}^n\frac{-J(\frac{w_j-W_{ij}}{h_j})\prod\limits_{k\neq j}{K}(\tfrac{w_k-W_{ik}}{h_k}) }{\fX(X_i)h_j\prod_{k=1}^d h_k}.
\end{align*}
Then, taking the expectation,
\begin{align*}
 \E[\bar{Z}_{hj}]
&= -\frac{1}{h_j} \int_{\R^d}J(z_j) \left(\prod\limits_{k\neq j}{K}(z_k)\right) f(w-h\cdot z)dz.
\end{align*}
To simplify the notations, we assume $\rel=1:r$. Then, by integration by part
\begin{align}
\mathbb{E}[\bar Z_{h,j}]
&=\int_{\R^d} \left(z_jK(z_j)\right)\left(\prod\limits_{k\neq j} K(z_k)\right)  \partial_j f(w-h\cdot z) dz \nonumber
\\
&=\int_{\R^r}  \left(\prod\limits_{k\in\mathcal{R}} K(z_k)\right) z_j \partial_jf_{\mathcal{R}}(w_{1:r}-(h.z)_{1:r})dz_{1:r}\label{eEZ},
\end{align}
where $f_{\mathcal{R}}$ is the restriction of $f$ to the first $r$ components (remember that for any $u\in\R^r$ and any $v\in\R^{d-r}$ $f_{\mathcal{R}}(u):=f_{\mathcal{R}}(u,v)$ does not depend on $v$).
%
%
%On peut notamment remarquer que pour tout $j$, la fonction $h\mapsto\mathbb{E}[Z_{h,j}]$ ne dépend que de ses composantes pertinentes. %(cf Remarque \ref{rmqDepZ}).
%On note, pour tous $j\in\mathcal{R}$, $z\in\R^d$ et $h\in(\R_+^*)^d$ :
%\begin{itemize}
%\item[-] $e_j$ le $j^e$ vecteur de la base canonique de $\R^r$,
%\item[-] $(w-Hz)_{1:r-\lbrace j\rbrace}:=w_{1:r}-(h_1z_1,h_2z_2,\dots,h_{j-1}z_{j-1},0,h_{j+1}z_{j+1},\dots,h_rz_r)$ et 
%\item[-] $G_{j,z,h}:[0,1]\rightarrow \R$  $t\mapsto \partial_jf_{\mathcal{R}}((w-Hz)_{1:r-\lbrace j\rbrace}-th_jz_je_j)$.
%\end{itemize}
%
%
Let us denote by $G_{j,z,h}:[0,1]\rightarrow \R$ the function 
$$t\mapsto \partial_jf_{\mathcal{R}}(w_1-h_1z_1,\dots,w_j-th_jz_j,\dots, w_r-h_rz_r).$$
Then 
\begin{align*}
\mathbb{E}[\bar Z_{h,j}]&=\int_{\R^r}  \left(\prod\limits_{k\in\mathcal{R}} K(z_k)\right) z_j G_{j,z,h}(1) dz_{1:r}
\\
&=\int_{\R^r}  \left(\prod\limits_{k\in\mathcal{R}} K(z_k)\right) z_j \lbrace G_{j,z,h}(1)-G_{j,z,h}(0)\rbrace dz_{1:r},
\end{align*}
since the order $p$ of $K$ satisfies: $p\geq s>q\geq 1$.
Next we use the  Taylor expansion given by Lemma~\ref{lmTaylor}: 
\begin{equation}\label{Tay}
G_{j,z,h}(1)-G_{j,z,h}(0)=\sum\limits_{l=1}^{q-1} \frac{G_{j,z,h}^{(l)}(0)}{l!} + R'_{j,z,h,q-1},
\end{equation}
where $R'_{j,z,h,q-1}:=\int_{t_1=0}^1 \int_{t_2=0}^{t_1}\dots \int_{t_{q-1}=0}^{t_{q-2}} (G_{j,z,h}^{(q-1)}(t_{q-1})-G_{j,z,h}^{(q-1)}(0)) dt_{q-1} dt_{q-2}\dots dt_1$.
But
 $$ G_{j,z,h}^{(l)}(t)=(\text{-}h_jz_j)^l \partial_j^{l+1} f_{\mathcal{R}}(w_1-h_1z_1,\dots,w_j-th_jz_j,\dots, w_r-h_rz_r).$$
 Then, the first $q-1$ terms in the r.h.s. of \eqref{Tay} vanish since 
 $\int z_j^{l+1}K(z_j)dz_j=0.$
%On reprend \eqref{eEZ}. En intégrant par rapport à $z_j$ et en remplaçant $G_{j,z,h}(1)-G_{j,z,h}(0)$ par son développement de Taylor, on obtient un terme contenant le reste intégral et $q-1$ termes (pour $l=1:(q-1)$) de la forme :
%$$ \int_{\R}  z_j K_j(z_j) \frac{(\text{-}h_jz_j)^l}{l!} \partial_j^{l+1} f_{\mathcal{R}}((w-Hz)_{1:r-\lbrace j\rbrace})dz_j=0,$$
%car $\partial_j^{l+1} f_{\mathcal{R}}((w-Hz)_{1:r-\lbrace j\rbrace})$ ne dépend pas de $z_j$, donc le terme s'annule par l'hypothèse K\ref{hypK} sur les moments nuls de $K$.
Now, we will bound the integral remainder of \eqref{Tay}.
Using that $f$ belongs to $\mathcal{H}_d(s, L)$, for all $t\in[0,1]$,
%$$\left|\partial_j^{q}f(w_1-h_1z_1,\dots,w_j-th_jz_je_j,\dots, w_r-h_rz_r) - \partial_j^{q}f(w_1-h_1z_1,\dots,w_j,\dots, w_r-h_rz_r)\right|\leq L |t|^{s- q}$$
$$\left| G_{j,z,h}^{(q-1)}(t)- G_{j,z,h}^{(q-1)}(0)\right|
\leq |h_jz_j|^{q-1} L |th_jz_j|^{s- q},$$
since $w-h\cdot z+(1-t)h_jz_je_j\in \UUn$.
Hence
\begin{align*}
\vert R'_{j,z,h,q-1}\vert 
&\leq \int_{t_1=0}^1 \int_{t_2=0}^{t_1}\dots \int_{t_{q-1}=0}^{t_{q-2}} \left\vert G_{j,z,h}^{(q-1)}(t_{q-1})-G_{j,z,h}^{(q-1)}(0)\right\vert dt_{q-1} dt_{q-2}\dots dt_1
\\
&\leq L (h_j\vert z_j\vert)^{s-1} \int_{t_1=0}^1 
\int_{t_2=0}^{t_1}\dots \int_{t_{q-1}=0}^{t_{q-2}} t_{q-1}^{s-q} dt_{q-1} dt_{q-2}\dots dt_1
%\\&
=\frac{L (h_j\vert z_j\vert)^{s-1}}{(s-1)!}, 
\end{align*}
denoting $(s-1)!:= (s-q+1)(s-q+2)\dots(s-1)$.
Finally,
\begin{align*}
\vert\mathbb{E}[\bar Z_{h,j}]\vert 
&=\left|\int_{\R^r}  \left(\prod\limits_{k\in\mathcal{R}} K(z_k)\right) z_j R'_{j,z,h,q-1} dz_{1:r}\right|
%\\&
\leq \int_{\R^r}  \left(\prod\limits_{k\in\mathcal{R}} \vert K(z_k)\vert \right) \vert z_j \vert \frac{L (h_j\vert z_j\vert)^{s-1}}{(s-1)!} dz_{1:r}
\\
& \leq \frac{L h_j^{s-1}}{(s-1)!}
\left(\prod\limits_{k\in\mathcal{R}\setminus\lbrace j\rbrace}\Vert K\Vert_1\right) 
\int_{\R} \vert z_j\vert^{s} \vert K(z_j)\vert  dz_{1:r} 
%\\&
\leq \CEZ h_j^{s-1},
\end{align*}
denoting $\CEZ:=\left(\int_{\R} \vert z\vert^{s} \vert K(z)\vert dz\right) \Vert K\Vert_1^{r-1}L/(s-1)!$.

\subsection{Proof of \Cref{lmDelta}}
Before establishing the upper bounds, let us control $\1_{\Atilde}\left\Vert\frac{\fX-\hfX}{\hfX}\right\Vert_{\infty,\Un}$.
First, using Assumption \AfXmin:
$$\delta:=\inf_{u \in \mathcal{U}_1}\fX(u) >0,$$
remark that: for any $u\in \Un$,
\begin{align*}
\1_{\Atilde} \hfX(u)
&\geq  \1_{\Atilde}\left( \fX(u) - \Vert \fX-\hfX \Vert_{\infty,\mathcal{U}_1}\right)\\
&\geq \1_{\Atilde}\left( \delta -  M_{X} \frac{(\log n)^{\frac{a}2}}{\sqrt{n}} \right)\quad \text{ by \ref{fXtildeAccuracy}},\\
&\geq \1_{\Atilde} \frac{\delta}2 \quad \text{(for }n\text{ large enough).}
\end{align*}
Therefore: 
$$\dX:= \inf_{u \in \Un}\hfX(u)\geq  \1_{\Atilde} \frac{\delta}2,$$
which leads to:
\begin{align}
\1_{\Atilde}\left\Vert\frac{\fX-\hfX}{\hfX}\right\Vert_{\infty,\Un}
&\leq  \1_{\Atilde}\frac{ \left\Vert \fX-\hfX\right\Vert_{\infty,\Un} }{\dX}
\nonumber\\
&\leq \frac{2 M_X}{\delta} \frac{(\log n)^{a/2}}{n^{1/2}}.\label{eMajRatiofXtilde}
\end{align}
Let us now prove the first upper bound.
\begin{enumerate}
\item {We still denote, for any bandwidth $h$, any component $k$ and any observation $i$,
$$ \bar{Z}_{hik}:=\frac{\partial}{\partial h_k} \left( \frac{\K_h(w-W_i)}{\fX(X_i)}\right),$$
such that $\bar{Z}_{hk}=\frac1n \sum\limits_{i=1}^n \bar{Z}_{hik}$, with $\lbrace \bar{Z}_{hik}\rbrace_{i=1}^n$ i.i.d..
Then we can write:
\begin{equation*}
\DZ{hk}:= Z_{hk}-\bar{Z}_{hk} 
= \frac1n \sum\limits_{i=1}^n \left( \tfrac{\fX}{\hfX}(X_i)-1\right) \bar{Z}_{hik} 
= \frac1n \sum\limits_{i=1}^n \left( \tfrac{\fX -\hfX}{\hfX}(X_i)\right) \bar{Z}_{hik}.
\end{equation*}
Note that since $K$ is compactly supported, if $X_i\notin\Un$,
$$\bar{Z}_{hik}= 0.$$
Hence:
\begin{align*}
\left\vert \DZ{hk}\right\vert
&\leq   \left\Vert \tfrac{\fX -\hfX}{\hfX}\right\Vert_{\infty,\Un} \times \frac1n \sum\limits_{i=1}^n \vert \bar{Z}_{hik}\vert
\\
&\leq   \left\Vert \tfrac{\fX -\hfX}{\hfX}\right\Vert_{\infty,\Un} 
\times \left( \E\left[\left\vert\bar{Z}_{h1k}\right\vert\right]
	+ \frac1n \sum\limits_{i=1}^n \left\vert \bar{Z}_{hik}\right\vert-\E\left[\left\vert\bar{Z}_{hik}\right\vert\right]
	\right).
\end{align*}
Using the above Inequality \eqref{eMajRatiofXtilde} and the upper bounds 1. and 4. of \Cref{lmZbar}:
\begin{align*}
\1_{\Atilde\cap \BernZabs{h,k}}\left\vert \DZ{hk}\right\vert
&\leq   \left( \frac{2 M_X}{\delta} \frac{(\log n)^{a/2}}{n^{1/2}}\right) \times 
2\CEZA h_k^{-1}\\
&\leq \frac14 \seuil{h,k}:=\frac{\Cl}{4} \frac{(\log n)^{a/2}}{n^{1/2}h_k\left(\prod\limits_{k'=1}^d h_{k'}\right)^{1/2}},
\end{align*}
if $\left(\prod\limits_{k'=1}^d h_{k'}\right)^{1/2}\leq \frac{\delta\Cl}{16M_X \CEZA}$.
Note that $M_X$ is determined in order to satisfy:
$$\frac{\delta\Cl}{16M_X \CEZA}= 1.$$
Hence the condition on the bandwidth becomes:
$$\left(\prod\limits_{k'=1}^d h_{k'}\right)^{1/2}\leq 1.$$
}
\item{ We still denote, for any bandwidth $h$ and any observation $i$,
$$ \fbar{hi}:= \frac{\K_h(w-W_i)}{\fX(X_i)},$$
such that $\fbar{h}=\frac1n \sum\limits_{i=1}^n \fbar{hi}$, with $\lbrace \fbar{hi}\rbrace_{i=1}^n$ i.i.d.
Then we can write:
\begin{equation*}
\Delta_{h}:= \hat{f}_{h}(w)-\fbar{h} 
= \frac1n \sum\limits_{i=1}^n \left( \tfrac{\fX}{\hfX}(X_i)-1\right) \fbar{hi} 
= \frac1n \sum\limits_{i=1}^n \left( \tfrac{\fX -\hfX}{\hfX}(X_i)\right) \fbar{hi}.
\end{equation*}
Note that since $K$ is compactly supported, if $X_i\notin\Un$,
$$\fbar{hi} = 0.$$
Hence:
\begin{align*}
\left\vert \Delta_{h}\right\vert
&\leq   \left\Vert \tfrac{\fX -\hfX}{\hfX}\right\Vert_{\infty,\Un} \times \frac1n \sum\limits_{i=1}^n \vert \fbar{hi}\vert
\\
&\leq   \left\Vert \tfrac{\fX -\hfX}{\hfX}\right\Vert_{\infty,\Un} 
\times \left( \E\left[\left\vert\fbar{h1}\right\vert\right]
	+ \frac1n \sum\limits_{i=1}^n \left\vert \fbar{hi}\right\vert-\E\left[\left\vert\fbar{hi}\right\vert\right]
	\right).
\end{align*}

Using the above Inequality \eqref{eMajRatiofXtilde} and the upper bounds 1. and 4. of \Cref{lmfbar}:
\begin{align*}
\1_{\Atilde\cap \Bernfabs{h}}\left\vert \Delta_{h}\right\vert
&\leq   \left( \frac{2 M_X}{\delta} \frac{(\log n)^{a/2}}{n^{1/2}}\right) \times 
2\CEbar\\
&=\frac{4 M_X\CEbar}{\delta\Cs} \sigma_h\left(\prod\limits_{k'=1}^d h_{k'}\right)^{1/2}\leq \CMD \sigma_h.
\end{align*}
}
\end{enumerate}
\subsection{Proof of \Cref{propfXtilde}}
The proof is very similar to the Proposition 1 of  \citep{Jeanne1}. The main modification is due to the tighter $\log$ exponent in \ref{fXtildeAccuracy} and the enlarged neighborhood $\Un$ of $x$.
We introduce the classical kernel density estimator $\fXK$: for any $u\in\R^{d_1}$ and a bandwidth $\hX\in\R_+^*$ to be specified later, 
\begin{equation}\label{eDef fXK}
  \fXK (u):= \frac{1}{\nX.\hX^{d_1}}\sum\limits_{i=1}^{\nX} \prod\limits_{j=1}^{d_1} \KX \left(\frac{u_j-\X_{ij}}{\hX}\right),
  \end{equation}
where $\KX:\R\rightarrow\R$ is a kernel which is compactly supported, of class $C^1$ and
of order $p_X\geq \frac{d_1}{2(c-1)},$ where we recall that $c>1$ is defined by $\nX=n^c$. We first show that there exists $C_X>0$ such that for any $\xi>0$:
\begin{equation}
\PP{\Vert \fX-\fXK\Vert_{\infty,\Un} > C_X\frac{(\log n)^{\frac{1+\xi}2}}{\sqrt{n}} }\leq  \mathcal{O}\left(\nX^{d_1+1} \exp\left(- (\log n)^{1+\xi}\right)\right).
%\exp(-(\log n)^{1+\frac{\xi}2}).
\label{eCond1fXK}
\end{equation} 
Then we set 
$$\hfX\equiv \fXK\vee n^{-\frac12},$$
and 
we shall prove that this estimator  satisfies \ref{fXtildemin} and \ref{fXtildeAccuracy} for $\hfX$.

Let us prove Inequality \eqref{eCond1fXK}.
Let us first explicit $\fXK$'s behaviour. Following \Cref{lmfXK} gives a pointwise concentration inequality and a control of the bias of $\fXK$ on $\Un$.
% using Bernstein's Inequality. 
%It is an extended version of  2. of Lemma 4 in \citep{Jeanne1} where we have generalized the $\log n$ exponent $\frac32$ to values of the form $1+\varepsilon$, $\varepsilon>0$ and enlarged the considered neighbourhood of $x$ to $\Un$.\\
We introduce an enlarged neighborhood of $\Un$:
$$\Up:=\left\lbrace u'= u- \hX z : u\in\Un, z\in\supp(\KX)\right\rbrace.$$

\begin{lemma}[$\fXK$ behaviour]\label{lmfXK}
The estimator $\fXK$ satisfies the following results:
  \begin{enumerate}
  \item If there exists $q_X\in\mathbb{N}$ such that $\fX$ is $C^{q_X}$ on $\Up$ and such that $\KX$ has $q_X-1$ zero moments, then there exists a positive constant 
${\normalfont \CbiasX'}$ such that
  \begin{align*}
\left\Vert \E{\fXK}- \fX\right\Vert_{\infty,\Un}
&\leq {\normalfont \CbiasX'} h_X^{q_X}.
\end{align*}
  \item For any $\xi>0$, any $u\in \Un$ and  any $\lambda>0$ such that:
  $$4 {\normalfont\CvX} \frac{(\log n)^{1+\xi}}{\nX\hX^{d_1}}
  \leq \lambda^2
  \leq \frac{9{{\normalfont\CvX}}^2}{\Vert\KX \Vert_{\infty}^{2d_1}},$$
  where ${\normalfont\CvX}:=\Vert \KX\Vert_2^{d_1}\Vert \fX\Vert_{\infty,\Up}^{\frac12}$, 
%$\lambda_X:=\sqrt{
%{\color{red}\tfrac{4\Vert \KX\Vert_2^{2d_1}\Vert \fX\Vert_{\infty,\Up}}{\hX^{d_1}\nX}
%\text{À VÉRIFIER}}  
%(\log n)^{1+\xi}}$ 
%and if the bandwidth $\hX$
%satisfies the condition 
%$$\text{Cond}_X(\hX): \quad \hX^{d_1}
%\geq 
%{\color{red}
%\frac{4\Vert \KX\Vert_{\infty}^{2d_1}}{9\Vert \KX\Vert_2^{2d_1}\Vert \fX\Vert_{\infty,\Up}}
%\text{À VÉRIFIER}}  
%\frac{(\log n)^{1+\xi}}{\nX},$$ 
\begin{align*}
\PP{\left\vert\fXK(u)-\E{\fXK(u)}\right\vert>\lambda}
\leq 2 \exp\left(- (\log n)^{1+\xi}\right).
\end{align*}
\end{enumerate} 
\end{lemma}
This lemma is proved in Section~\ref{sec:prooflemma5}. 
% Notice that if $s> p_X+1$, then $f_X$ is of class $\mathcal{C}^{p_X+1}$ on $\Up$ (defined in \eqref{edefUp}), thus $f_X$ is $(s', L')$-Hölder, with $s'=p_X+1$ and $L'=\max\limits_{j=1:d_1} \Vert\partial_j^{p_X+1}f_X\Vert_{\infty,}$. 
%We denote 
%$$(s',L'):=\left\lbrace\begin{tabular}{@{}l}
%$(s,L)$\quad if $s\leq p_X+1$\\ 
%$(p_X+1, \max\limits_{j=1:d_1} \Vert\partial_j^{p_X+1}f_X\Vert_{\infty,})$\quad  otherwise.\\ 
%\end{tabular} 
%\right.
%$$
We define $p'_X=\min(p',p_X)$, so that: $\fX$ is of class $C^{p'_X}$ and the first $p'_X-1$ moments of $\KX$ vanish. Therefore, we can apply 1. of \Cref{lmfXK}: 
\begin{align*}
\left\Vert \E{\fXK}- \fX\right\Vert_{\infty,\Un}
&\leq \CbiasX' h_X^{p'_X}.
\end{align*}
Therefore:
\begin{align*}
\left\Vert \fXK-\fX \right\Vert_{\infty,\Un}
&\leq \left\Vert \fXK-\E{\fXK} \right\Vert_{\infty,\Un} + \left\Vert \E{\fXK}-\fX \right\Vert_{\infty,\Un}
\\
&\leq \left\Vert \fXK-\E{\fXK} \right\Vert_{\infty,\Un} + \CbiasX' h_X^{p'_X},
\end{align*}
and we have for any threshold $\lambda$:
\begin{align}
\PP{\left\Vert \fXK-\fX \right\Vert_{\infty,\Un}\geq \lambda}
&\leq \PP{\left\Vert \fXK-\E{\fXK} \right\Vert_{\infty,\Un}\geq \lambda-\CbiasX' h_X^{p'_X}}. \label{eProbSansBias}
\end{align}
We have then reduced the problem to a concentration inequality of $\fXK$ in sup norm.
In order to move from a supremum on $\Un$ to a maximum on a finite set of elements of $\Un$, let us construct an $\epsilon$-net $\lbrace u_{(l)}\rbrace_{l}$ of $\Un$, in the meaning that for any $u\in\Un$, there exists $l$ such that $\Vert u-u_{(l)}\Vert_{\infty}:=\max\limits_{k=1:d_1}\vert u_k- u_{(l)k}\vert\leq \epsilon$. We denote $A>0$ such that:
$$\supp(\KX)\cup \supp(K)\subset \left[-\tfrac{A}2,\tfrac{A}2\right].$$
Set $N(\epsilon)$ is the smallest integer such that 
 $ 2\epsilon N(\epsilon)\geq A,$ % and $N:=N_1^{d_1}$
% \textit{i.e.}: 
% $$N(\epsilon):=\left\lceil\frac{A}{2\epsilon}\right\rceil,$$
and for $l\in\left(1:N(\epsilon)\right)^{d_1}$,  $u_{(l)}$ such that its $j$-th component is equal to:
 $$u_{(l)j}:=x_j-\frac{A}{2}+(2l_j-1)\epsilon.$$
 Then $\lbrace u_{(l)}\rbrace_{l\in\left(1:N(\epsilon)\right)^{d_1}}$ is an $\epsilon$-net of $\Un$. Therefore in order to obtain Inequality \eqref{eCond1fXK}, we only need to obtain the concentration inequality for each point of $\lbrace u_{(l)}: l\in(1:N(\epsilon))^{d_1}\rbrace$ and to control the difference of the function $\fXK -\E{\fXK}$ evaluated at the point $u$ and at the nearest point of $u$ in the $\epsilon$-net. More formally, we have to control the following supremum $$\sup\limits_{u\in\Un}\min\limits_{l\in(1:N(\epsilon))^{d_1}} \left\vert \fXK(u) -\E{\fXK(u) }- \fXK(u_{(l)})+\E{\fXK(u_{(l)}) } \right\vert.$$ 
For this purpose, we obtain (from Taylor's Inequality): for any $u,v\in\R^{d_1}$,
$$\left\vert\prod\limits_{k=1}^{d_1} \KX(u_k) - \prod\limits_{k=1}^{d_1} \KX(v_k) \right\vert 
\leq d_1 \Vert \KX'\Vert_{\infty} \Vert \KX\Vert_{\infty}^{d_1-1}\Vert u-v\Vert_{\infty}.$$
Therefore, for any $u,v\in\Un$:
\begin{align*}
\left\vert \fXK(u) - \fXK(v) \right\vert 
&\leq \frac{1}{\nX .\hX^{d_1}}\sum\limits_{i=1}^{\nX} \left\vert\prod\limits_{k=1}^{d_1} \KX(\tfrac{u_k-\X_{ik}}{\hX}) - \prod\limits_{k=1}^{d_1} \KX(\tfrac{v_k-\X_{ik}}{\hX}) \right\vert
\\
&\leq d_1 \Vert \KX'\Vert_{\infty} \Vert \KX\Vert_{\infty}^{d_1-1}\frac{\Vert u-v\Vert_{\infty}}{\hX^{d_1+1}}.
\end{align*}
Since $\lbrace u_{(l)}: l\in(1:N(\epsilon))^{d_1}\rbrace$ is an $\epsilon$-net of $\Un$:
\begin{align*}
\sup\limits_{u\in\Un}\min\limits_{l\in(1:N(\epsilon))^{d_1}} \left\vert \fXK(u) - \fXK(u_{(l)}) \right\vert
&\leq d_1 \Vert \KX'\Vert_{\infty} \Vert \KX\Vert_{\infty}^{d_1-1}\frac{\epsilon}{\hX^{d_1+1}},
\end{align*}
and also:
\begin{align*}
\sup\limits_{u\in\Un}\min\limits_{l\in(1:N(\epsilon))^{d_1}} \left\vert \E{\fXK(u)} - \E{\fXK(u_{(l)})} \right\vert
&\leq d_1 \Vert \KX'\Vert_{\infty} \Vert \KX\Vert_{\infty}^{d_1-1}\frac{\epsilon}{\hX^{d_1+1}}.
\end{align*}
Therefore:
\begin{align*}
\sup\limits_{u\in\Un}\min\limits_{l\in(1:N(\epsilon))^{d_1}} \left\vert \fXK(u)-\E{\fXK(u)} - \fXK(u_{(l)})+ \E{\fXK(u_{(l)})} \right\vert
&\leq 2d_1 \Vert \KX'\Vert_{\infty} \Vert \KX\Vert_{\infty}^{d_1-1}\frac{\epsilon}{\hX^{d_1+1}}.
\end{align*}
We denote $\text{C}_{\text{diff}}:=2d_1 \Vert \KX'\Vert_{\infty} \Vert \KX\Vert_{\infty}^{d_1-1}$.
We then obtain the following inequality:
\begin{align*}
\left\Vert \fXK -\E{\fXK}\right\Vert_{\infty,\Un} 
&\leq \max\limits_{l\in(1:N(\epsilon))^{d_1}} \left\vert \fXK(u_{(l)}) -\E{\fXK(u_{(l)})}\right\vert
\\
&\qquad + \sup\limits_{u\in\Un}\min\limits_{l\in(1:N(\epsilon))^{d_1}} \left\vert \fXK(u)-\E{\fXK(u)} - \fXK(u_{(l)})+ \E{\fXK(u_{(l)})} \right\vert
\\
&\leq \max\limits_{l\in(1:N(\epsilon))^{d_1}} \left\vert \fXK(u_{(l)}) -\E{\fXK(u_{(l)})}\right\vert
+\text{C}_{\text{diff}}\frac{\epsilon}{\hX^{d_1+1}}.
\end{align*}
Then the inequality \eqref{eProbSansBias} becomes: for any threshold $\lambda$,
\begin{align}
\PP{\left\Vert \fXK-\fX \right\Vert_{\infty,\Un}\geq \lambda}
&\leq \PP{\left\Vert \fXK-\E{\fXK} \right\Vert_{\infty,\Un}\geq \lambda-\CbiasX' h_X^{p'_X}} \nonumber
\\
&\leq \PP{\max\limits_{l\in(1:N(\epsilon))^{d_1}}\left\vert \fXK(u_{(l)})-\E{\fXK(u_{(l)})} \right\vert \geq \lambda-\CbiasX' h_X^{p'_X}-\text{C}_{\text{diff}}\frac{\epsilon}{\hX^{d_1+1}}}\nonumber
\\
&\leq  N(\epsilon)^{d_1} \max\limits_{l\in(1:N(\epsilon))^{d_1}}\PP{\left\vert \fXK(u_{(l)})-\E{\fXK(u_{(l)})} \right\vert \geq \lambda-\CbiasX' h_X^{p'_X}-\text{C}_{\text{diff}}\frac{\epsilon}{\hX^{d_1+1}}}.
 \label{eProbSansBiasSansSup}
\end{align}
It then remains to apply 2. of \Cref{lmfXK} for each $u_{(l)}$, $l\in(1:N(\epsilon))^{d_1}$. 
We set the following settings: 
\begin{itemize}
\item $\hX:= \nX^{-\frac{c-1}{c.d_1}}$;%(\log n)^{\frac{-d(c-1)}{d_1}};
\item $\epsilon:=\hX^{1+\frac{d_1}{2}}\nX^{-\frac12}$;
\item $\lambda:=2\lambda_X$, where $\lambda_X$ is defined by:
\begin{align*}
\lambda_X
%:=\sqrt{\tfrac{4\Vert \KX\Vert_2^{2d_1}\Vert f_X\Vert_{\infty,\Up}}{\hX^{d_1}\nX}(\log n)^{1+\xi}}
:= 2\sqrt{\CvX} (\log n)^{\frac{1+\xi}2}\hX^{-\frac{d_1}2}\nX^{-\frac12}
= 2\sqrt{\CvX} (\log n)^{\frac{1+\xi}2}\nX^{-\frac1{2c}},
\end{align*}
where we recall that $\CvX:=\Vert \KX\Vert_2^{d_1}\Vert \fX\Vert_{\infty,\Up}^{\frac12}$.
\end{itemize}
In particular, since we take $p_X\geq \frac{d_1}{2(c-1)}$ and we assume $p'\geq \frac{d_1}{2(c-1)}$, then $p'_X=\min(p',p_X)\geq \frac{d_1}{2(c-1)}$. Hence we obtain for $n$ large enough:
\begin{align*}
\CbiasX' h_X^{p'_X} 
&= \CbiasX' %(\log n)^{\frac{p'_X.d(c-1)}{d_1}} 
\nX^{-\frac{p'_X(c-1)}{c.d_1}}
\\
&\leq \CbiasX'    \nX^{-\frac{1}{2c}}
\\
&\leq \frac12 \lambda_X 
= \sqrt{ \CvX} (\log n)^{\frac{1+\xi}2}\nX^{-\frac1{2c}}
\end{align*}
and also, since $c>1$:
\begin{align*} 
\text{C}_{\text{diff}}\frac{\epsilon}{\hX^{d_1+1}}
&= \text{C}_{\text{diff}} \hX^{-\frac{d_1}2}\nX^{-\frac12}
= \text{C}_{\text{diff}}\ \nX^{-\frac{1}{2c}}
\\
&\leq \frac12 \lambda_X
= \sqrt{\CvX }(\log n)^{\frac{1+\xi}2}\nX^{-\frac1{2c}}.
\end{align*}
%for $n$ large enough, more precisely when $(\log n)^{\frac34}
%\geq \frac{2\text{C}_{\text{diff}}}{\text{C}_{\lambda X}}.$
Hence, we have
$$\lambda-\CbiasX' h_X^{p'_X}-\text{C}_{\text{diff}}\frac{\epsilon}{\hX^{d_1+1}}\geq \lambda_X,$$ and the inequality \eqref{eProbSansBiasSansSup}  becomes:
\begin{align}\label{eProbSansConcentration}
\PP{\left\Vert \fXK-\fX \right\Vert_{\infty,\Un}\geq \lambda}
&\leq  N(\epsilon)^{d_1} \max\limits_{l\in(1:N(\epsilon))^{d_1}}\PP{\left\vert \fXK(u_{(l)})-\E{\fXK(u_{(l)})} \right\vert \geq \lambda_X}.
\end{align}
We apply 2. of \Cref{lmfXK}: we verify (since $n_X=n^c$)
\begin{align*}
4 \CvX \frac{(\log n)^{1+\xi}}{\nX\hX^{d_1}}
  = \lambda^2_X
  &= 4\CvX (\log n)^{1+\xi}n^{-1}\\
  &\leq \frac{9{\CvX}^2}{\Vert\KX \Vert_{\infty}^{2d_1}}, \quad \text{(for }n \text{ large enough),}
\end{align*} 
then we obtain
%\begin{align*}
% \hX^{d_1}
% &= \nX^{-\frac{c-1}{c}}
% \\
% &\geq \frac{4\Vert \KX\Vert_{\infty}^{2d_1}}{9\Vert \KX\Vert_2^{2d_1}\Vert f_X\Vert_{\infty,\Up}} \frac{(\log n)^{\frac32}}{\nX}.
%\end{align*}
\begin{align*}
\PP{\left\vert\fXK(u_{(l)})-\E{\fXK(u_{(l)})}\right\vert>\lambda_X}
\leq 2 \exp\left(- (\log n)^{1+\xi}\right).
\end{align*}
Thus the inequality \eqref{eProbSansConcentration} becomes:
\begin{align}\label{eProbAvecNeps}
\PP{\left\Vert \fXK-\fX \right\Vert_{\infty,\Un}\geq \lambda}
&\leq  2N(\epsilon)^{d_1} \exp\left(- (\log n)^{1+\xi}\right).
\end{align}
Let us control $ 2N(\epsilon)^{d_1}$:
 \begin{align*}
 2N(\epsilon)^{d_1}
 =2\left\lceil\frac{A}{2\epsilon}\right\rceil^{d_1}
 = 2\left\lceil\frac{A}{2\hX^{1+\frac{d_1}{2}}\nX^{-\frac12}}\right\rceil^{d_1}  
% &=2\left\lceil\frac{A}{\hX^{1+\frac{d_1}{2}}\nX^{-\frac12}\log n}\right\rceil^{d_1}  
 %\\
= o\left( \nX^{d_1+1}\right)
 \end{align*}
Therefore, we have obtained the desired concentration inequality \eqref{eCond1fXK}.
 Now we consider $\hfX\equiv \fXK\vee n^{-1/2}$, therefore $\hfX$ satisfies \ref{fXtildemin}. 
Let us show it also satisfies \ref{fXtildeAccuracy}, for $n$ large enough. We first show: 
\begin{equation}\label{eImpl fXK fXtilde}
\left\lbrace \left\Vert \fXK-\fX \right\Vert_{\infty,\Un}
< \lambda\right\rbrace 
\quad  \Rightarrow\quad  
\left\lbrace\left\Vert \hfX-\fX \right\Vert_{\infty,\Un}
< \lambda\right\rbrace.
\end{equation}
Assume that for any $u\in\Un$, 
$\left\vert \fXK(u)-\fX(u) \right\vert 
< \lambda$. 
Let us fix $u\in\Un$. Three cases occurs:
\begin{enumerate}[label=(\alph*)]
\item When $\fXK(u) \geq n^{-\frac12}$, then $\hfX(u):=\fXK(u)$, and obviously: 
$$\left\vert \hfX(u)-\fX(u) \right\vert 
< \lambda.$$
\item When $\fXK(u) < n^{-\frac12}$ and $\fX(u)\geq n^{-\frac12}$, then since $\hfX(u)=n^{-\frac12}>\fXK(u)$,
$$ \left\vert \hfX(u)-\fX(u) \right\vert 
\leq \left\vert \fXK(u)-\fX(u) \right\vert 
< \lambda.$$
\item When $\fXK(u) < n^{-\frac12}$ and $\fX(u)< n^{-\frac12}$, then $\hfX(u)=n^{-\frac12}$, so for $n$ large enough:
$$ \left\vert \hfX(u)-\fX(u) \right\vert 
\leq n^{-\frac12}
< \lambda .$$
%=4\sqrt{\CvX }(\log n)^{\frac{1+\xi}2}n^{-\frac1{2}}.$$
\end{enumerate}
Therefore these three cases show Implication \eqref{eImpl fXK fXtilde}, and thus, from \Cref{eProbAvecNeps}, we obtain:
\begin{align*}
\PP{\left\Vert \hfX-\fX \right\Vert_{\infty,\Un}\geq \lambda}
\leq \PP{\left\Vert \fXK-\fX \right\Vert_{\infty,\Un}\geq \lambda}
\leq  2N(\epsilon)^{d_1} \exp\left(- (\log n)^{1+\xi}\right).
\end{align*}
Now, to obtain \ref{fXtildeAccuracy}, for $\xi$ such that $1+\frac{a-1}2<1+\xi<a$,
\begin{equation}\label{xia}
\lambda
=4\sqrt{\CvX} (\log n)^{\frac{1+\xi}2}n^{-\frac1{2}}
\leq 
M_X (\log n)^{\frac{a}2}n^{-\frac1{2}} \text{ (for }n\text{ large enough)}.
\end{equation}
Therefore:
 \begin{align*}
 \PP{\left\Vert \hfX-\fX \right\Vert_{\infty,\Un}\geq M_X (\log n)^{\frac{a}2}n^{-\frac1{2}}}
 &\leq  \PP{\left\Vert \hfX-\fX \right\Vert_{\infty,\Un}\geq \lambda}\\
&\leq  2N(\epsilon)^{d_1} \exp\left(- (\log n)^{1+\xi}\right) \\
&\leq \exp\left(- (\log n)^{1+\frac{a-1}2}\right),
 \end{align*}
that is \ref{fXtildeAccuracy}.
%%%%%%%%%%%%%%%%%%%%%
\subsection{Proof of \Cref{lmfXK}}\label{sec:prooflemma5}
The result 1. of \Cref{lmfXK} is proved in Lemma 4 of \cite{Jeanne1}.
To prove 2. of \Cref{lmfXK}, let us fix $\xi>0$. Then, we simply apply Bernstein's Inequality (see Lemma 10 in \cite{Jeanne1}). 
We define for any $u\in\Un$ and for $i \in 1:n$
$$\fXKi{i}(u):=\frac{1}{\hX^{d_1}}\prod\limits_{j=1}^{d_1} \KX \left(\frac{u_j-\X_{ij}}{\hX}\right).$$
Observe that the $\fXKi{i}(u)$'s are \textit{i.i.d.}
Then we pick up the following bounds from \cite[p. 23]{Jeanne1}:
\begin{align*}
\left\vert\fXKi{1}(u)\right\vert
&\leq  \text{M}_{\hX}:=\Vert \KX\Vert_{\infty}^{d_1} \hX^{-d_1}.\\
\var\left(\fXKi{1}(u)\right)
&\leq \text{v}_{\hX}:=\CvX\hX^{-d_1},
\end{align*} 
(we recall
$\CvX:=\Vert \KX\Vert_2^{2d_1} \Vert \fX\Vert_{\infty,\Up}$). Therefore: for any $\lambda>0$,
\begin{align*}
\PP{\left\vert\fXK(u)-\E{\fXK(u)}\right\vert>\lambda}
\leq 2 \exp\left(- \min\left(\frac{\nX \lambda^2}{4 \text{v}_{\hX}}, \frac{3\nX\lambda}{4\text{M}_{\hX}}\right)\right).
\end{align*}
Let us show that when 
$$4 \CvX \frac{(\log n)^{1+\xi}}{\nX\hX^{d_1}}
  \leq \lambda^2
  \leq \frac{9{\CvX}^2}{\Vert\KX \Vert_{\infty}^{2d_1}},$$
  then, we have
  $$ (\log n)^{1+\xi}
  \leq  \frac{\nX \lambda^2}{4 \text{v}_{\hX}}
  \leq \frac{3\nX\lambda}{4\text{M}_{\hX}} .$$
  Indeed,
  \begin{align*}
 \frac{\nX \lambda^2}{4 \text{v}_{\hX}}  
  \leq  \frac{3\nX\lambda}{4\text{M}_{\hX}} 
  &\quad \Leftrightarrow \quad 
   \lambda 
   \leq \frac{3 \text{v}_{\hX}}{\text{M}_{\hX}} 
   =\frac{3\CvX}{\Vert \KX\Vert_{\infty}^{d_1}}
   \\
   &\quad \Leftrightarrow \quad 
   \lambda^2 
   \leq \frac{9\CvX^2}{\Vert \KX\Vert_{\infty}^{2d_1}}
  \end{align*}
  and 
\begin{align*}
 (\log n)^{1+\xi}
  \leq \frac{\nX \lambda^2}{4 \text{v}_{\hX}}  
  &\quad \Leftrightarrow \quad 
  \frac{4 \CvX (\log n)^{1+\xi}}{\nX\hX^{d_1}}
  \leq \lambda^2.
  \end{align*}
  Therefore when $$4 \CvX \frac{(\log n)^{1+\xi}}{\nX\hX^{d_1}}
  \leq \lambda^2
  \leq \frac{9{\CvX}^2}{\Vert\KX \Vert_{\infty}^{2d_1}},$$
  \begin{align*}
\PP{\left\vert\fXK(u)-\E{\fXK(u)}\right\vert>\lambda}
&\leq 2 \exp\left(- \min\left(\frac{\nX \lambda^2}{4 \text{v}_{\hX}}, \frac{3\nX\lambda}{4\text{M}_{\hX}}\right)\right)
=2 \exp\left(-\frac{\nX \lambda^2}{4 \text{v}_{\hX}}\right)
\\
&\leq 2 \exp\left(- (\log n)^{1+\xi}\right).
\end{align*}

%% file: JCV.bbl
\begin{thebibliography}{}

\bibitem[Nguyen et~al., 2021]{JCV}
Nguyen, M.-L., Lacour, C., and Rivoirard, V. (2021).
\newblock Adaptive greedy algorithm for moderately large dimensions in kernel
  conditional density estimation.
\newblock {\em Submitted}.

\end{thebibliography}


\begin{thebibliography}{}

\bibitem[Bashtannyk and Hyndman, 2001]{BH01}
Bashtannyk, D.~M. and Hyndman, R.~J. (2001).
\newblock Bandwidth selection for kernel conditional density estimation.
\newblock {\em Comput. Statist. Data Anal.}, 36(3):279--298.

\bibitem[Bertin et~al., 2016]{BLR16}
Bertin, K., Lacour, C., and Rivoirard, V. (2016).
\newblock Adaptive pointwise estimation of conditional density function.
\newblock {\em Ann. Inst. H. Poincar\'e Probab. Statist.}, 52(2):939--980.

\bibitem[Bouaziz and Lopez, 2010]{BouazizLopez10}
Bouaziz, O. and Lopez, O. (2010).
\newblock Conditional density estimation in a censored single-index regression
  model.
\newblock {\em Bernoulli}, 16(2):514--542.

\bibitem[Brunel et~al., 2007]{brunelcomtelacour07}
Brunel, E., Comte, F., and Lacour, C. (2007).
\newblock Adaptive estimation of the conditional density in the presence of
  censoring.
\newblock {\em Sankhy\=a}, 69(4):734--763.

\bibitem[Chagny, 2013]{Chagny}
Chagny, G. (2013).
\newblock Warped bases for conditional density estimation.
\newblock {\em Mathematical Methods of Statistics}, 22(4):253--282.

\bibitem[Comminges and Dalalyan, 2012]{CD12}
Comminges, L. and Dalalyan, A.~S. (2012).
\newblock Tight conditions for consistency of variable selection in the context
  of high dimensionality.
\newblock {\em The Annals of Statistics}, 40(5):2667--2696.

\bibitem[De~Gooijer and Zerom, 2003]{DgZ03}
De~Gooijer, J.~G. and Zerom, D. (2003).
\newblock On conditional density estimation.
\newblock {\em Statist. Neerlandica}, 57(2):159--176.

\bibitem[Delyon et~al., 2016]{DP16}
Delyon, B., Portier, F., et~al. (2016).
\newblock Integral approximation by kernel smoothing.
\newblock {\em Bernoulli}, 22(4):2177--2208.

\bibitem[Efromovich, 2010]{efromovich10b}
Efromovich, S. (2010).
\newblock Dimension reduction and adaptation in conditional density estimation.
\newblock {\em Journal of the American Statistical Association},
  105(490):761--774.

\bibitem[Fan et~al., 1996]{FYT96}
Fan, J., Yao, Q., and Tong, H. (1996).
\newblock Estimation of conditional densities and sensitivity measures in
  nonlinear dynamical systems.
\newblock {\em Biometrika}, 83(1):189--206.

\bibitem[Fan and Yim, 2004]{FYim04}
Fan, J. and Yim, T.~H. (2004).
\newblock A crossvalidation method for estimating conditional densities.
\newblock {\em Biometrika}, 91(4):819--834.

\bibitem[Fan et~al., 2009]{FPYZ09}
Fan, J.-q., Peng, L., Yao, Q.-w., and Zhang, W.-y. (2009).
\newblock Approximating conditional density functions using dimension
  reduction.
\newblock {\em Acta Mathematicae Applicatae Sinica, English Series},
  25(3):445--456.

\bibitem[Hall et~al., 2004]{HRL04}
Hall, P., Racine, J., and Li, Q. (2004).
\newblock Cross-validation and the estimation of conditional probability
  densities.
\newblock {\em J. Amer. Statist. Assoc.}, 99(468):1015--1026.

\bibitem[Holmes et~al., 2010]{HGI10}
Holmes, M.~P., Gray, A.~G., and Isbell, C.~L. (2010).
\newblock Fast kernel conditional density estimation: A dual-tree monte carlo
  approach.
\newblock {\em Computational Statistics \& Data Analysis}, 54(7):1707 -- 1718.

\bibitem[Hyndman et~al., 1996]{HBG96}
Hyndman, R.~J., Bashtannyk, D.~M., and Grunwald, G.~K. (1996).
\newblock Estimating and visualizing conditional densities.
\newblock {\em J. Comput. Graph. Statist.}, 5(4):315--336.

\bibitem[Ichimura and Fukuda, 2010]{IchimuraFukuda10}
Ichimura, T. and Fukuda, D. (2010).
\newblock A fast algorithm for computing least-squares cross-validations for
  nonparametric conditional kernel density functions.
\newblock {\em Computational Statistics \& Data Analysis}, 54(12):3404--3410.

\bibitem[Izbicki and Lee, 2016]{IzbickiLee16}
Izbicki, R. and Lee, A.~B. (2016).
\newblock Nonparametric conditional density estimation in a high-dimensional
  regression setting.
\newblock {\em Journal of Computational and Graphical Statistics},
  25(4):1297--1316.

\bibitem[Izbicki and Lee, 2017]{IzbickiLee17}
Izbicki, R. and Lee, A.~B. (2017).
\newblock Converting high-dimensional regression to high-dimensional
  conditional density estimation.
\newblock {\em Electron. J. Statist.}, 11(2):2800--2831.

\bibitem[Izbicki et~al., 2018]{IzbickiLeePospisil18}
Izbicki, R., Lee, A.~B., and Pospisil, T. (2018).
\newblock Abc-cde: Towards approximate bayesian computation with complex
  high-dimensional data and limited simulations.
\newblock {\em arXiv preprint arXiv:1805.05480}.

\bibitem[Lafferty and Wasserman, 2008]{LW08}
Lafferty, J. and Wasserman, L. (2008).
\newblock Rodeo: Sparse, greedy nonparametric regression.
\newblock {\em Ann. Statist.}, 36(1):28--63.

\bibitem[Le~Pennec and Cohen, 2013]{CohenLePennec}
Le~Pennec, E. and Cohen, S. (2013).
\newblock Partition-based conditional density estimation.
\newblock {\em ESAIM: Probability and Statistics}, eFirst.

\bibitem[Lincheng and Zhijun, 1985]{LinchengZhijun85}
Lincheng, Z. and Zhijun, L. (1985).
\newblock Strong consistency of the kernel estimators of conditional density
  function.
\newblock {\em Acta Mathematica Sinica}, 1(4):314--318.

\bibitem[Liu et~al., 2007]{LLW07}
Liu, H., Lafferty, J.~D., and Wasserman, L.~A. (2007).
\newblock Sparse nonparametric density estimation in high dimensions using the
  rodeo.
\newblock In {\em International Conference on Artificial Intelligence and
  Statistics}, pages 283--290.

\bibitem[Nguyen, 2018]{Jeanne1}
Nguyen, M.-L. (2018).
\newblock Nonparametric method for sparse conditional density estimation in
  moderately large dimensions.
\newblock {\em arXiv:1801.06477}.

\bibitem[Nguyen, 2019]{Jeanne-these}
Nguyen, M.-L. (2019).
\newblock {\em Estimation non param\'etrique de densit\'es conditionnelles :
  grande dimension, parcimonie et algorithmes gloutons}.
\newblock PhD thesis, Universit\'e Paris-Saclay.

\bibitem[Nguyen et~al., 2021]{JCV-supp}
Nguyen, M.-L., Lacour, C., and Rivoirard, V. (2021).
\newblock Supplementary material of adaptive greedy algorithm for moderately
  large dimensions in kernel conditional density estimation.
\newblock {\em Submitted}.

\bibitem[Otneim and Tj{\o}stheim, 2018]{OtneimTjostheim18}
Otneim, H. and Tj{\o}stheim, D. (2018).
\newblock Conditional density estimation using the local gaussian correlation.
\newblock {\em Statistics and Computing}, 28(2):303--321.

\bibitem[Raynal et~al., 2018]{RaynalMPR18}
Raynal, L., Marin, J.-M., Pudlo, P., Ribatet, M., Robert, C.~P., and Estoup, A.
  (2018).
\newblock Abc random forests for bayesian parameter inference.
\newblock {\em Bioinformatics}, 35(10):1720--1728.

\bibitem[Rebelles, 2015]{rebelles15}
Rebelles, G. (2015).
\newblock Pointwise adaptive estimation of a multivariate density under
  independence hypothesis.
\newblock {\em Bernoulli}, 21(4):1984--2023.

\bibitem[Rosenblatt, 1969]{rosenblatt69}
Rosenblatt, M. (1969).
\newblock Conditional probability density and regression estimators.
\newblock In {\em Multivariate {A}nalysis, {II} ({P}roc. {S}econd {I}nternat.
  {S}ympos., {D}ayton, {O}hio, 1968)}, pages 25--31. Academic Press, New York.

\bibitem[Sart, 2017]{sart17}
Sart, M. (2017).
\newblock Estimating the conditional density by histogram type estimators and
  model selection.
\newblock {\em ESAIM: Probability and Statistics}, 21:34--55.

\bibitem[Shiga et~al., 2015]{Shiga2015}
Shiga, M., Tangkaratt, V., and Sugiyama, M. (2015).
\newblock Direct conditional probability density estimation with sparse feature
  selection.
\newblock {\em Machine Learning}, 100(2):161--182.

\bibitem[Tsybakov, 1998]{tsybakov98}
Tsybakov, A.~B. (1998).
\newblock Pointwise and sup-norm sharp adaptive estimation of functions on the
  {S}obolev classes.
\newblock {\em Ann. Statist.}, 26(6):2420--2469.

\bibitem[Wasserman and Lafferty, 2006]{LW06}
Wasserman, L. and Lafferty, J.~D. (2006).
\newblock Rodeo: Sparse nonparametric regression in high dimensions.
\newblock In {\em Advances in Neural Information Processing Systems}, pages
  707--714.

\end{thebibliography}
